\documentclass{amsart}

\usepackage{graphicx}
\usepackage{amsmath}
\usepackage{amscd}
\usepackage{amsthm}
\usepackage{diagrams}
\usepackage{amsfonts}
\usepackage{amssymb}

\newtheorem{theorem}{Theorem}[section]

\newtheorem*{acknowledgements}{Acknowledgements}

\newtheorem{conjecture}[theorem]{Conjecture}
\newtheorem{corollary}[theorem]{Corollary}

\newtheorem{definition}[theorem]{Definition}

\theoremstyle{definition}
\newtheorem{example}[theorem]{Example}
\newtheorem{remark}[theorem]{Remark}
\theoremstyle{plain}

\newtheorem{lemma}[theorem]{Lemma}

\newtheorem{proposition}[theorem]{Proposition}

\numberwithin{equation}{section} 

\begin{document}

\title[Higher-Order Polynomial Invariants of 3-Manifolds]{Higher-Order Polynomial Invariants of 3-Manifolds Giving Lower Bounds for the
Thurston Norm}

\author[Shelly L. Harvey]{Shelly L. Harvey$^\dag$}

\thanks{$^{\dag}$The author was
supported by NSF DMS-9803694 and NSF DMS-0104275 as well as by the
Bob E. and Lore Merten Watt Fellowship}
\date{}

\maketitle

\section{Introduction}

\subsection{Summary of Results}

In this paper we define new 3-manifold invariants and show that
they give new information about the topology of the 3-manifold.
Given a 3-dimensional manifold $X$ and a cohomology class $\psi
\in H^1\left(X;\mathbb{Z}\right)$ we define a sequence of
invariants $\delta_{n}\left(\psi\right)$ which arise as degrees of
``higher-order Alexander polynomials.'' These integers measure the
``size'' of the successive quotients, $G_{r}^{\left( n+1\right)
}/G_{r}^{\left( n+2\right) }$, of the terms of the (rational)
derived series of $G = \pi _{1}\left( X\right) $. Loosely
speaking, $\delta_{n}\left( \psi\right)  $ is the degree of a
polynomial that kills the elements of the first homology of the
regular $G/G_{r}^{\left( n+1\right)  }$-cover of $X$.
The precise definitions are given in Section \ref{defs}. In the
case of knot exteriors and zero surgery on knots, these covering
spaces were studied by T. Cochran, K. Orr, and P. Teichner in
\cite{COT} and \cite{Co}. They defined similar generalized
Alexander modules and were able to obtain important new results on
knot concordance.

Although these invariants are defined algebraically, they have
many exciting topological applications. We show that the degree
$\delta_n$ of each of our family of polynomials gives a lower bound for
the Thurston norm of a class in $H_2\left(X,\partial X
;\mathbb{Z}\right)\cong H^1\left(X;\mathbb{Z}\right)$ of a
3-dimensional manifold. We show that these invariants can give
much more precise estimates of the Thurston norm than previously
known computable invariants. We also show that the $\delta_n$ give
obstructions to a 3-manifold fibering over a circle and to a
3-manifold being Seifert fibered. Moreover, we show that the
$\delta_n$ give computable algebraic obstructions to a 4-manifold
of the form $X \times S^1$ admitting a symplectic structure even
when the obstructions given by the Seiberg-Witten invariants fail.
Some other applications are to the minimal ropelength and genera
of knots and links in $S^3$.

Note that $G_{r}^{\left( n+1\right) }/G_{r}^{\left( n+2\right) }$
is a module over $\mathbb{Z}[G/G_r^{\left(n+1\right)}]$. When
$n=0$, $G_{r}^{\left( 1\right) }/G_{r}^{\left( 2\right) }$ is the
classical Alexander module.  Since $G/G_{r}^{\left( 1\right)}$ is
the (torsion-free) abelianization of G, $G_{r}^{\left( 1\right)
}/G_{r}^{\left( 2\right) }$ is a module over the commutative
polynomial ring in several variables
$\mathbb{Z}[G/G_r^{\left(1\right)}]$. These modules have been
studied thoroughly and with much success. For general $n$,
however, these ``higher-order Alexander modules'' are modules over
non-commutative rings. Very little was previously known in this
case due to the difficulty of classifying such modules.

Let $X$ be a compact, connected, orientable 3-manifold and let
$\psi\in H^{1}\left( X;\mathbb{Z}\right)$.  There is a
Poincar\'{e} duality isomorphism $H^1\left(X;\mathbb{Z}\right)
\cong H_2\left(X,\partial X;\mathbb{Z}\right)$.  If an oriented
surface $F$ in $X$ represents a class $\left[F\right] \in
H_2\left(X,\partial X;\mathbb{Z}\right)$ that corresponds to
$\psi$ under this isomorphism, we say that $F$ is dual to $\psi$
and vice versa.  We measure the complexity of $X$ via the Thurston
norm which is defined in \cite{Th} as follows. If $F$ is any
compact connected surface, let $\chi\left(X\right)$ be its Euler
characteristic and let $\chi\_\left( F\right) =\left| \chi\left(
F\right) \right| $ if $\chi\left( F\right) \leq0$ and equal $0$
otherwise. For a surface $F=\amalg F_{i}$ with multiple
components, let $\chi\_\left( F\right) =\Sigma\chi\_\left(
F_{i}\right)  $. Note that $-\chi\left(F\right) \leq
\chi\_\left(F\right) $ in all cases.  The \emph{Thurston norm} of
$\psi\in H^{1}\left( X;\mathbb{Z}\right)$ is
\[
\left\|  \psi\right\|  _{T}=\inf\left\{  \chi\_\left(  F\right)
|F\text{ is a properly embedded oriented surface dual to
}\psi\right\}  .
\]
This norm extends continuously to all of $H^{1}\left(
X;\mathbb{R}\right)  $. This norm is difficult to compute except
for in the simplest of examples because it is a minimum over an
unknown set.

Thurston showed that the unit ball of the norm is a finite sided
polyhedron and that the set of classes of $H_2\left(X,\partial
X;\mathbb{R}\right)$ representable by a fiber of a fibration over
$S^1$ corresponds to lattice points lying in the cone of the union
of some open faces of this polyhedron \cite{Th}.  This norm has
been useful in the resolution of many open questions in
3-manifolds. D. Gabai used the Thurston norm to show the existence
of taut, finite-depth codimension one foliations 3-manifolds (see
\cite{G1,G2,G3}). In particular, he shows that if $X$ is a
compact, connected, irreducible and oriented 3-manifold and $F$ is
any norm minimizing surface then there is a taut foliation of
finite depth containing $F$ as a compact leaf. Corollaries of
Gabai's existence theorems are that the Property R and Poenaru
conjectures are true.

In a recent paper, C. McMullen defined the Alexander norm of a
cohomology class of a 3-manifold via the Alexander polynomial and
proved that it is a lower bound for the Thurston norm \cite{Mc}.
This theorem has also been recently proved by S. Vidussi \cite{Vi}
using Seiberg-Witten theory  and the work of P. Kronheimer, T.
Mrowka, G. Meng, and C. Taubes \cite{K1,K2,KM,MT}. We prove in
Section \ref{sec:reltothurstonnorm} that the (unrefined)
higher-order degrees $\overline{\delta}_n$ also give lower bounds
for the Thurston norm.  When $n=0$,
$\overline{\delta}_0\left(\psi\right)=\left\|\psi\right\|_A$ hence
this gives another proof of McMullen's theorem.
\newtheorem*{delbartheorem}{Theorem \ref{delbarthm}}
\begin{delbartheorem}
Let $X$ be a compact, orientable 3-manifold (possibly with
boundary). For all $\psi\in H^{1}\left( X;\mathbb{Z}\right)  $ and
$n\geq0$
\[
\overline{\delta}_{n}\left(  \psi\right)  \leq\left\| \psi\right\|
_T
\]
except for the case when $\beta_{1}\left(  X\right) =1$, $n=0$,
$X\ncong S^1 \times S^2$, and $X \ncong S^1 \times D^2$. In this
case,  $\overline{\delta}_{0}\left( \psi\right) \leq\left\|
\psi\right\| _{T}+1+\beta_{3}\left( X\right)  $ whenever $\psi$ is
a generator of $H^{1}\left( X;\mathbb{Z}\right) \cong\mathbb{Z}$.
Moreover, equality holds in all cases when $\psi:\pi_{1}\left(
X\right) \twoheadrightarrow \mathbb{Z}$ can be represented by a
fibration $X\rightarrow S^{1}$.
\end{delbartheorem}

This theorem generalizes the classical result that for a knot
complement, the degree of the Alexander polynomial is less than or
equal to twice the genus of the knot.  It also generalizes
McMullen's theorem.  We remark that $\delta_n=\overline{\delta}_n$
except for some cases where $\overline{\delta}_n=0$.  In fact, for
most of the cases that we are interested in, the
$\overline{\delta}_n$ in Theorem \ref{delbarthm} can be replaced
with $\delta_n$.

Not only do the $\delta_n$ give lower bounds for the Thurston
norm, but we construct 3-manifolds for which $\delta_{n}$ give
much sharper bounds for the Thurston norm than bounds given by the
Alexander norm.  In Theorem \ref{realization}, we start with a
3-manifold $X$ and subtly alter it to obtain a new 3-manifold
$X^{\prime}$. The resulting $X^{\prime}$ cannot be distinguished
from $X$ using the $i^{th}$-order Alexander modules for $i<n$ but
the $n^{th}$-order degrees of $X^{\prime}$ are strictly greater
than those of $X$. We alter a fibered 3-manifold in this manner to
obtain the following result.

\newtheorem*{increasetheorem}{Theorem \ref{increase}}
\begin{increasetheorem}
For each $m\geq1$ and $\mu\geq2$ there exists a 3-manifold $X$
with $\beta_{1}\left(  X\right)  =\mu$ such that
\[
\left\|  \psi\right\|  _{A}=\delta_{0}\left(  \psi\right)
<\delta_{1}\left( \psi\right)  <\cdots<\delta_{m}\left(
\psi\right)  \leq\left\|  \psi\right\| _{T}
\]
for all $\psi\in H^{1}\left(  X;\mathbb{Z}\right)$.  Moreover, $X$
can be chosen so that it is closed, irreducible and has the same
classical Alexander module as a 3-manifold that fibers over
$S^1$.\end{increasetheorem}

An interesting application of Theorem \ref{delbarthm} is to show
that the $\delta_{n}$ give new obstructions to a 3-manifold
fibering over $S^{1}$.  The previously known algebraic
obstructions to a 3-manifold fibering over $S^{1}$ are that the
Alexander module is finitely generated and (when $\beta_{1}\left(
X\right)  =1$) the Alexander polynomial is monic.  For $i,j,n \geq
0$ let $d_{ij}=\delta_i - \delta_j$ and let $r_n\left(X\right)$ be
the $n^{th}$-order rank of the module
$G_r^{\left(n+1\right)}/G_r^{\left(n+2\right)}$ (see Section
\ref{defs}).

\newtheorem*{obsfibtheorem}{Theorem \ref{obsfib}}
\begin{obsfibtheorem}Let $X$ be a compact, orientable 3-manifold.
 If at least one of the following conditions is satisfied then $X$ does not
 fiber over $S^1$.  \newline
$(1)$ $r_n\left(X\right)\neq 0$ for some $n \geq 0$,
\newline
$(2)$ $\beta_{1}\left( X\right) \geq2$ and there exists $i,j\geq
0$ such that $d_{ij}\left(\psi\right) \neq 0$ for all $\psi\in
H^{1}\left( X;\mathbb{Z}\right)$,
\newline
$(3)$ $\beta_{1}\left( X\right) =1$ and $d_{ij}\left(\psi\right)
\neq 0$ for some $i,j\geq 1$ and $\psi \in
H^1\left(X;\mathbb{Z}\right)$,
\newline$(4)$ $\beta_{1}\left( X\right) =1$, $X\ncong S^1 \times S^2$, $X \ncong S^1 \times D^2$ and $d_{0j}\left( \psi\right)\neq 1+\beta_{3}\left( X\right)$ for some
$j\geq 1$ where $\psi$ is a generator of
$H^1\left(X;\mathbb{Z}\right)$.
\end{obsfibtheorem}

As a corollary, we see that the examples in Theorem \ref{increase}
cannot fiber over $S^1$ but have the same classical Alexander
module and polynomial as a fibered 3-manifold.

\newtheorem*{notfibercor}{Corollary \ref{cor:notfiber}}
\begin{notfibercor}For each $\mu\geq 1$, Theorem \ref{increase} gives an infinite family of
closed irreducible 3-manifolds $X$ where
$\beta_1\left(X\right)=\mu$, $X$ does not fiber over $S^1$, and
$X$ cannot be distinguished from a fibered 3-manifold using the
classical Alexander module.
\end{notfibercor}

A second application of Theorem \ref{delbarthm} is to show that
the $\overline{\delta}_{n}$ give obstructions to a 4-manifold of
the form $X\times S^{1}$ admitting a symplectic structure.
Recently, S. Vidussi has extended the work of P. Kronheimer to
show that if a 4-manifold of the form $X\times S^{1}$ ($X$
irreducible) admits a symplectic structure then there is a face of
the Thurston norm ball of $X$ that is contained in a face of the
Alexander norm ball of $X$.  We use his work to prove the
following.

\newtheorem*{sympltheorem}{Theorem \ref{sympl}}
\begin{sympltheorem}
Let $X$ be a closed irreducible 3-manifold.  If at least one of
the following conditions is satisfied then $X \times S^1$ does not
admit a symplectic structure.
\newline$(1)$ $\beta_1\left(X\right) \geq 2$ and there exists an $n \geq 1$ such that
$\overline{\delta}_{n}\left( \psi\right)
>\overline{\delta}_{0}\left(  \psi\right) $ for all $\psi\in H^{1}\left(  X;\mathbb{Z}\right)$.
\newline$(2)$ $\beta_1\left(X\right) =1$,  $\psi$ is a generator of $H^{1}\left(  X;\mathbb{Z}\right)$,  and
$\overline{\delta}_{n}\left( \psi\right)
>\overline{\delta}_{0}\left(  \psi\right) - 2$ for some $n \geq 1$.
\end{sympltheorem}

Hence if $X$ is one of the examples in Theorem \ref{increase},
then $X\times S^{1}$ cannot admit a symplectic structure.  We note
that $X$ has the same Alexander module as a fibered 3-manifold
hence $X \times S^1$ cannot be distinguished from a symplectic
4-manifold using the Seiberg-Witten invariants.

\newtheorem*{notsymplcor}{Corollary \ref{cor:notsympl}}
\begin{notsymplcor}For each $\mu\geq 1$, Theorem \ref{increase} gives an infinite family of
4-manifolds $X \times S^1$ where $\beta_1\left(X\right)=\mu$, $X
\times S^1$ does not admit a symplectic structure, and $X$ cannot
be distinguished from fibered 3-manifold using the classical
Alexander module.
\end{notsymplcor}

Another application of Theorem \ref{delbarthm} is to give
computable lower bounds for the \emph{ropelength} of knots and
links. The ropelength of a link is the quotient of its length by
it's thickness.  In \cite{CKS}, J. Cantarella, R. Kusner, and J.
Sullivan show that the minimal ropelength $R\left(L_i\right)$ of
the $i^{th}$ component of a link $L=\coprod L_i$ is bounded from
below by $2\pi \left(1+\sqrt{\left\| \psi_i \right\| _T}\right)$.
Here $\psi_i$ is the cohomology class that evaluates to 1 on the
meridian of $L_i$ and 0 on the meridian of every other component
of $L$.  In Example \ref{ex:ropelength}, we use Corollary
\ref{cor:rl} to estimate ropelength for a specific link from
\cite{CKS}.

\newtheorem*{rlcor}{Corollary \ref{cor:rl}}
\begin{rlcor}Let $X=S^3-L$ and $\psi_i$
be as defined above.  For each $n \geq 0$,
\[ R\left(L_i\right) \geq 2\pi
\left(1+\sqrt{\delta_n\left(\psi_i\right)-1}\right).
\]Moreover, if $\beta_1\left(X\right) \geq 2$ or $n \geq 1$ (or both) then
\[ R\left(L_i\right) \geq 2\pi
\left(1+\sqrt{\overline{\delta}_n\left(\psi_i\right)}\right).
\]
\end{rlcor}

Lastly, we remark that the higher-order degrees give obstructions
to a 3-manifold admitting a Seifert fibration.

\newtheorem*{SFSprop}{Proposition \ref{pro:SFS}}
\begin{SFSprop}
Let $X$ be a compact, orientable Seifert fibered manifold that
does not fiber over $S^1$.  If $\beta_1\left(X\right) \geq 2$ or
$n \geq 1$ then
\[\overline{\delta}_{n}\left(
\psi\right) =0
\]
 for all $\psi\in H^{1}\left( X;\mathbb{Z}\right)
$.
\end{SFSprop}

\subsection{Outline of Paper}

In Section \ref{sec:alexpoly}, we review the classical Alexander
module, the multivariable Alexander polynomial, and the Alexander
norm of a 3-manifold.  In Section \ref{sec:rds} we define the
rational derived series of a group. This series is a slight
modification of the derived series so that successive quotients
are torsion free. This series will be used to define the
higher-order covers of a 3-manifold, the first homology of which
will be the chief object of study in this paper.

In Section \ref{skewpolyrings} we define certain skew Laurent
polynomial rings $\mathbb{K}_{n}\left[ t^{\pm1}\right]  $ which
contain $\mathbb{Z}\Gamma_{n}$ and depend on a class in the first
cohomology of the 3-manifold. Here, $\Gamma_{n}$ is the group of
deck translations of the higher-order covers. These will be
extremely important in our investigations. Of particular
importance is the fact that they are noncommutative (left and
right) principal ideal domains. Similar rings were used in the
work of \cite{COT}, where it was essential that the rings were
PIDs.

In Section \ref{defs} we define the new higher-order invariants.
If $X$ is any topological space, we define the higher-order
Alexander module and rank of $X$. Finally, if $\psi\in H^{1}\left(
X;\mathbb{Z}\right)  $ we define the higher-order degrees
$\delta_{n}\left(  \psi\right)  $ and $\overline{\delta}_{n}\left(
\psi\right)$.

Section \ref{sec:foxcalc} is devoted to the computation of these
invariants using Fox's Free Calculus. That is, the higher-order
invariants can be computed directly from a finite presentation of
$\pi_{1}\left( X\right)$. The reader familiar with Fox's Free
Calculus should be aware that the classical definitions must be
slightly altered since we are using right instead of left modules.

In Section \ref{sec:presmatrixsection}, we give a finite
presentation of the homology of $X$ with coefficients in
$\mathbb{K}_{n}\left[ t^{\pm1}\right] $. This will be crucial
 to prove that $\overline{\delta}_n$ is bounded above by
the Thurston norm.  In Section \ref{calcexamples}, we compute the
higher-order invariants of some well known 3-manifolds and give
some topological properties of the invariants. The most important
computation in this section is the computation of the higher-order
degrees and ranks for 3-manifolds that fiber over $S^{1}$.

Section \ref{sec:ranktorsionmodules} contains the algebra
concerning the rank of a torsion module over a skew (Laurent)
polynomial ring. Proposition \ref{boundrank} will be used in the
proof of Theorem \ref{delbarthm}. In Section
\ref{sec:reltothurstonnorm}, we show that the higher-degrees are
lower bounds for the Thurston norm. We also prove a theorem
relating higher-order degrees of a cohomology class $\psi$ to the
first Betti number of a surface dual to $\psi$, and prove a result
for links in $S^3$.

In Section \ref{realizationsection} we prove the Realization
Theorem and construct examples of 3-manifolds whose higher degrees
increase. We finish the paper by investigating the applications of
Theorem \ref{delbarthm} to 3-manifolds that fiber over $S^{1}$ and
symplectic 4-manifolds of the form $X\times S^{1}$ in Section
\ref{applications}.

\begin{acknowledgements}I would like to thank Tim Cochran
for many helpful conversations.  I would also like to thank
Cameron Gordon for his advice concerning the irreducibility of my
examples.
\end{acknowledgements}

\section{The Alexander Polynomial\label{sec:alexpoly}}

In this section, we define the Alexander polynomial, the Alexander
module, and the Alexander norm of a 3-manifold. For more
information about the Alexander polynomial we refer the reader to
\cite{CS} , \cite{Hi}, \cite{Ka} and \cite{Ro}.

Let $G$ be a finitely presented group and let $X$ be a finite CW
complex with $\pi_{1}\left( X,x_{0}\right)  $ isomorphic to $G$.
We can assume that $X$ has one 0-cell, $x_{0}$. Let $X_{0}$ be the
universal torsion free abelian cover of $X$ and $\tilde{x}_{0}$ be
the inverse image of $x_{0}$ in $X_{0}$. That is, $X_{0}$ is the
cover induced by the homomorphism from $G$ onto
$\operatorname{ab}\left( G\right) $. Here,
$\operatorname{ab}\left( G\right) =\left( G/\left[ G,G\right]
\right) /\left\{ \mathbb{Z}\!-\!\operatorname{torsion}\right\}  $
which is isomorphic to $\mathbb{Z}^{\mu}$ where
$\mu=\beta_{1}\left( X\right)  $ is the first betti number of $X$.
(The reason for the ``0'' in $X_{0}$ will become apparent later in
the paper.)

The \emph{Alexander Module} of $X$ is defined to be
\[
A_{X}=H_{1}\left(  X_{0},\tilde{x}_{0};\mathbb{Z}\right)
\]
considered as a $\mathbb{Z}\left[  \operatorname{ab}\left(
G\right) \right] $-module. \ After choosing a basis $\left\{
x_{1},\ldots,x_{\mu}\right\}  $ for $H_{1}\left(  X\right)  $ the
ring $\mathbb{Z}\left[  \operatorname{ab}\left(  G\right) \right]
$ can be identified with the ring of Laurent polynomials in
several variables $x_{1},\ldots,x_{\mu}$ with integral
coefficients. The ring $\mathbb{Z} \left[  \operatorname{ab}\left(
G\right) \right] $ has no zero divisors and is in fact a unique
factorization domain. We note that $A_{X}$ is finitely presented
as
\[
\mathbb{Z}\left[  \operatorname{ab}\left(  G\right)  \right]
^{r}\overset{\tilde{\partial }_{2}}{\rightarrow}\mathbb{Z}\left[
\operatorname{ab}\left(  G\right)  \right] ^{s}\rightarrow A_{X}
\]
where $r$ is the number of relations and $s$ is the number of
generators of a presentation of $G$. This presentation is obtained
by lifting each cell of $X$ to $\operatorname{ab}\left(  G\right)
$ cells of the torsion free abelian cover, $X_{0}$.

Let $\Lambda$ be a finitely generated free abelian group. For any
finitely presented $\mathbb{Z}\left[  \Lambda\right]  $-module $M$
with presentation
\[
\mathbb{Z}\left[  \Lambda\right]
^{r}\overset{P}{\rightarrow}\mathbb{Z}\left[ \Lambda\right]
^{s}\rightarrow M
\]
we define the i$^{th}$ elementary ideal $E_{i}\left(  M\right)
\subseteq\mathbb{Z}\left[  F\right]  $ to be the ideal generated
by the $\left(  s-i\right)  \times\left(  s-i\right)  $ minors of
the matrix $P$.  This ideal is independent of the presentation of
$M$. The \emph{Alexander ideal} is $I\left(  X\right) =E_{1}\left(
A_{X}\right)  $, the first elementary ideal of $A_{X}$. The
\emph{Alexander polynomial} $\Delta_{X}$ of $X$ is the greatest
common divisor of the elements of the Alexander ideal.
Equivalently, we could have defined $\Delta_{X}$ to be a generator
of the smallest principal ideal containing $I\left(  X\right)  $.
Note that $\Delta_{X}\in\mathbb{Z}\left[  \operatorname{ab}\left(
G\right) \right] $ and is well-defined up to units in
$\mathbb{Z}\left[ \operatorname{ab}\left( G\right) \right] $. We
point out the necessity that $\mathbb{Z}\left[
\operatorname{ab}\left(  G\right) \right]  $ be a UFD in the
definition of $\Delta_{X}$.

Now let $\psi \in H^1\left(X;\mathbb{Z}\right)$.  Let
$\Delta_{X}=\sum_{i=1}^{m}a_{i}g_{i}$ for $a_{i}
\in\mathbb{Z}\setminus\left\{  0\right\}  $ and $g_i \in
\operatorname{ab} \left( G\right)  $.  The \emph{Alexander norm}
of $\psi\in H^{1}\left( X;\mathbb{R}\right)  $ is defined to be
\[
\|  \psi\|_A=\underset{i,j}{\sup}\,\psi\left( g_{i} -g_{j}\right)
\]
where $\psi$ is a homomorphism from $G$ to $\mathbb{Z}$. In this
paper, we view $\mathbb{Z}$ as the multiplicative group generated
by $t$.  Hence the Alexander norm is equal to the degree of the
one-variable polynomial $\sum_{i=1}^{m}\psi\left( g_{i}\right) $
corresponding to $\psi$.

We note that the Alexander (as well as the Thurston) norm is
actually semi-norms since it can be zero on a non-zero vector of
$H^{1}\left( X;\mathbb{R}\right) $.

\section{Rational Derived Series\label{sec:rds}}

This paper investigates the homology of the covering spaces of a
3-manifold corresponding to the rational derived series of a
group. We begin by defining the rational derived series of $G$ and
proving some of the properties of the quotient $G/G_{r}^{\left(
n+1\right) }$. The most important for our purposes will be that
$G/G_{r}^{\left( n+1\right)  }$ is solvable and its successive
quotients $G_{r}^{\left(  i\right)  }/G_{r}^{\left( i+1\right) }$
are $\mathbb{Z}$-torsion-free and abelian.

\begin{definition}
Let $G_{r}^{\left(  0\right)  }=G$.  For $n\geq1$ define
$G_{r}^{\left( n\right)  }=\left[  G_{r}^{\left(  n-1\right)
},G_{r}^{\left(  n-1\right) }\right]  P_{n-1}$ where
$P_{n-1}=\left\{  g\in G_{r}^{\left(  n-1\right) }\mid
g^{k}\in\left[  G_{r}^{\left(  n-1\right)  },G_{r}^{\left(
n-1\right) }\right]  \text{ for some }k\in Z-\left\{  0\right\}
\right\}  $ to be the \emph{$n^{th}$ term of the rational derived
series} of $G$.
\end{definition}

We denote by $\Gamma_{n}$ the quotient $G/G_{r}^{\left( n+1\right)
}$ and by $\phi_{n}$ the quotient map
$G\twoheadrightarrow\Gamma_{n}$. By the following lemma,
$\Gamma_{n}$ is a group. Note that if $G$ is a finite group then
$G_{r}^{\left(  n\right)  }=G$ hence $\Gamma_{n}=\left\{ 1\right\}
$ for all $n\geq0$. Hence, in this paper we will only be
interested in groups with $\beta_{1}\left(  G\right)  \geq1$.

\begin{lemma}
$G_{r}^{\left(  n\right)  }$ is a normal subgroup of
$G_{r}^{\left(  i\right) }$for $0\leq i\leq n$.
\end{lemma}

\begin{proof}
We show that $\left[  G_{r}^{\left(  n-1\right)  },G_{r}^{\left(
n-1\right) }\right]  $ and $P_{n-1}$ are both normal subgroups of
$G$. Since $G_{r}^{\left(  i+1\right)  }\subseteq G_{r}^{\left(
i\right)  }$ for all $i\geq0$, $\left[  G_{r}^{\left(  n-1\right)
},G_{r}^{\left(  n-1\right) }\right]  \subseteq G_{r}^{\left(
n-1\right)  }$, and $P_{n-1}\subseteq G_{r}^{\left(  n-1\right) }$
it follows that $\left[  G_{r}^{\left( n-1\right) },G_{r}^{\left(
n-1\right)  }\right]  $ and $P_{n-1}$ are normal in $G_{r}^{\left(
i\right)  }$ for $0\leq i\leq n$. Therefore $G_{r}^{\left(
n\right)  }$ is a normal subgroup of $G_{r}^{\left(  i\right)
}$for $0\leq i\leq n$. Let $N$ be a normal subgroup of $G$.   Then
$\left[  N,N\right]  $ is normal in $G $ since $g\left(
\prod\left[  n_{1},n_{2}\right]  \right)  g^{-1} =\prod\left[
gn_{1}g^{-1},gn_{2}g^{-1}\right]  $. Therefore $\left[
G_{r}^{\left(  n-1\right)  },G_{r}^{\left(  n-1\right) }\right]  $
is normal in $G$ by induction on $n$. Now we show that $P_{n-1}$
is a closed under multiplication. Let $p_{1},p_{2}\in P_{n-1}$
then for some $k_{1},k_{2}\neq0$, $p_{1}^{k_{1}
},p_{2}^{k_{2}}\in\left[  G_{r}^{\left(  n-1\right)
},G_{r}^{\left( n-1\right)  }\right]  $. Now for any two elements
$w_{1},w_{2}\in G_{r}^{\left(  n-1\right)  }$, we have
$w_{1}w_{2}=w_{2}w_{1}c$ where $c=\left[
w_{1}^{-1},w_{2}^{-1}\right]  \in\left[  G_{r}^{\left(  n-1\right)
},G_{r}^{\left(  n-1\right)  }\right]  $. Hence $\left(
p_{1}p_{2}\right) ^{k_{1}k_{2}}=\left(  p_{1}^{k_{1}}\right)
^{k_{2}}\left(  p_{2}^{k_{2} }\right)  ^{k_{1}}\prod c_{i}$ where
$c_{i}\in\left[ G_{r}^{\left( n-1\right)  },G_{r}^{\left(
n-1\right)  }\right]  $ so $p_{1}p_{2}\in P_{n-1}$, which shows
that $P_{n-1}$ is a subgroup of $G$. Now if $g\in G$ then $\left(
gp_{1}g^{-1}\right) ^{k_{1}}=gp_{1}^{k_{1}}g^{-1}\in\left[
G_{r}^{\left(  n-1\right) },G_{r}^{\left(  n-1\right)  }\right]  $
since $\left[ G_{r}^{\left(  n-1\right)  },G_{r}^{\left(
n-1\right)  }\right] $ is normal in $G.$ Therefore $P_{n-1}$ is a
normal subgroup of $G$.
\end{proof}

\begin{definition}
A group $\Gamma$ is poly-torsion-free-abelian (PTFA) if it admits
a normal series $\left\{  1\right\}  =G_{0}\vartriangleleft
G_{1}\vartriangleleft \cdots\vartriangleleft G_{n}=\Gamma$ such
that each of the factors $G_{i+1}/G_{i}$ is torsion-free abelian.
(In the group theory literature only a subnormal series is
required.)
\end{definition}

\begin{remark}
\label{PTFA}If $A\vartriangleleft G$ is torsion-free-abelian and
$G/A$ is PTFA then $G$ is PTFA. Any PTFA group is torsion-free and
solvable (the converse is not true). Any subgroup of a PTFA group
is a PTFA group \cite[Lemma 2.4, p.421]{Pa}.
\end{remark}

We show that the successive quotients of the rational derived
series are torsion-free abelian. In fact, the following lemma
implies that if $N$ is a normal subgroup of $G_{r}^{\left(
i\right)  }$ with $G_{r}^{\left(  i\right) }/N$
torsion-free-abelian then $G_{r}^{\left(  i+1\right)  }\subseteq
N$ .
\begin{lemma}
\label{derivedquotient}$G_{r}^{\left(  i\right)  }/G_{r}^{\left(
i+1\right) }$ is isomorphic to $\left(  G_{r}^{\left(  i\right)
}\left/  \left[ G_{r}^{\left(  i\right)  },G_{r}^{\left(  i\right)
}\right]  \right. \right)  \left/  \left\{
\mathbb{Z}\!-\!\operatorname{torsion}\right\} \right.  $ for all
$i\geq0$.
\end{lemma}

\begin{proof}
Since $\left[  G_{r}^{\left(  i\right)  },G_{r}^{\left(  i\right)
}\right] \subseteq G_{r}^{\left(  i+1\right)  }$, we can extend
the natural projection $p:G_{r}^{\left(  i\right)
}\twoheadrightarrow G_{r}^{\left(  i\right) }/G_{r}^{\left(
i+1\right)  }$ to a surjective map $p_{1}:G_{r}^{\left( i\right)
}\left/  \left[  G_{r}^{\left(  i\right)  },G_{r}^{\left( i\right)
}\right]  \right.  \twoheadrightarrow G_{r}^{\left( i\right)  }
/G_{r}^{\left(  i+1\right)  }$. If $\left[  g\right] $ is a
torsion element in $G_{r}^{\left(  i\right)  }\left/  \left[
G_{r}^{\left(  i\right)  } ,G_{r}^{\left(  i\right)  }\right]
\right.  $ then $\left[ g\right] ^{k}=\left[  g^{k}\right]  =1 $
so $g\in P_{i}\subseteq G_{r}^{\left( i+1\right)  }$. Hence we can
extend $p_{1}$ to $p_{2}:\left(  G_{r}^{\left( i\right)  }\left/
\left[ G_{r}^{\left(  i\right)  },G_{r}^{\left(  i\right) }\right]
\right.  \right)  /T\twoheadrightarrow G_{r}^{\left(  i\right)
}/G_{r}^{\left(  i+1\right)  }$ where $T$ is the torsion subgroup
of $G_{r}^{\left(  i\right)  }\left/  \left[  G_{r}^{\left(
i\right)  } ,G_{r}^{\left(  i\right)  }\right]  \right.  $. We
show that $p_{2}$ is injective hence is an isomorphism. Suppose
$p_{2}\left( g_{2}\right)  =1$ then $p\left(  g\right)  =1$ for
any $g$ such that $q_{2}\left(  q_{1}\left( g\right)  \right)
=g_{2}$. Hence $g=fh$ where $f\in\left[  G_{r}^{\left( i\right)
},G_{r}^{\left( i\right)  }\right]  $ and $h^{k}\in\left[
G_{r}^{\left(  i\right) },G_{r}^{\left(  i\right)  }\right]  $ for
some $k\neq0$. Therefore $\left(  g_{2}\right)  ^{k}=\left(
q_{2}\left( q_{1}\left(  fh\right)  \right)  \right)  ^{k}=\left(
q_{2}\left( q_{1}\left(  f\right)  q_{1}\left(  h\right)  \right)
\right) ^{k}=\left( q_{2}\left(  q_{1}\left(  h\right)  \right)
\right) ^{k}=q_{2}\left( q_{1}\left(  h^{k}\right)  \right)
=q_{2}\left( 1\right)  =1$ which implies that $g_{2}=1$.
\end{proof}

If $G=\pi_{1}\left(  X\right)  $ this shows that $G_{r}^{\left(
n\right) }/G_{r}^{\left(  n+1\right)  }\cong H_{1}\left(
X_{\Gamma_{n-1}}\right) /\left\{
\mathbb{Z}\!-\!\operatorname{torsion}\right\}  $ where
$X_{\Gamma_{n-1}}$ is the regular $\Gamma_{n-1}$ cover of $X$.
When $n=0$, note that $G/G_{r}^{\left( 1\right)  }=G_{r}^{\left(
0\right)  }/G_{r}^{\left(  1\right)  }\cong H_{1}\left(  X\right)
\left/  \left\{\mathbb{Z}\!-\!\operatorname{torsion}\right\}
\right. \cong\mathbb{Z}^{\beta_{1}\left( X\right)  }$.

\begin{corollary}
\label{gammaisptfa}$\Gamma_{n}$ is a PTFA group.
\end{corollary}

\begin{proof}
Consider the subnormal series
\[
1=\frac{G_{r}^{\left(  n+1\right)  }}{G_{r}^{\left(  n+1\right)  }
}\vartriangleleft\frac{G_{r}^{\left(  n\right)  }}{G_{r}^{\left(
n+1\right)
}}\vartriangleleft\cdots\vartriangleleft\frac{G_{r}^{\left(
i\right)  } }{G_{r}^{\left(  n+1\right)
}}\vartriangleleft\cdots\vartriangleleft \frac{G_{r}^{\left(
1\right)  }}{G_{r}^{\left(  n+1\right)  }}
\vartriangleleft\frac{G_{r}^{\left(  0\right)  }}{G_{r}^{\left(
n+1\right) }}=\Gamma_{n}.
\]
$G_{r}^{\left(  i\right)  }$ is a normal subgroup of
$G_{r}^{\left(  j\right) }$ for $0\leq j\leq i$ hence
$G_{r}^{\left(  i\right)  }/G_{r}^{\left( n+1\right)  }$ is a
normal subgroup of $G_{r}^{\left(  j\right)  } /G_{r}^{\left(
n+1\right)  }$. From the lemma above, $\left( \frac {G_{r}^{\left(
i\right)  }}{G_{r}^{\left(  n+1\right) }}\right)  \left/ \left(
\frac{G_{r}^{\left(  i+1\right)  }}{G_{r}^{\left(  n+1\right)  }
}\right)  \right.  =G_{r}^{\left(  i\right)  }/G_{r}^{\left(
i+1\right)  }$ is isomorphic to $\left(  G_{r}^{\left(  i\right)
}\left/  \left[ G_{r}^{\left(  i\right)  },G_{r}^{\left(  i\right)
}\right]  \right. \right)  \left/  \left\{
\mathbb{Z}\!-\!\operatorname{torsion}\right\} \right.  $ hence is
torsion free and abelian.
\end{proof}

We next show that if the successive quotients of the derived
series of $G$ are torsion-free then the rational derived series
agrees with the derived series.  In general we only know that $
G^{\left(  i\right)  }\subseteq G_{r}^{\left(  i\right)  } $ for
all $i\geq0$.

\begin{corollary}
If $G^{\left(  i\right)  }/G^{\left(  i+1\right)  }$ is
$\mathbb{Z} $-torsion-free for all $i$ then $G_{r}^{\left(
i\right)  }=G^{\left( i\right)  } $ for all $i.$
\end{corollary}

\begin{proof}
We prove this by induction on $i$. First, we know that
$G_{r}^{\left( 0\right)  }=G^{\left(  0\right)  }=G$. Now assume
that $G_{r}^{\left( i\right)  }=G^{\left(  i\right)  }$, then by
assumption
\begin{align*}
G_{r}^{\left(  i\right)  }\left/  \left[  G_{r}^{\left(  i\right)
} ,G_{r}^{\left(  i\right)  }\right]  \right.
&=G^{\left(  i\right)  }/G^{\left(  i+1\right)  }
\end{align*}
is $\mathbb{Z}$-torsion-free. Hence Lemma \ref{derivedquotient}
gives us $G_{r}^{\left(  i\right)  }/G_{r}^{\left(  i+1\right)
}=G_{r}^{\left( i\right)  }\left/  \left[  G_{r}^{\left(  i\right)
},G_{r}^{\left(  i\right) }\right]  \right.  $ hence
$G_{r}^{\left(  i+1\right)  }=\left[ G_{r}^{\left(  i\right)
},G_{r}^{\left(  i\right)  }\right]
=G^{\left(  i+1\right) }$.
\end{proof}

R. Strebel showed that if $G$ is the fundamental group of a
(classical) knot exterior then the quotients of successive terms
of the derived series are torsion-free abelian \cite{Str}. Hence
for knot exteriors we have $G_{r}^{\left(  i\right)  }=G^{\left(
i\right)  }$.  This is also well known to be true for free groups.
Since any non-compact surface has free fundamental group, this
also holds for all orientable surface groups.

\section{Skew Laurent Polynomial Rings}\label{skewpolyrings}

In this section we define some skew Laurent polynomial rings,
$\mathbb{K}_{n}\left[  t^{\pm1}\right]  $, which are obtained from
$\mathbb{Z}\Gamma_{n}$ by inverting elements of the ring that are
``independent'' of $\psi\in H^{1}\left(  G;\mathbb{Z}\right)  $.
Very similar rings were used in the work of Cochran, Orr and
Teichner \cite[Definition 3.1]{COT}.  Skew polynomial rings with
coefficients in a (skew) field are known to be left and right
principal ideal domains as is discussed herein.

Let $\Gamma$ be a PTFA group as defined in the previous section. A
crucial property of $\mathbb{Z}\Gamma$ is that is has a (skew)
quotient field. Recall that if $R$ is a commutative integral
domain then $R$ embeds in its field of quotients. However, if $R$
is non-commutative domain then this is no longer always possible
(and is certainly not as trivial if it does exist).  We discuss
conditions which guarantee the existence of such a (skew) field.

Let $R$ be a ring and $S$ be a subset of $R$.  $S$ is a
\emph{right divisor set} of $R$ if the following properties hold.
\begin{enumerate}
    \item $0 \notin S$, $1 \in S$.
    \item $S$ is multiplicatively closed.
    \item Given $r \in R$, $s \in S$ there exists $r_1 \in R$,
    $s_1 \in S$ with $rs_1 = sr_1$.
\end{enumerate}
It is known that if $S\subseteq R$ is a right divisor set then the
\emph{right quotient ring} $RS^{-1}$ exists (\cite[p.146]{Pa} or
\cite[p.52]{Ste}).  By $RS^{-1}$ we mean a ring containing $R$
with the property that
\begin{enumerate}
    \item Every element of $S$ has an inverse in $RS^{-1}$.
    \item Every element of $RS^{-1}$ is of the form $rs^{-1}$ with $r \in R$, $s \in S$.
\end{enumerate}

If $R$ has no zero-divisors and $S=R-\left\{ 0\right\}$ is a right
divisor set then $R$ is called an \emph{Ore domain}. If $R$ is an
Ore domain, $RS^{-1}$ is a skew field, called the \emph{classical
right ring of quotients} of $R$ (see \cite{Ste}).  It is observed
in \cite[Proposition 2.5]{COT} that the group ring of a PTFA group
has a right ring of quotients.

\begin{proposition}
[pp. 591-592,611 of \cite{Pa}]If $\Gamma$ is PTFA then
$\mathbb{Q}\Gamma$ is a right (and left) Ore domain; i.e.
$\mathbb{Q}\Gamma$ embeds in its classical right ring of quotients
$\mathcal{K}$, which is a skew field.
\end{proposition}

If $\mathcal{K}$ is the (right) ring of quotients of
$\mathbb{Z}\Gamma$, it is a $\mathcal{K}$-bimodule and a
$\mathbb{Z}\Gamma$-bimodule. Note that
$\mathcal{K}=\mathbb{Z}\Gamma\left(\mathbb{Z}\Gamma-\{0\}\right)^{-1}$
as above.  We list a some properties of $\mathcal{K}$.
\begin{remark}\label{flatremark}
If $R$ is an Ore domain and $S$ is a right divisor set then then
$RS^{-1}$ is flat as a left $R$-module  \cite[Proposistion
II.3.5]{Ste}.  In particular, $\mathcal{K}$ is a flat left
$\mathbb{Z} \Gamma$-module, i.e. $\cdot \otimes_{\mathbb{Z}
\Gamma} \mathcal{K}$ is exact.
\end{remark}

\begin{remark}
\label{fact2}Every module over $\mathcal{K}$ is a free module
\cite[Proposistion I.2.3]{Ste} and such modules have a well
defined rank $\operatorname{rk}_{\mathcal{K}}$ which is additive
on short exact sequences \cite[p. 48]{Coh}.
\end{remark}

If $M$ is a right $R$-module with $R$ an Ore domain then the
\emph{rank of $M$}  is defined as $\operatorname{rank} M =
rk_{\mathcal{K}} \left(M \otimes_R \mathcal{K}\right)$. Combining
Remarks \ref{flatremark} and \ref{fact2} we have the following

\begin{remark}\label{eulerchar}If $\mathcal{C}$ is a non-negative
finite chain complex of finitely generated free right
$\mathbb{Z}\Gamma$-modules then the Euler characteristic
$\chi\left(\mathcal{C}\right)=\sum_{i=0}^{\infty}\left(-1\right)^i
\operatorname{rank} C_i$ is defined and is equal to
$\sum_{i=0}^{\infty}\left(-1\right)^i \operatorname{rank}
H_i\left(\mathcal{C}\right)$.
\end{remark}

The rest of this section will be devoted to the rings
$\mathbb{K}_n^{\psi}\left[t^{\pm 1}\right]$.  Consider the group
$\Gamma_{n}=G/G_{r}^{\left( n+1\right) }$ for $n\geq 0$. Since
$\Gamma_{n}$ is PTFA (Corollary \ref{gammaisptfa}),
$\mathbb{Z}\Gamma_{n}$ embeds in its right ring of quotients,
which we denote by $\mathcal{K}_{n}$. Let $\psi\in H^{1}\left(
G;\mathbb{Z}\right)$ be primitive.  Since $H^{1}\left(
G;\mathbb{Z}\right) \simeq \operatorname{Hom}_{\mathbb{Z}}\left(
G,\mathbb{Z}\right)$, $\psi$ can be considered as an epimorphism
from $G$ to $\mathbb{Z}$. In particular, $\psi$ is trivial on
$G_{r}^{\left( n+1\right)  }$ so it induces a well defined
homomorphism $\overline{\psi}:\Gamma_{n}\twoheadrightarrow
 \mathbb{Z}$. Let $\Gamma_{n}^{\prime}$ be
the kernel of $\overline{\psi}$. Since $\Gamma_{n}^{\prime}$ is a
subgroup of $\Gamma_{n}$, $\Gamma_{n}^{\prime}$ is PTFA by Remark
\ref{PTFA}. Therefore $\mathbb{Z}\Gamma_{n}^{\prime}$ is an Ore
domain and $S_{n}=\mathbb{Z}\Gamma_{n}^{\prime}-\left\{ 0\right\}
$ is a right divisor set of $\mathbb{Z}\Gamma_{n}$
\cite[p.609]{Pa}.  Let $\mathbb{K}_{n}=\left(
\mathbb{Z}\Gamma_{n}^{\prime}\right) S_{n}^{-1}$ be the right ring
of quotients of $\mathbb{Z}\Gamma_{n}^{\prime}$,
$g_{\psi}:\mathbb{Z}\Gamma_{n}^{\prime}\rightarrow \mathbb{K}_{n}$
be the embedding of $\mathbb{Z}\Gamma_{n}^{\prime} $ into
$\mathbb{K}_{n}$, and $ R_{n}=\left( \mathbb{Z}\Gamma _{n}\right)
S_{n}^{-1}$. We will show that $ R_{n}$ is isomorphic to a certain
skew Laurent polynomial ring $\mathbb{K}_{n}\left[ t^{\pm1}\right]
$ (defined below) which is a non-commutative principal right and
left ideal domain by \cite[2.1.1. p.49]{Coh2}.  That is,
$\mathbb{K}_n\left[t^{\pm 1}\right]$ has no zero divisors and
every right and left ideal is principal.

We recall the definition of a skew Laurent polynomial ring.  If
$K$ is a skew field, $\alpha$ is an automorphism of $K$ and $t$ is
an indeterminate, the \emph{skew (Laurent) polynomial ring in} $t$
\emph{over }$K$ associated with $\alpha$ is the ring consisting of
all expressions of the form
\[
t^{-m}a_{-m}+\cdots+t^{-1}a_{-1}+a_{0}+ta_{1}+\cdots+t^{n}a_{n}
\]
where $a_{i}\in K$. The operations are coordinate-wise addition
and multiplication defined by the usual multiplication for
polynomials and the rule $at=t\alpha\left(  a\right)  $ \cite[p.
54]{Coh}.

Consider the split short exact sequence
\[ 0\longrightarrow\ker\left(
\overline{\psi}\right)
\longrightarrow\Gamma_{n}\overset{\overset{\_}{\phi}}
{\longrightarrow} \mathbb{Z} \longrightarrow 0. \] Choose a
splitting $\xi:\mathbb{Z}\rightarrow\Gamma_{n}$. Then $\xi$
induces an automorphism of $\Gamma_n^{\prime} = \ker\left(
\overline{\psi} \right) $ by $g\mapsto\xi\left(  t\right)^{-1}
g\xi\left( t\right)$.  This induces a ring automorphism of
$\mathbb{Z}\Gamma_n^{\prime}$ and hence a field automorphism
$\alpha$ of $\mathbb{K}_n$ by $\alpha\left(rs^{-1}\right) =
\alpha\left(r\right)\alpha\left(s\right)^{-1}$.  This defines
$\mathbb{K}_n \left[t^{\pm1}\right]$ as above.

\begin{proposition}
The embedding $g_{\psi}:\mathbb{Z}\Gamma_{n}^{\prime}\rightarrow
\mathbb{K}_{n}$ extends to an isomorphism $ R_{n}\rightarrow
\mathbb{K}_{n}\left[  t^{\pm1}\right] $
\end{proposition}

\begin{proof}
Any element of $\Gamma_n$ can be written uniquely as $\xi\left(
t\right)^m a_m$ for some $m \in \mathbb{Z}$ and $a_m \in
\mathbb{Z}\Gamma_{n}^{\prime}$.  It follows that
$\mathbb{Z}\Gamma_n$ is isomorphic to the skew (Laurent)
polynomial ring $\mathbb{Z}\Gamma_{n}^{\prime}\left[
x^{\pm1}\right]$ by sending $\xi\left(  t\right)^m a_m$ to $x^m
a_m$.  The automorphism of $\mathbb{Z}\Gamma_{n}^{\prime}$ is
induced by conjugation, $a \rightarrow x^{-1} a x$ since $a
\xi\left(  t\right) = \xi\left(  t\right) \left( \xi\left(
t\right)^{-1} a \xi\left(  t\right)\right)$.

Hence there is an obvious ring homomorphism of
$\mathbb{Z}\Gamma_{n}\rightarrow\mathbb{K}_{n}\left[
t^{\pm1}\right]  $ extending $g_{\psi}$. Note that the
automorphism $g\mapsto\xi\left(  t\right)^{-1}  g\xi\left(
t\right)$ defining $\mathbb{K}_{n}\left[  t^{\pm1}\right]  $
agrees with conjugation in $\Gamma$ so this map is a ring
homomorphism. The non-zero elements of $\mathbb{Z}
\Gamma_{n}^{\prime}$ map to invertible elements in
$\mathbb{K}_{n}\left[ t^{\pm1}\right]  $. It is then easy to show
that $ R_{n}\cong \mathbb{K}_{n}\left[ t^{\pm1}\right] $.
\end{proof}

We note that $\left(  \mathbb{Z}\Gamma_{n}\right)  S^{-1}$ depends
on the (primitive) class $\psi\in H^{1}\left( G;\mathbb{Z}\right)
$. Moreover, the isomorphism of $\left(
\mathbb{Z}\Gamma_{n}\right) S^{-1}$ with $\mathbb{K}_{n}\left[
t^{\pm1}\right]  $depends on the splitting
$\xi:\mathbb{Z}\rightarrow\Gamma_{n}$. For any $\psi\in
H^{1}\left( X;\mathbb{Z}\right)  $ we have
\[
\mathbb{Z}\Gamma_{n}\hookrightarrow\mathbb{K}_{n}\left[
t^{\pm1}\right] \hookrightarrow\mathcal{K}_{n}.
\]
One should note that the first and last rings only depend on the
group $G$ while the middle ring $\mathbb{K}_{n}\left[
t^{\pm1}\right]  $ depends on the homomorphism
$\psi:G\rightarrow\mathbb{Z}$ and splitting $\xi:\mathbb{Z}
\rightarrow\Gamma_{n}$. Often we write
$\mathbb{K}_{n}^{\psi}\left[ t^{\pm1}\right]  $ to emphasize the
class $\psi$ on which $\mathbb{K}_{n} ^{\psi}\left[
t^{\pm1}\right]  $ is dependent.  From Remark \ref{flatremark} we
have the following.

\begin{remark}
$\mathbb{K}_{n}^{\psi}\left[  t^{\pm1}\right]  $ and
$\mathcal{K}_{n}$ are flat left $\mathbb{Z}\Gamma_{n}$-modules.
\end{remark}

\section{Definition of Invariants\label{defs}}

Suppose $X$ is a connected CW-complex with $A\subseteq X$ and $x_0
\in A$ a basepoint.  Let $\phi:\pi_{1}\left( X\right)
\rightarrow\Gamma$ be a homomorphism and $X_{\Gamma}
\overset{p}{\rightarrow} X$ denote the regular $\Gamma$-cover of
$X$ associated to $\phi$. That is, $X_\Gamma$ is defined to be the
pullback of the universal cover of $K\left(\Gamma,1\right)$. We
note that there is an induced coefficient system on $A$,
$\phi\circ i_{\ast} :\pi_{1}\left( A\right) \rightarrow\Gamma$
where $i$ is the inclusion map of $A$ into $X$. Thus, we have a
regular covering map of the pair $\left(X_{\Gamma},A_\Gamma\right)
\overset{\operatorname{p}}{\rightarrow} \left(X,A\right)$. If
$\phi$ is not surjective then $X_{\Gamma} $ is the disjoint union
of $\Gamma/\operatorname{Im}\left( \phi\right) $ copies of the
regular $\operatorname{Im}\left( \phi\right) $-cover corresponding
to $\phi:\pi_{1}\left(  X\right)
\twoheadrightarrow\operatorname{Im}\left( \phi\right)$.

There is a natural homomorphism $\tau : \Gamma \rightarrow
G\left(X_\Gamma\right)$ where $G\left(X_\Gamma\right)$ is the
group of deck transformations of $X_\Gamma$ (see \cite{Ha1} or
\cite{Ma}) for more details). We note that $\tau$ is an
isomorphism when $\phi$ is surjective.  This map is defined by
sending $\gamma$ to the deck transformation $\tau_\gamma$ that
takes $x_0$ to $\tilde{\gamma}\left(1\right)$ where
$\tilde{\gamma}$ is the unique lift of $\gamma$ starting at $x_0$.
$\tau$ gives us a left $\Gamma$ action on $X_\Gamma$ by
$\tilde{x}\gamma = \tau_\gamma\left(\tilde{x}\right)$. We make
this into a right action by defining $\tilde{x}\gamma=
\gamma^{-1}\tilde{x}$. Hence $\tilde{x}\gamma=
\tau_{\gamma^{-1}}\left(\tilde{x}\right)$.

The right action of $\Gamma$ on $X_\Gamma$ induces a right action
on the group $C_{\star }\left(  X_{\Gamma}\right)$ of singular
$n$-chains on $X_\Gamma$, by sending a singular $n$-simplex
$\sigma:\Delta^n \rightarrow X_\Gamma$ to the composition
$\Delta^n \rightarrow X_\Gamma \overset{\gamma}{\rightarrow}
X_\Gamma$.  The action of $\Gamma$ on $C_{\star }\left(
X_{\Gamma}\right)$ makes $C_{\star }\left(  X_{\Gamma}\right)$ a
right $\mathbb{Z}\Gamma$-module.

Let $\mathcal{M}$ be a $\mathbb{Z}\Gamma$-bimodule$.$ The
equivariant homology of $X$ and $\left(  X,A\right)  $ are defined
as follows.

\begin{definition}
Given $X$, $A$, $\phi$, $\mathcal{M}$ as above, let
\[
H_{\star}\left(  X;\mathcal{M}\right)  \equiv H_{\star}\left(  C_{\star
}\left(  X_{\Gamma}\right)  \otimes_{\mathbb{Z}\Gamma}\mathcal{M}\right)
\]
and
\[
H_{\star}\left(  X,A;\mathcal{M}\right)  \equiv H_{\star}\left(
C_{\star }\left(  X_{\Gamma},A_{\Gamma}\right)
\otimes_{\mathbb{Z}\Gamma} \mathcal{M}\right)
\]
as right $\mathbb{Z}\Gamma$-modules
\end{definition}
 These are well-known to be isomorphic to the homology of $X$ and
$\left( X,A\right)  $ with coefficient system induced by $\phi$
\cite[Theorem VI 3.4]{W}.

We now restrict to the case when $\Gamma$ is PTFA.  We state the
following useful proposition.  A proof of this can be found in
\cite[Section 3]{Co}.  We remark that the finiteness condition in
Proposition \ref{upperboundrank} is necessary.

\begin{proposition}\label{upperboundrank}Suppose
$\pi_1\left(X\right)$ is finitely generated and $\phi :
\pi_1\left(X\right) \rightarrow \Gamma$ is non-trivial.  Then
\[\operatorname{rank} H_1\left(X;\mathbb{Z}\Gamma\right)\leq \beta_1\left(X\right)-1.\]
\end{proposition}

Let $G = \pi_{1}\left( X,x_{0}\right)$.  Define the $n^{th}$-order
cover $X_{n} \overset{\operatorname{p_n}}{\longrightarrow} X$ of
$X$ to be the regular $\Gamma_{n}$-cover corresponding to the
coefficient system $\phi_{n}:G\twoheadrightarrow \Gamma_{n}$ where
$\Gamma_{n}=G/G_{r}^{\left(  n+1\right)  }$ is as defined in
Section \ref{sec:rds}. Recall that $\mathbb{Z}\Gamma_{n}$ has a
(skew) field of quotients $\mathcal{K}_{n}$. If $ R$ is any ring
with $\mathbb{Z}\Gamma_{n}\subseteq R\subseteq\mathcal{K}_{n}$
then $ R$ is a $\mathbb{Z}\Gamma_{n}$-bimodule. Moreover, $H_{\ast
}\left(  X; R\right)  $ can be considered as a right $ R $-module.
We will be interested in the cases when $ R$ is
$\mathbb{Z}\Gamma_{n}$, $\mathcal{K}_{n}$ and
$\mathbb{K}_{n}\left[  t^{\pm 1}\right]  $ where
$\mathbb{K}_{n}\left[  t^{\pm1}\right]  $ is as described in the
Section \ref{skewpolyrings}. If $M$ is a right (left) $ R$-module
where $ R$ is an Ore domain then we let $T_{ R}M$ be the $
R$-torsion submodule of $M$. When there is no confusion we
suppress the $ R$ and just write $TM$.

We define the higher-order modules.  The integral invariants that
we can extract from these modules will be our chief interest for
the rest of this paper.
\begin{definition}
The $n^{th}$-order Alexander module of a CW-complex $X$ is
\[
\mathcal{A}_{n}\left(  X\right)  \equiv T_{\mathbb{Z}\Gamma_{n}}
H_{1}\left( X;\mathbb{Z}\Gamma_{n}\right)
\]
considered as right $\mathbb{Z}\Gamma_{n}$-module. Similarly, we
define
\[ \overline{\mathcal{A}}_{n}\left(  X\right)  \equiv H_{1}\left(
X;\mathbb{Z}\Gamma_{n}\right)
\]considered as right $\mathbb{Z}\Gamma_{n}$-module.
\end{definition}

Let $X$ and $Y$ be 3-manifolds with $G=\pi_1\left(X\right)$ and
$H=\pi_1\left(Y\right)$.  Suppose that $G$ is isomorphic to $H$.
We would like for their higher-order Alexander modules to be ``the
same.'' However, they are modules over different (albeit
isomorphic) rings.  We remedy this dilemma with the following
definition.  It is easy to verify that the following defines an
equivalence relation.

\begin{definition}Let $M$ and $N$ be right (left) $R$ and $S$-modules
respectively and $f:R \rightarrow S$ be an isomorphism.  $N$ can
be made into a right (left) $R$-module via $f$.  We say that $M$
is equivalent to $N$ provided $N$ is isomorphic to $M$ as a right
(left) $R$-module.
\end{definition}

Let $X$ be a topological space with $G=\pi_1\left(X\right)$.  The
higher-order Alexander modules can be defined group theoretically.
Define a right
$\mathbb{Z}\left[G/G_r^{\left(n+1\right)}\right]$-module structure
on $G_r^{\left(n+1\right)} \left/
\left[G_r^{\left(n+1\right)},G_r^{\left(n+1\right)}\right]
\right.$ by $\left[h\right]\left[g\right] = \left[g^{-1}hg\right]$
for $h \in G_r^{\left(n+1\right)}$ and $g \in G$.  We see that
\[ \overline{\mathcal{A}}_{n}\left( X\right)
\cong
\frac{G_r^{\left(n+1\right)}}{\left[G_r^{\left(n+1\right)},G_r^{\left(n+1\right)}
\right]}
\]
as a right
$\mathbb{Z}\left[\frac{G}{G_r^{\left(n+1\right)}}\right]$-module.
 We also
note that \[\overline{\mathcal{A}}_{n}\left(
X\right)\left/\{\mathbb{Z}\!-\!\operatorname{torsion}\}\right.
\cong G_r^{\left(n+1\right)}/G_r^{\left(n+2\right)}\] by
Proposition \ref{derivedquotient}. Suppose that $Y$ is
homeomorphic to $X$, then $\pi_1\left(Y\right)$ is isomorphic to
$G$.  It is easy to verify that the isomorphism of groups leads to
an equivalence of $\overline{\mathcal{A}}_{n}\left( X\right)$ and
$\overline{\mathcal{A}}_{n}\left( Y\right)$. Therefore the
equivalence classes of the higher-order Alexander modules are
topological invariants. Similarly, one can easily verify that the
rest of the definitions in this section are invariants of $X$ or a
pair $\left(X,\psi\right)$.

\begin{definition}
The $n^{th}$-order rank of $X$ is
\[
r_{n}\left(  X\right)
=\operatorname{rk}_{\mathcal{K}_{n}}H_{1}\left( X;\mathcal{K}
_{n}\right)  .
\]
\end{definition}

In the literature, the classical Alexander module of a 3-manifold
is often defined as $H_{1}\left(
X,x_{0};\mathbb{Z}\Gamma_{0}\right)$ (see Section
\ref{sec:alexpoly}) and $\alpha\left( X\right) =\operatorname{rk}
H_1\left(X,x_0;\mathbb{Z}\Gamma_0\right)$ is called the nullity of
$X$ \cite{Hi}. We will now show that $H_{1}\left(
X;\mathbb{Z}\Gamma_{n}\right) $ and $H_{1}\left(
X,x_{0};\mathbb{Z}\Gamma_{n}\right)  $ are related by $r_{n}\left(
X\right) =\operatorname{rk}_{\mathcal{K}_{n}}H_{1}\left(  X,x_{0}
;\mathcal{K}_{n}\right) -1$ and
$T_{\mathbb{Z}\Gamma_{n}}H_{1}\left( X;\mathbb{Z}\Gamma_{n}\right)
=T_{\mathbb{Z}\Gamma_{n}}H_{1}\left(
X,x_{0};\mathbb{Z}\Gamma_{n}\right)$. Hence, we could have defined
$\mathcal{A} _{n}\left( X\right) $ and $r_{n}\left( X\right)  $
using homology rel basepoint as well.

\begin{proposition}
\label{comp}Let $\Gamma$ be PTFA, and $\phi:\pi_{1}\left(
X,x_{0}\right) \rightarrow\Gamma$ be non-trivial. Then
\begin{equation}
\operatorname{rk}_{\mathcal{K}}H_{1}\left(  X;\mathcal{K}\right)
=\operatorname{rk}_{\mathcal{K}}H_{1}\left(
X,x_{0};\mathcal{K}\right) -1\label{comprank}
\end{equation}
and
\begin{equation}
T_{ R}H_{1}\left(  X;R\right)  \cong T_{ R}H_{1}\left( X,x_{0};
R\right) \label{comptor}
\end{equation}
for any ring $ R$ such that $\mathbb{Z}\Gamma\subseteq R
\subseteq\mathcal{K}$, where $\mathcal{K}$ is the (skew) field of
quotients of $\mathbb{Z}\Gamma$.
\end{proposition}

\begin{proof}
To prove $\left(  \ref{comprank}\right)$, consider the long exact
sequence the pair $\left(X,x_0\right)$,
\[
0\rightarrow H_{1}\left(  X;\mathcal{K}\right)
\overset{\rho}{\rightarrow }H_{1}\left( X,x_{0};\mathcal{K}\right)
\overset{\partial}{\rightarrow} H_{0}\left(
x_{0};\mathcal{K}\right)  \rightarrow H_{0}\left(  X;\mathcal{K}
\right).
\]
Since $\phi:\pi_{1}\left(  X\right)  \rightarrow\Gamma$ is
non-trivial, $H_{0}\left(  X;\mathcal{K}\right)  =0$ by the
following Lemma. The first result follows since $H_{0}\left(
x_{0};\mathcal{K}\right) \cong\mathcal{K}$.

\begin{lemma}\label{lem:H0}Suppose $X$ is a connected CW complex.  If
$\phi: \pi_1\left(X\right) \rightarrow \Gamma$ is a non-trivial
coefficient system and $\Gamma$ is PTFA then
$H_0\left(X;\mathcal{K}\right)=0$.
\end{lemma}
\begin{proof}By \cite[p. 275]{W} and \cite[p.34]{B},
$H_0\left(X;\mathcal{K}\right)$ is isomorphic to the cofixed set
$\mathcal{K}/\mathcal{K}I$ where $I$ is the augmentation ideal of
$\mathbb{Z}\left[\pi_1\left(X\right)\right]$ acting via
$\mathbb{Z}\left[\pi_1\left(X\right)\right] \rightarrow
\mathbb{Z}\Gamma \rightarrow \mathcal{K}$.  If $\phi$ is
non-trivial then the composition is non-trivial.  Thus $I$
contains an element that is a unit hence
$\mathcal{K}I=\mathcal{K}$.
\end{proof}

We show that the map $\rho:H_{1}\left(  X; R\right) \rightarrow
H_{1}\left(  X,x_{0}; R\right)  $ restricts to an isomorphism from
$TH_{1}\left(  X; R\right)$ onto $TH_{1}\left(  X,x_{0} ; R\right)
$. Certainly $\rho:TH_{1}\left( X; R\right) \rightarrow
TH_{1}\left( X,x_{0}; R\right)  $ is a monomorphism. Let
$\sigma\in TH_{1}\left( X,x_{0}; R\right)  $ with $\sigma r=0$
where $r\neq0$. Since $H_{0}\left( x_{0}; R\right) \cong R$ is $
R$-torsion-free, $\partial\left( \sigma\right)  =0$ so there
exists $\theta\in H_{1}\left( X; R\right)  $ with $\rho\left(
\theta\right)  =\sigma$. We see that $\theta$ is $ R$-torsion
since $\rho\left( \theta r\right) =\rho\left(  \theta\right)
r=\sigma r=0$ and $\rho$ is a monomorphism. Therefore
$\rho:TH_{1}\left( X; R \right) \rightarrow TH_{1}\left( X,x_{0};
R\right)  $ is surjective.
\end{proof}

For any primitive class $\psi\in H^{1}\left( X;\mathbb{Z}\right) $
and splitting $\xi:\mathbb{Z}\rightarrow\Gamma_{n}$ we consider
the skew Laurent polynomial ring $\mathbb{K}_{n}\left[
t^{\pm1}\right]$.  We note that $TH_{1}\left( X;\mathbb{K}
_{n}\left[ t^{\pm1}\right]  \right)  $ is a finitely generated
right $\mathbb{K}_{n}$-module. Moreover, any module over
$\mathbb{K}_{n}$ has a well defined rank which is additive on
short exact sequences by Remark \ref{fact2}.

\begin{definition}
Let $X$ be a finite CW-complex. For each primitive $\psi\in
H^{1}\left( X;\mathbb{Z}\right)  $ and $n\geq0$ we define the
\emph{refined $n^{th}$-order Alexander module corresponding to
$\psi$} to be $\mathcal{A}_{n}^{\psi}\left(  X\right)
=T_{\mathbb{K}_{n}\left[ t^{\pm1}\right]  }H_{1}\left(
X;\mathbb{K}_{n}\left[  t^{\pm1}\right] \right)$ viewed as a right
$\mathbb{K}_n\left[t^{\pm 1}\right]$-module.
\end{definition}
Since $\mathcal{A}_{n}^{\psi}\left(  X\right)$ is a finitely
generated module over the principal ideal domain
$\mathbb{K}_n\left[t^{\pm 1}\right]$,
$$\mathcal{A}_{n}^{\psi}\left(  X\right)\cong \bigoplus_{i=1}^{m}
\frac{\mathbb{K}_n\left[t^{\pm 1}\right]}{p_i\left(t\right)
\mathbb{K}_n\left[t^{\pm 1}\right]}$$ for some nonzero
$p_i\left(t\right) \in \mathbb{K}_n\left[t^{\pm 1}\right]$
\cite[Theorem 16, p. 43]{Ja}.  We define the refined
$n^{th}$-order degree of $\psi$ to be the degree of the polynomial
$\prod p_i\left(t\right)$.  One can verify that this is equal to
the rank of $\mathcal{A}_{n}^{\psi}\left( X\right)$ as a
$\mathbb{K}_n$-module.
\begin{definition}Let $X$ be a finite CW-complex. For each primitive $\psi\in
H^{1}\left( X;\mathbb{Z}\right)  $ and $n\geq0$ we define the
\emph{refined $n^{th}$-order degree of $\psi$} to be
\[
\delta_{n}\left(  \psi\right)
=\operatorname{rk}_{\mathbb{K}_{n}}\mathcal{A}_{n}^{\psi }\left(
X\right)
\]
We extend by linearity to define $\delta_{n}\left( \psi\right)$
for non-primitive classes $\psi$.
\end{definition}

Similarly we define the unrefined higher-order Alexander modules
and degrees.

\begin{definition}
Let $X$ is a finite CW-complex.  For each primitive $\psi\in
H^{1}\left( X;\mathbb{Z}\right)  $ and $n\geq0$ we define the
\emph{unrefined $n^{th}$-order Alexander module corresponding to
$\psi$} to be $\overline{\mathcal{A}}_{n}^{\psi}\left(  X\right)
=H_{1}\left( X;\mathbb{K}_{n}\left[  t^{\pm1}\right] \right)$
viewed as a right $\mathbb{K}_n\left[t^{\pm 1}\right]$-module. The
\emph{unrefined $n^{th}$-order degree of $\psi$} is
\[
\overline{\delta}_{n}\left(  \psi\right)
=\operatorname{rk}_{\mathbb{K}_{n}}\overline{\mathcal{A}}_{n}^{\psi
}\left( X\right)
\]if $\operatorname{rk}_{\mathbb{K}_{n}}\overline{\mathcal{A}}_{n}^{\psi
}\left( X\right)$ is finite and 0 otherwise.  We extend by
linearity to define $\overline{\delta}_{n}\left( \psi\right)$ for
non-primitive classes $\psi$.
\end{definition}

We note that $$\overline{\mathcal{A}}_{n}^{\psi }\left(
X\right)\cong \left(\bigoplus_{i=1}^{m}
\frac{\mathbb{K}_n\left[t^{\pm 1}\right]}{p_i\left(t\right)
\mathbb{K}_n\left[t^{\pm 1}\right]} \right) \bigoplus
\mathbb{K}_n\left[t^{\pm 1}\right]^{r_n\left(X\right)}.$$ Hence
$\operatorname{rk}_{\mathbb{K}_n}\overline{\mathcal{A}}_{n}^{\psi
}\left( X\right)$ is finite if and only if $r_{n}\left(  X\right)
=0$.
\begin{remark}\label{rem:del_eq_delbar} If
$r_{n}\left( X\right) =0$ then $\overline{\delta}_{n}\left(
\psi\right) =\delta_{n}\left( \psi\right)  $ otherwise
$0=\overline{\delta}_{n}\left( \psi\right)  \leq\delta_{n}\left(
\psi\right)  $. \end{remark} We now show that
$\overline{\delta}_{0}\left( \psi\right)  $ is equal to the
Alexander norm of $\psi$ hence $\overline{\delta}_{0}\left(
\psi\right) $ is a convex function.

\begin{proposition}\label{pro:del0}
$\overline{\delta}_{0}\left(  \psi\right)  =\left\|  \psi\right\|
_{A}$ for all $\psi\in H^{1}\left(  X;\mathbb{Z}\right)  $.
\end{proposition}
\begin{proof}
Recall that $\Gamma_0=\mathbb{Z}^{\beta_1\left(X\right)}$ hence
$\mathbb{Z}\Gamma_0$ is isomorphic to the polynomial ring in
several variables.  Let $\nu:\mathbb{Z} \Gamma_{0}
\hookrightarrow\mathbb{K} _{0}\left[ t^{\pm1}\right] $ be the
embedding of $\mathbb{Z} \Gamma_{0}  $ into the principal ideal
domain $\mathbb{K}_{0}\left[ t^{\pm1}\right]  $ and
$\Delta_{X}=\sum n_{g}g$ be the Alexander polynomial of $X$. We
begin by showing that $\left\|  \psi\right\| _{A}=\deg\nu\left(
\Delta_{X}\right) $. For all $j$ consider the polynomial
$\Delta_{X} ^{j}=\sum_{\psi\left(  g\right) =t^{j}}n_{g}g$. Note
that any such $g$ can be written (using the splitting $\xi$) as
$h_g \tau^j$ where $\psi\left(\tau\right)=t$. We see that
$\nu\left( \Delta_{X}^{j}\right) =\nu\left( \sum_{\psi\left(
g\right) =t^{j}}n_{g}g\right)  =\left( \sum n_{g}h_{g}\right)
t^{j}$ where $c_{j}\equiv\sum n_{g}h_{g}$ $\in\mathbb{Z}\left[
\ker\overline{\psi} _{0}\right]  $. Since $\nu$ is a monomorphism
we have $c_{j}\neq0$ unless $n_{g}=0$ for all $g$ with $\psi\left(
g\right)  =t^{j}$. \ It follows that $\deg\nu\left(
\Delta_{X}\right) =\deg\nu\left( \sum\Delta_{X}^{j}\right)
=\deg\sum c_{j}t^{j}=\left\| \psi\right\|  _{A}$. After choosing a
group presentation for $G$, Fox's Free Calculus (Section
\ref{sec:foxcalc}) gives us a presentation matrix $M$ for
$H_{1}\left( X,x_{0};\mathbb{Z}\left[ \Gamma _{0}\right]  \right)
=H_{1}\left( \widetilde{X},\widetilde{x}_{0}\right)  $ where
$\widetilde{X}$ is the torsion-free abelian cover of $X$. Moreover
a presentation of $H_{1}\left( X,x_{0};\mathbb{K}_{0}\left[
t^{\pm1}\right] \right) $ is also given by $M^{\nu}$, that is we
consider each entry in as an element of $\mathbb{K}_{0}\left[
t^{\pm1}\right]  $. If $s$ is the number of generators in the
presentation of $G$ then $\Delta_{X}=\gcd\left( E_{1}\left(
H_{1}\left( X,x_{0};\mathbb{Z}\Gamma_{0} \right)  \right) \right)
=\gcd\left\{ d_{1},\ldots,d_{k}\right\}  $ where $\left\{
d_{1},\ldots,d_{k}\right\}  $ is the set of determinants of the
$\left( s-1\right) \times\left(  s-1\right)  $ minors of $M$
(Section \ref{sec:alexpoly}). Note that $\Gamma_0$ is free abelian
so $\mathbb{K}_0$ is a commutative field and hence and
$\mathbb{Z}\Gamma_0$ and $\mathbb{K}_0\left[t^{\pm 1}\right]$ are
unique factorization domains, since any principal ideal domain is
a unique factorization domain.  We compute $\gcd\left( E_{1}\left(
H_{1}\left( X,x_{0};\mathbb{K}_{0}\left[  t^{\pm 1}\right] \right)
\right) \right)  =\gcd\left\{  \nu\left( d_{1}\right)
,\ldots,\nu\left( d_{k}\right)  \right\}  $.

Since $\nu$ is an embedding  can checks that the degrees of
$\nu\left( \gcd\left\{ d_{1},\ldots,d_{k}\right\}  \right)$ and
$\gcd\left\{ \nu\left( d_{1}\right)  ,\ldots,\nu\left(
d_{k}\right) \right\}$ are equal. It follows that
\begin{align*}
\left\|  \psi\right\|  _{A} &= \deg\nu\left(  \Delta_{X}\right) \\
&=\deg\nu\left(  \gcd\left\{  d_{1},\ldots,d_{k}\right\}  \right) \\
&=\deg\gcd\left\{  \nu\left(  d_{1}\right)  ,\ldots,\nu\left(
d_{k}\right)
\right\} \\
&=\deg\gcd\left\{  E_{1}\left(  H_{1}\left(
X,x_{0};\mathbb{K}_{0}\left[ t^{\pm1}\right]  \right)  \right)
\right\}
\end{align*}
so to complete the proof it suffices to show that $\deg\left(
\gcd\left\{ \nu\left(  d_{1}\right)  ,\ldots,\nu\left(
d_{k}\right)  \right\}  \right) =\overline{\delta}_{0}\left(
\psi\right)  $. Since $\mathbb{K}_{0}\left[  t^{\pm1}\right]$ is a
principal ideal domain, $H_1\left(X,x_0;\mathbb{K}_{0}\left[
t^{\pm1}\right] \right)$ is isomorphic to a direct sum of cyclic
$\mathbb{K}_{0}\left[ t^{\pm1}\right]$-modules.  That is,
$M^{\nu}$ is equivalent to a matrix of the form
\[
\left(
\begin{tabular}[c]{cccc}
  $p_{1}\left(t\right)$ &  &   \\
    & $\ddots$ & \\
    & & $p_{s-1}\left(  t\right)$\\
  0 & $\cdots$ & 0
\end{tabular}\right)
\]
where $p_{i}\left(  t\right)$ is zero for some $i$ if and only if
$r_{0}\left(  X\right)  >0$. We note that the last row of the
matrix can be assumed to be zero since the
$\operatorname{rk}_{\mathcal{K}_{0}}H_{1}\left(
X,x_{0};\mathcal{K} _{0}\right)
=\operatorname{rk}_{\mathcal{K}_{0}}H_{1}\left(
X;\mathcal{K}_{0}\right) +1$. Hence if $r_{0}\left( X\right) =0$,
$H_{1}\left( X;\mathbb{K} _{0}\left[ t^{\pm1}\right] \right)
=TH_{1}\left( X;\mathbb{K}_{0}\left[ t^{\pm1}\right] \right)
=TH_{1}\left( X,x_{0};\mathbb{Z}\Gamma _{0} \right) $ so
$\overline{\delta}_{0}\left( \psi\right) =\deg\left( p_{1}\left(
t\right)  \cdots p_{s-1}\left(  t\right)  \right)  $. Otherwise
$p_{i}\left( t\right)  =0$ for some $i$ so we have $\overline
{\delta}_{0}\left(  \psi\right)  =0=\deg\left( p_{1}\left(
t\right)  \cdots p_{s-1}\left(  t\right)  \right)  $. Using the
latter presentation of $H_{1}\left( X,x_{0};\mathbb{K}_{0}\left[
t^{\pm1}\right]  \right)  $ we compute $\gcd\left(  E_{1}\left(
H_{1}\left( X,x_{0};\mathbb{K}_{0}\left[ t^{\pm1}\right]  \right)
\right) \right)  =p_{1}\left(  t\right)  \cdots p_{s-1}\left(
t\right)  $ so
\begin{align*}
\left\|  \psi\right\|  _{A}
&=\deg\gcd\left\{  E_{1}\left( H_{1}\left(
X,x_{0};\mathbb{K}_{0}\left[ t^{\pm1}\right]  \right)
\right) \right\} \\
&=\deg\left(  p_{1}\left(  t\right)  \cdots p_{s-1}\left(  t\right)  \right) \\
&=\overline{\delta}_{0}\left(  \psi\right)  .
\end{align*}
\end{proof}

\section{Computing $\delta_{i}$ and $\mathcal{A}_n^{\psi}$ via Fox's Free Calculus\label{sec:foxcalc}}

 We will describe a method of computing the higher-order
invariants using Fox's Free Calculus.  We remark that this is
slightly different than the classically defined free derivatives
because we are working with right (instead of the usual left)
modules.  We refer the reader to \cite{Fo1}, \cite{Fo2}, \cite{CF}
and \cite{He} for more on the Free Calculus (for left modules).

Let $G$ be any finitely presented group with presentation
\[P=\left\langle
x_{1},\ldots,x_{l}|r_{1},\ldots,r_{m}\right\rangle,
\] $F=\left\langle x_{1},\ldots,x_{l}\right\rangle $ be the free
group on $l$ generators and $\chi:F\twoheadrightarrow G$. For each
$x_{i}$ there is a mapping $\frac{\partial}{\partial
x_{i}}:F\twoheadrightarrow\mathbb{Z}F$ called the $i^{th}$ free
derivative. This map is determined by the two conditions
\begin{align*}
\frac{\partial x_{j}}{\partial x_{i}}  &  =\delta_{i,j}\\
\frac{\partial\left(  uv\right)  }{\partial x_{i}}  &
=\frac{\partial u}{\partial x_{i}}+u\frac{\partial v}{\partial
x_{i}}.
\end{align*}
From these, one can prove that
\[
\frac{\partial u^{-1}}{\partial x_{i}}=-u^{-1}\frac{\partial
u}{\partial x_{i}}.
\]
The map $\chi:F\twoheadrightarrow G$ extends by linearity to a map
$\chi:\mathbb{Z}F\twoheadrightarrow\mathbb{Z}G$. The matrix
\[
\left(  \frac{\partial r_{j}}{\partial x_{i}}\right)
^{\chi}=\left(
\begin{array}
[c]{ccc} \chi\left(  \frac{\partial r_{1}}{\partial x_{1}}\right)
& \cdots &
\chi\left(  \frac{\partial r_{n}}{\partial x_{1}}\right) \\
\vdots & \ddots & \vdots\\
\chi\left(  \frac{\partial r_{1}}{\partial x_{m}}\right)  & \cdots &
\chi\left(  \frac{\partial r_{n}}{\partial x_{m}}\right)
\end{array}
\right)
\]
with entries in $\mathbb{Z}G$ is called the Jacobian of the
presentation $P$.  We note that this matrix is dependent on the
presentation.

Suppose $X$ be a finite CW-complex, $G\cong\pi _{1}\left(
X,x_{0}\right)$ and $\phi:G\rightarrow\Gamma$.  We can assume that
$X$ has one 0-cell, $x_0$.  Hence the chain complex of $\left(
X_{\Gamma},\widetilde{x}_{0}\right) $ is
\[ \cdots
\rightarrow\mathbb{Z}\Gamma^{m}\overset{\widetilde{\partial}_{2}}{\rightarrow
}\mathbb{Z}\Gamma^{l}\rightarrow0
\]
where $l$ and $m$ are the number of one and two cells of $X$
respectively.  We define an involution on the group ring
$\mathbb{Z}F$ by
\[\overline{\sum{m_i f_i}}= \sum{m_i f^{-1}}\]
and extend $\phi:G\twoheadrightarrow\Gamma$ to
$\phi:\mathbb{Z}G\twoheadrightarrow \mathbb{Z}\Gamma$ by
linearity.  It is straightforward to verify that
$\widetilde{\partial}_{2}=\left( \overline{\frac{\partial
r_{j}}{\partial x_{i} }}\right) ^{\chi\phi}$.  Hence $H_{1}\left(
X,x_{0};\mathbb{Z}\Gamma\right) $ is finitely presented as $\left(
\overline{ \frac{\partial r_{j}}{\partial x_{i}}}\right)
^{\chi\phi}$.  We remark that the existence of the involution in
the presentation of $H_1$ is necessary since we chose to work with
right rather than left modules.  In the case that $\Gamma$ is
abelian, the involution is not necessary.

Let $\iota:\mathbb{Z}\Gamma\rightarrow R$ be a ring homomorphism.
\ Then $ R$ is a $\mathbb{Z}\Gamma$-$ R$-bimodule and we can
consider the right $ R$-module $H_{1}\left(  X,x_{0} ; R\right) $.
The chain complex for $\left(  X,x_{0} ; R\right)  $ is
\[ \cdots
\rightarrow\mathbb{Z}\Gamma^{n}\otimes_{\mathbb{Z}\Gamma} R
\overset{\widetilde{\partial}_{2}\otimes
id_{M}}{\rightarrow}\mathbb{Z}
\Gamma^{m}\otimes_{\mathbb{Z}\Gamma} R\rightarrow0.
\]
Since $\mathbb{Z}\Gamma^{k}\otimes_{\mathbb{Z}\Gamma} R \cong
R^{k}$, it follows that $H_{1}\left( X,x_{0}; R \right) $ is
finitely presented as
\begin{equation}
\left( \overline{ \frac{\partial r_{j}}{\partial x_{i}}}\right)
^{\chi\phi\iota} .\label{foxpres}
\end{equation}

Now let $\phi=\phi_{n}:G\twoheadrightarrow\Gamma_{n}$ be as
defined in Section \ref{sec:rds} and
$\psi:G\twoheadrightarrow\mathbb{Z}$. Choose a splitting
$\xi:\mathbb{Z} \rightarrow\Gamma_{n}\ $\ and let $
R=\mathbb{K}_{n}\left[  t^{\pm 1}\right]  $. We can use
$\left(\ref{foxpres}\right)$ to show that\ $H_{1}\left(
X,x_{0};\mathbb{K}_{n}\left[  t^{\pm1}\right]  \right)  $ is
finitely presented as $\left( \overline{ \frac{\partial
r_{j}}{\partial x_{i}}}\right) ^{\chi\phi_{n}\iota _{\xi}}$ where
$\iota_{\xi}:$
$\mathbb{Z}\Gamma_{n}\hookrightarrow\mathbb{K}_{n}\left[
t^{\pm1}\right]  $ is the embedding of $\mathbb{Z}\Gamma_{n}$ $\
$into $\mathbb{K}_{n}\left[  t^{\pm1}\right]  $.

Moreover, we compute the $\delta_n\left(\psi\right)$ as follows.
Since $\mathbb{K}_{n}\left[ t^{\pm1}\right] $ is a principal ideal
domain, $\left( \overline{ \frac{\partial r_{j}}{\partial
x_{i}}}\right)  ^{\chi\phi_{n} \iota_{\xi}}$ is equivalent to a
diagonal presentation matrix of the form $ \{p_{1}\left(
t\right),\ldots,p_{\lambda}\left( t\right)
 ,0_{\left( r,s\right) }
\}$
%
\cite[Theorem 16, p. 43]{Ja} and $0_{\left( r,s\right)  }$ is a
$r\times s$ size matrix of zeros.  Proposition \ref{comp} implies
that $r_{n}\left( X\right) =\operatorname{rk}_{\emph{K}
_{n}}H_{1}\left( X,x_{0};\emph{K}_{n}\right)  -1$ hence
\[
r_{n}\left(  X\right)  =r-1.
\]

The above presentation implies that $TH_{1}\left(
X,x_{0};\mathbb{K}_{n}\left[ t^{\pm1}\right] \right)  $ has a
diagonal presentation matrix of the form $\left\{ p_{1}\left(
t\right) ,\ldots,p_{\lambda}\left( t\right) \right\} $. Moreover,
Proposition \ref{comp} gives $TH_{1}\left( X;\mathbb{K}_{n}\left[
t^{\pm1}\right] \right)  \cong TH_{1}\left(
X,x_{0};\mathbb{K}_{n}\left[ t^{\pm1}\right] \right) $. Thus we
have used Fox's Free Calculus to derive a presentation matrix for
$\mathcal{A}_n^{\psi}$ and we have shown that
\[
\delta_{n}\left(  \psi\right)  =\deg \prod_{1\leq
i\leq\lambda}p_{i}\left( t\right)  .
\]

\section{A Presentation of $\overline{\mathcal{A}}_{n}^{\psi}
$ in Terms of a Surface Dual to
$\psi$}\label{sec:presmatrixsection}

In the previous section, we used Fox's Free Calculus to find a
presentation matrix of the higher-order Alexander module,
$\mathcal{A}_{n}^{\psi}\left(  X\right) $.  When $X$ is a
3-manifold, we will show that the localized modules
$\overline{\mathcal{A}}_{n}^{\psi}\left( X\right) $ are finitely
presented and that the presentation matrix has topological
significance.  The matrix will depend on the surface dual to a
cohomology class. The presentation will be the higher-order analog
of the presentation obtained from a Seifert matrix for knot
complements.  The presentation obtained will be the main tool that
we use in Section \ref{sec:reltothurstonnorm} to prove that the
higher-order degrees give lower bounds for the Thurston norm.

Let $X^{3}$ be a compact, orientable 3-manifold (possibly with
boundary), $G=\pi_{1}\left( X,x_0\right) $ and $\psi \in
H^1\left(X;\mathbb{Z}\right)$. Let $\phi:G\rightarrow\Gamma$ be a
non-trivial coefficient system and
$X_{\Gamma}\overset{p}{\rightarrow}X$ be the regular $\Gamma$
cover of $X$.  For any $\psi$ as above, there exists a properly
embedded surface $F$ in $X$ such that the class $\left[F\right]
\in H_2\left(X,\partial X;\mathbb{Z}\right)$ is Poincare dual to
$\psi$.  We say that $F$ is dual to $\psi$.

\begin{figure}[htbp]
\begin{center}
\includegraphics
{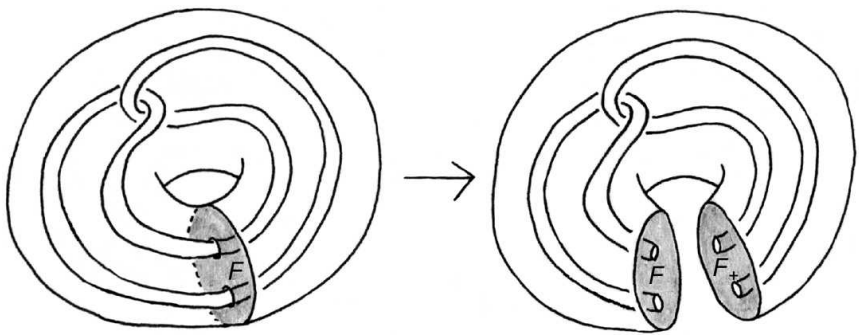}
 \caption{The Whitehead manifold cut open along $F$}
\label{whitehead}
\end{center}
\end{figure}

Let $F$ be a surface dual to $\psi$, $Y=X-\left(F \times
\left(0,1\right)\right)$, $F_+=F \times \{1\}$ (see Figure
\ref{whitehead} for an example), and $x_0$ be a point of $F=F
\times \{0\} \subset Y$.
Let $R$ be a ring and $\tau : \mathbb{Z} \Gamma \rightarrow R$ be
a ring homomorphism defining $R$ as a $\mathbb{Z}\Gamma$-bimodule.
We will exhibit a presentation of $H_{1}\left( X;R\right) $ in
terms of $H_{1}\left( F;R\right) $ and $H_{1}\left( Y;R\right)$.
First we remark that it makes sense to speak of the homology of
$F$ with coefficients in $R$. By this, we mean the homology
corresponding to the coefficient system $\pi_{1}\left( F,x_0
\right) \overset{i_{\ast}}{\rightarrow}\pi_{1}\left( X,x_0 \right)
\overset{\phi }{\rightarrow}\Gamma$. Similarly, we have
coefficient systems for $\pi_1 \left(Y, x_0 \right)$ and the other
terms that are involved in Proposition \ref{genseq} below.

Before we state Proposition \ref{genseq}, we need the following
notation. Let $c$ be a path in $Y$ with initial point
$c\left(0\right) = \left(x_0,1\right)$ and endpoint
$c\left(1\right) = \left(x_0,0\right)$. Let $c_+\left(s\right) =
\left(x_0,s\right)$ and $\alpha$ be the closed curve $c_0 \cdot c$
based at $x_0$. Let $i_\pm:F \rightarrow Y$ include $F$ into $Y$
by $i_-\left(f\right)=\left(f,0 \right)$ and
$i_+\left(f\right)=\left(f,1 \right)$.  Finally let $j:Y
\rightarrow X$ be the inclusion of $Y$ into $X$.

\begin{proposition}
\label{genseq}Suppose $\gamma=\phi\left(\alpha\right)$ is a
non-zero element of $\Gamma$ and either some element of the
augmentation ideal of $\mathbb{Z}\left[\pi_1\left(F\right)\right]$
is invertible (under $\tau \circ \phi \circ i_{\ast}$) in $R$ or
$\pi_1\left(F\right)=1$.  Then the sequence
\[H_{1}\left(  F;R\right) \overset{\eta}{\rightarrow}H_{1}\left(
Y;R\right) \overset{j_{\ast}}{\rightarrow}H_{1}\left( X;R\right)
\rightarrow H_{0}\left(  F \cup F_{+} \cup c;R \right)
\]
is exact where $\eta=\left(i_{+}\right)_{\ast}  -
\left(i_-\right)_{\ast} \gamma$ .
\end{proposition}

\begin{proof}
For convenience, we will omit the $R$ in $H_1\left(-;R\right)$ in
this proof. Let $U= \left(F \times I \right) \cup \alpha$ where
$I=\left[0,1\right]$. Then $X = U \cup Y$ and $U \cap Y = F \cup
F_+ \cup c$. Consider the homology Mayer-Vietoris sequence for
$\left(U,Y\right)$ with coefficients in $R$ \cite{W},
\begin{align}
H_{1}\left(  F \cup F_{+} \cup c \right) \rightarrow H_{1}\left(
\left(F\times I \right) \cup \alpha  \right)  \oplus H_{1}\left( Y
\right) \rightarrow &H_{1}\left(  X  \right)  \rightarrow
\nonumber
\\ \label{mvseq}
&H_{0}\left(  F \cup F_{+} \cup c  \right) \rightarrow.
\end{align}
We examine the $H_1$ terms involving $F$ in (\ref{mvseq}). We will
compute the homology of these term using the Mayer-Vietoris
sequences for $\left(F,F_+ \cup c\right)$,
\begin{align}
0 \rightarrow H_{1}\left( F \right) \oplus H_{1}\left( F_+ \cup c
\right) \rightarrow H_{1}\left(  F \cup F_+ \cup c  \right)
\rightarrow  &H_{0}\left(  x_0 \right) \rightarrow \nonumber
\\ \label{mvseq1}
&H_0\left( F \right) \oplus H_0\left( F_+ \cup c \right)
\end{align}
and $\left(F\times I ,\alpha \right)$,
\begin{align}
0 \rightarrow H_{1}\left( F \times I \right) \oplus H_{1}\left(
\alpha \right) \rightarrow H_{1}\left(  F \times I \cup \alpha
\right) \rightarrow &
 H_{0}\left(  x_0 \right) \rightarrow \nonumber
 \\ \label{mvseq2}
 &H_0\left( F \times
I \right) \oplus H_0\left( \alpha \right).
\end{align}

The idea behind the rest of the proof in both of the cases (stated
in the hypothesis) are similar however the proof when
$\pi_1\left(F\right)\neq 1$ is more technical.  Hence we will
first consider the special case when $\pi_1\left(F\right)=1$.
Since $\pi_1\left(F\right)$ is trivial, both $H_{0}\left( x_0
\right) \rightarrow H_0\left( F \right) \oplus H_0\left( F_+ \cup
c \right)$ and $H_{0}\left( x_0 \right) \rightarrow H_0\left( F
\times I \right) \oplus H_0\left( \alpha \right)$ are injective.
Hence, $H_{1}\left( F \cup F_{+} \cup c;\mathcal{M}\right) \cong
H_{1}\left(  F ;\mathcal{M}\right) \oplus H_{1}\left( F_{+} \cup
c;\mathcal{M}\right)$ and $H_{1}\left( \left(F\times I \right)
\cup \alpha
 \right) \cong H_{1}\left( \left(F\times I \right)
 \right) \oplus H_1\left(\alpha \right)$.

Since $\gamma$ is non-trivial in $\Gamma$, the curve $\alpha$ does
not lift to the $\Gamma$-cover of $X$.  Therefore $H_1\left(\alpha
\right)=0$ and hence $H_{1}\left( \left(F\times I \cup \alpha
\right)  \right) \cong H_{1}\left( F\times I \right)$.  Moreover,
$H_{1}\left( \left(F\times I \right)  \right) \cong H_{1}\left( F
 \right)$ where the isomorphism is induced by the map which sends
$\left(f,s\right)$ to $\left(f,0\right)$.

We analyze the first term in the sequence.  The isomorphism
$\pi_1\left(F,x_0\right) \rightarrow \pi_1\left(F_{+} \cup
c,x_0\right)$ given by $\left[\beta\right] \mapsto \left[c \cdot
i_+\left(\beta\right) \cdot \overline{c}\right]$ induces an
isomorphism $g:H_{1}\left( F \right) \rightarrow H_{1}\left(F_{+}
\cup c \right)$. By $\overline{c}$ we mean the curve defined by
$\overline{c}\left(s\right)=c\left(1-s\right)$. Note that $\left[c
\cdot i_+\left(\beta\right) \cdot \overline{c}\right]=\alpha^{-1}
\left[\beta\right] \alpha$ in $\pi_1\left(F \times
 I ,x_0\right)$. Therefore $H_{1}\left( F \cup F_{+}
\cup c \right) \cong H_1\left(F \right) \oplus H_1\left(F
\right)$.
 We note that the composition
\[ H_1\left(F \right) \rightarrow
H_1\left(F_+ \cup c \right) \rightarrow H_1\left(F \times I \cup
\alpha \right) \rightarrow H_1\left(F \right)
\]
sends $\sigma$ to $\sigma \gamma$ and the composition $H_1\left(F
\right) \rightarrow H_1\left(F \times I \cup \alpha \right)
\rightarrow H_1\left(F \right)$ is the identity.

Using the isomorphisms above, we rewrite
$\left(\ref{mvseq}\right)$ as
\begin{eqnarray*}
H_{1}\left(  F \right)  \oplus H_{1}\left( F \right)
\overset{\left(  f_{F},f_{Y}\right) ,}{\rightarrow }H_{1}\left( F
\right)  \oplus H_{1}\left(  Y \right)
\rightarrow  H_{1}\left( X \right) \rightarrow
H_{0}\left( F \cup F_+ \cup c \right) \rightarrow
\end{eqnarray*}
where $f_{F}\left(  \sigma_{1},\sigma_{2}\right) =\sigma_{1} +
\sigma_{2} \gamma $ and $f_{Y}\left( \sigma_{1},\sigma_{2}\right)
=-\left( \left(  i_{-}\right) _{\ast}\left(  \sigma_{1}\right)
+\left( i_{+}\right) _{\ast}\left( \sigma_{2}\right)  \right)  $.
It follows from Lemma \ref{new_mv} that the sequence
\[
H_{1}\left(  F \right) \overset{\eta}{\rightarrow}H_{1}\left( Y
\right) \overset{j_{\ast}}{\rightarrow}H_{1}\left( X \right)
\rightarrow H_{0}\left(  F_{-}\cup F_{+}  \right)
\]
is exact with $\eta\left(  \sigma\right)
=-\left(  \left(  i_{-}\right)
_{\ast}\left( \sigma \gamma\right) +\left(  i_{+}\right)
_{\ast}\left( -\sigma  \right) \right)  =\left( \left(  i_{+}
\right) _{\ast}  -\left( i_{-}\right) _{\ast} \gamma \right)
\left( \sigma\right) $.  The proof of Lemma \ref{new_mv} is
straightforward hence omitted.

\begin{lemma}
\label{new_mv}Suppose $A\oplus A\overset{\left( f_{A},f_{B}\right)
}{\longrightarrow}A\oplus
B\overset{g}{\rightarrow}C\overset{h}{\rightarrow}D$ is an exact
sequence of right $ R$-modules with $f_{A}\left(
a_{1},a_{2}\right)  =a_{1}+a_{2} r$ where $r \in R$, then $
A\overset{f^{\prime}}{\rightarrow}B\overset{g^{\prime}}{\rightarrow}
C\overset{h}{\rightarrow}D $ is an exact sequence of right
$R$-modules with $f^{\prime }\left( a\right)  \equiv f_{B}\left( a
r,-a r\right)  $ and $g^{\prime}\left( b\right)  \equiv g\left(
0,b\right)  $.
\end{lemma}

Now we assume that some element of the augmentation ideal of
$\mathbb{Z}\left[\pi_1\left(F\right)\right]$ is invertible in $R$.
By \cite[p. 275]{W} and \cite[p.34]{B}, $H_0\left(F \right)$ is
isomorphic to the cofixed set $R/RJ$ where $J$ is the augmentation
ideal of $\mathbb{Z}\left[\pi_1\left(F\right)\right]$.  Therefore
$H_0\left(F \right)=H_0\left(F_+ \cup c \right)=0$.  We note that
$H_0\left(x_0 \right)$ is the free $R$-module of rank one
generated by $\left[x_0\right]$.  Choose a splitting $\xi_1$ for
short exact sequence in (\ref{mvseq1}) to get
\begin{equation}\label{iso1}
H_1\left(F \cup F_+ \cup c  \right) \cong M \oplus H_1\left(F
\right) \oplus H_1\left(F_+ \cup c \right)
\end{equation}
where $M$ is the free $R$-module of rank one generated by
$\xi_1\left(\left[x_0\right]\right)$. Let $\beta$ be a curve in $F
\cup F_+ \cup c$ representing
$\xi_1\left(\left[x_0\right]\right)$.

Consider the sequence in (\ref{mvseq2}).  We note that
$H_0\left(\alpha \right)=R/{\langle \tau\left(\gamma\right) - 1
\rangle}$. Since $\gamma$ is non-trivial,
$\operatorname{Im}\left(H_{1}\left( F \times I \cup \alpha
\right)  \rightarrow
 H_{0}\left(  x_0 \right)\right)
$ is the free $R$-module of rank one generated by
$\left(\tau\left(\gamma\right)-1\right)\left[x_0\right]$.
Moreover, $\left[\beta\right] \mapsto
u\left(\tau\left(\gamma\right)-1\right)\left[x_0\right]$ where $u$
is a unit of $R$ under the boundary homomorphism.  We choose the
splitting
$\xi_2:u\left(\tau\left(\gamma\right)-1\right)\left[x_0\right]
\mapsto \left[\beta\right]$ to get
\begin{equation}\label{iso2}
H_1\left(F \times I \cup \alpha  \right) \cong N \oplus H_1\left(F
\times I \right) \oplus H_1\left(\alpha \right)
\end{equation}
where $N$ is the free $R$-module of rank one generated by
$\left[\beta\right]$.

Using the isomorphisms in (\ref{iso1}) and (\ref{iso2}), we can
rewrite  (\ref{mvseq}) as
\begin{align*}
M \oplus H_{1}\left(F \right) \oplus H_1\left(F_{+} \cup c \right)
\rightarrow N \oplus H_{1}\left(F\times I \right) \oplus
H_1\left(\alpha \right)\oplus & H_1\left(Y  \right) \rightarrow
H_{1}\left( X  \right) \rightarrow \\ & H_{0}\left(  F \cup F_{+}
\cup c  \right) \rightarrow.
\end{align*}
We use Lemma \ref{lem:tech} to get an exact sequence without the
$M$ and $N$ terms as in the case when $\pi_1\left(F\right)=1$. To
complete the proof of the proposition, we use the same argument as
in the case when $\pi_1\left(F\right)=1$.
\end{proof}

\begin{lemma}\label{lem:tech}
Let $M$ and $N$ be free right $R$-modules of rank one generated by
$m$ and $n$ respectively.  Suppose $M \oplus A
\overset{f}{\rightarrow}  N \oplus B \overset{g}{\rightarrow}C
\overset{h^{\prime}}{\rightarrow} D$ is an exact sequence of right
$R$-modules with
$f\left(rm,a\right)=\left(rn,f_2\left(rm,a\right)\right)$ for some
$f_2 : M \oplus A \rightarrow B$.  Let $\eta : 0 \oplus B
\rightarrow B$ be the isomorphism defined by $\left(0,b\right)
\mapsto b$.  Then $A \overset{f^{\prime}}{\rightarrow}  B
\overset{g^{\prime}}{\rightarrow}C
\overset{h^{\prime}}{\rightarrow} D $ is an exact sequence of
right $R$-modules where $f^{\prime}$ and $g^{\prime}$ are defined
by $f^{\prime}\left(a\right)=\eta\left(f\left(0,a\right)\right)$
and $g^{\prime}\left(b\right) = \left(g\left(0,b\right)\right)$.
\end{lemma}
\begin{proof}The proof is straightforward hence omitted.
\end{proof}

Now we consider the presentation of $H_{1}\left(
X;\mathbb{K}_{n}\left[ t^{\pm1}\right]  \right)  $ where
$\mathbb{K}_{n}$ is the skew field of fractions of
$\mathbb{Z}\Gamma_{n}^{\prime}$ as defined before. Let $F$ and $Y$
be as defined above.  Since $\pi_{1}\left(  F,x_0\right)  $ and
$\pi_{1}\left(Y,x_0\right)$ are contained in the kernel of $\psi$,
we can consider the homology of $F$ and $Y$ with coefficients in
$\mathbb{Z}\Gamma^{\prime}$ and $\mathbb{K}_n$. Since $F$ and $K$
are finite CW-complexes, $H_1\left(F;\mathbb{K}_n\right)$ and
$H_1\left(Y;\mathbb{K}_n\right)$ are finitely generated free
modules hence are isomorphic to $\mathbb{K}_n^l$ and
$\mathbb{K}_n^m$ respectively.  Thus the $\mathbb{K}_n$-module
homomorphisms $i_\pm :H_1\left(F;\mathbb{K}_n\right) \rightarrow
H_1\left(Y;\mathbb{K}_n\right)$ can be represented by $l \times m$
matrices $V_\pm$ with coefficients in $\mathbb{K}$.  We will show
that the higher-order module corresponding to $\psi$ is presented
by $V_+ - V_- t$.

\begin{proposition}
\label{alexpres} $\operatorname{Im} \left(  j_{\ast}\right)  $ is
finitely presented as $\mathbb{K}_{n}\left[ t^{\pm1}\right]
^{l}\overset{P}{\rightarrow}\mathbb{K}_{n}\left[ t^{\pm 1}\right]
^{m}\overset{j_{\ast}}{\twoheadrightarrow}\operatorname{Im}\left(
j_{\ast}\right)$ where $P=V_+ - V_- t$.  Moreover, if
$\psi|_{\pi_{1}\left( F^{j}\right) }$ is nontrivial for each
component $F^{j}$ of $F$ then $\overline{\mathcal{A}}_{n}^{\psi
}\left( X\right)$ is finitely presented as
\[
\mathbb{K}_{n}\left[  t^{\pm1}\right] ^{l}\overset{P}{\rightarrow}
\mathbb{K}_{n}\left[ t^{\pm1}\right] ^{m}\overset{j_{\ast}}
{\twoheadrightarrow}\overline{\mathcal{A}}_{n}^{\psi}\left(
X\right) .
\]
\end{proposition}

\begin{proof}
Let $\gamma=\phi_n\left(\alpha\right)$ be as in Proposition
\ref{genseq}.  We note that $\psi\left(\gamma\right)=t$.  Choose
the splitting $\xi: t \mapsto \gamma$.  Since
$\pi_1\left(F,x_0\right) \subset \Gamma_n^{\prime}$ we have
\[ C_{\ast}\left(F_{\Gamma_n}\right) \otimes_{\mathbb{Z}\Gamma_n}
\mathbb{K}_n\left[t^{\pm 1}\right] \cong \left(
C_{\ast}\left(F_{\Gamma_n^{\prime}}\right)
\otimes_{\mathbb{Z}\Gamma_n^{\prime}} \mathbb{K}_n \right)
\otimes_{\mathbb{K}_n} \mathbb{K}_n\left[t^{\pm 1}\right] .
\]
Moreover, $\mathbb{K}_n\left[t^{\pm 1}\right]$ is a direct sum of
free $\mathbb{K}_n$-modules $\left( \mathbb{K}_n\left[t^{\pm
1}\right] \cong \oplus_{i=-\infty}^{\infty}\mathbb{K}_n \right)$.
Therefore $\mathbb{K}_n\left[t^{\pm 1}\right]$ is a flat left
$\mathbb{K}_n$-module. Thus
\begin{align*}H_1\left(F;\mathbb{K}_n\left[t^{\pm 1}\right]\right) &\cong
H_1\left(F;\mathbb{K}_n\right) \otimes_{\mathbb{K}_n}
\mathbb{K}_n\left[t^{\pm 1}\right] \\ &\cong \mathbb{K}_{n}^{l}
\otimes_{\mathbb{K}_n} \mathbb{K}_n\left[t^{\pm 1}\right] \\ &
\cong \mathbb{K}_n\left[t^{\pm 1}\right]^l.
\end{align*}
Similarly, we have $H_1\left(F;\mathbb{K}_n\left[t^{\pm
1}\right]\right) \cong \mathbb{K}_n\left[t^{\pm 1}\right]^m$. The
first result follows from Proposition \ref{genseq}. If
$\psi|_{\pi_{1}\left( F^{j}\right) }$ is nontrivial for all $j$
then $H_{0}\left( F \cup F_{+} \cup c;\mathbb{K}_{n}\left[  t^{\pm
1}\right]  \right) =0$ so $\operatorname{Im}\left( j_{\ast}\right)
=\overline{\mathcal{A}}_{n}^{\psi}\left(  X\right) $.
\end{proof}

We use the following lemma to show that it suffices to use the
presentation matrix $V_+ - V_- t$ when computing
$\mathcal{A}_{n}^{\psi }\left( X\right)$.

\begin{lemma}
\label{equaltor}Suppose
$B\overset{g}{\rightarrow}C\overset{h}{\rightarrow}D$ is an exact
sequence of right $ R$-modules where $D$ is $ R$-torsion free and
$ R$ is an Ore domain then $T_R C=T_R \operatorname{Im}\left(
g\right)  $.
\end{lemma}

\begin{proof}
Since $\operatorname{Im}\left(  g\right)  \subseteq C$, it is easy
to verify $T\operatorname{Im}\left(  g\right)  \subseteq TC$. Let
$c\in TC$ then there exists a non-zero $r\in R$ such that $cr=0$.
This says that $h\left(  c\right)  r=h\left(  cr\right) =0$ so
that $h\left(  c\right)  $ is $ R$-torsion in $D$ hence $h\left(
c\right)  =0$. By exactness at $C$ we see that
$c\in\operatorname{Im}\left(  g\right)  $ and $cr=0$ so
$TC\subseteq T\operatorname{Im}\left(  g\right)  $.
\end{proof}

\begin{proposition}\label{pro:pres_A}
$\mathcal{A}_{n}^{\psi }\left( X\right)$ is isomorphic to the
$\mathbb{K}_n\left[t^{\pm 1}\right]$-torsion submodule of \\
$\operatorname{cok}\left(V_+ - V_- t\right)$.
\end{proposition}
\begin{proof}Recall that $\mathcal{A}_{n}^{\psi }\left( X\right) \cong
T_{\mathbb{K}_n\left[t^{\pm 1}\right]}
\overline{\mathcal{A}}_{n}^{\psi }\left( X\right)$. The result
follows immediately from Lemma \ref{equaltor} and Proposition
\ref{alexpres}.\end{proof}

\section{Examples}\label{calcexamples}

In this section we will compute $r_{n}$, $\delta_{n}$, and
$\overline{\delta}_{n}$ for some well known 3-manifolds and relate
their values to those given by the Thurston norm. In each of the
examples we denote the fundamental group of $X$ by $G$.  Of
particular importance will be the 3-manifolds which fiber over
$S^1$ and those which are Seifert fibered.  We start with the
standard examples.

\begin{example}3-torus

Let $X=S^{1}\times S^{1}\times S^{1}$ then $G_{r}^{\left( 1\right)
}=\left\{  1\right\}  $ hence $r_{n}\left(  X\right)
=\delta_{n}\left( \psi\right)  =0$ for all $\psi$ and $n\geq0$.
Note that since $H_{2}\left( X;\mathbb{Z}\right)  $ is generated
by tori, the Thurston norm is zero for all $\psi \in
H^1\left(X,;\mathbb{Z}\right)$. More generally, if $G$ is any
finitely generated abelian group (with $\beta_{1}\left( G\right)
\geq1$) then $G_{r}^{\left( n\right) }=T$ where $T$ is the torsion
subgroup of $G$.  Hence $r_{n}\left( X\right) =\delta_{n}\left(
\psi\right) =0$.
\end{example}

\begin{example}
$X=\underset{i=1}{\overset{m}{\#}}S^{2}\times S^{1}$

Let $X=\#_{i=1}^{m}S^{2}\times S^{1}$ for $m\geq1$ then $G=F_{m}$,
the free group on $m$ generators. Since $H_{2}\left(
X;\mathbb{Z}\right)  $ is generated by spheres $\left\|
\psi\right\|  _{T}=0$ for all $\psi$. Moreover, every class in
$H_2\left(X;\mathbb{Z}\right)$ can be represented by a disjoint
union of embedded spheres.  Hence there exists a surface $F$ dual
to $\psi$ such that $H_{1}\left( F;\mathbb{K}_{n}\left[
t^{\pm1}\right] \right) =0$. By Proposition \ref{pro:pres_A} we
have $\mathcal{A}^{\psi}_n\left(X\right)
=T\operatorname{cok}\left(\mathbf{0}\right)$. Therefore
$\delta_{n}\left( \psi\right)  =0$ for all $\psi$ and $n\geq0 $.

Since $r_{n}$ only depends on the group $G$ we can assume that $X$
is a wedge of $m$ circles.  By Remark \ref{eulerchar},
$1-m=\chi\left(X\right)=\sum_{i=0}^{1} rk_{\mathcal{K}_n}
H_i\left(X;\mathcal{K}_n\right)$. Since $\phi_{n}$ is a
non-trivial homomorphism, by Lemma \ref{lem:H0},
$H_0\left(X;\mathcal{K}_n\right)=0$. Therefore $r_{n}\left(
X\right)=m-1$.
\end{example}

There is a large class of 3-manifolds for which $r_0\left(X\right)
\geq 1$. Recall that a \emph{boundary link} is a link $L$ in $S^3$
such that the components admit mutually disjoint Seifert surfaces.
It is easy to see that each of these surfaces lifts to the
universal abelian cover of $X=S^3-L$. By Proposition
\ref{positiverank}, $r_0\left(X\right) \geq 1$.  Hence the
Alexander norm for $X$ is always trivial. It is often true that
the refined Alexander norm, $\delta_0$, is non-trivial.  We
compute an example below where this is the case.

\begin{figure}[h]
\begin{center}
\includegraphics{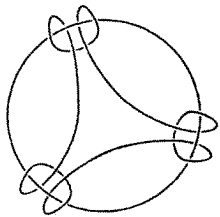}
\caption{Each component of $L$ has minimal ropelength at least
$2\pi\left(1+\sqrt 3\right)$.} \label{ropelengthex}
\end{center}
\end{figure}

Let $L$ be the link pictured in Figure \ref{ropelengthex}.
Let $F$
be a Seifert surface of one of the components of $L$ as in Figure
representing the minimal $\chi\_$ and $\psi$ be dual to $F$. We
will show in Example \ref{ex:ropelength} that $\delta_0 \left(
\psi \right)=4$. Moreover, each component bounds a once punctured
genus two surface hence $\left\| \psi \right\|_T \leq 3$.  Hence
by Corollary \ref{cor:linkbound}, we get
\[\delta_0 \leq \left\| \psi \right\|_T + 1.\]  We conclude that $\left\| \psi
\right\|_T=3$.  Hence, even for $n=0$, $\delta_n$ gives a sharper
bound for the Thurston norm that the Alexander norm.  This also
shows that the minimal ropelength of each of the components of $L$
is at least $2\pi\left(1+\sqrt 3\right)$. For more information on
the ropelength of knots and links, see \cite[Corollary 22]{CKS}.
\begin{example}\label{ex:ropelength}
Let $L$ be the link in Figure \ref{ropelengthex} and $X=S^3-L$.
 We
use the techniques of Section \ref{sec:foxcalc} to compute
$\delta_0\left(\psi\right)$. Using a Wirtinger presentation, we
present $G=\pi\left(X\right)$ as
\begin{align*}
\left< a,b,c,d,e,f,g,h,i,j,k,l | \right. b g^{-1} i c^{-1} i^{-1}
g, c j^{-1} l a^{-1} l^{-1} j,  f e^{-1} h g^{-1} h^{-1} e,
\\ i h^{-1} k j^{-1} k^{-1} h, l k^{-1} e d^{-1}e^{-1} k,
d a^{-1} e^{-1} a, e b f^{-1} b^{-1},\\ g b^{-1} h^{-1} b,  h c
i^{-1} c^{-1}, j c^{-1} k^{-1} c, k a l^{-1} a^{-1} \left.
\right>.  \end{align*} Using Fox's Free Calculus we obtain a
presentation matrix $M$ for $H_1\left(X_0,\tilde{x}_0\right)$ (see
below). Here, $x$ is the abelianization of $a$ and $y$ is the
abelianization of $d$. Since we used a Wirtinger presentation for
$G$, $x$ and $y$ represent the meridians of $L$.
\[\left(
\begin{array}{c@{\;\;}c@{\;\;}c@{\;\;}c@{\;\;}c@{\;\;}c@{\;\;}c@{\;\;}c@{\;\;}c@{\;\;}c@{\;\;}c}
0 &  -y  &   0 &  0 &  0 &  1-y &  0 &  0 &  0 &  0  & y-1 \\
y & 0 &  0 &  0 &  0 &  0 &  y-1 &  1-y &  0 &  0 &  0\\-y & y & 0
& 0 &  0 &  0 &  0 &  0 &  y-1 &  1-y &  0\\0 &  0 &  0 & 0 & -y &
x &  0 &  0 &  0 &  0 &  0\\0 &  0 &  1-y &  0 & y-1 & -1 & 1 & 0
& 0 &  0 &  0\\0 &  0 &  y &  0 &  0 & 0 & -x & 0 & 0 & 0 & 0\\1-x
&  0 &  -y &  0 &  0 &  0 & 0 & x &  0 &  0 & 0\\0 & 0 & y-1 & 1-y
&  0 &  0 &  0 & -1 & 1 &  0 &  0\\x-1 & 0 &  0 & y & 0 & 0 &  0 &
0 &  -x & 0 &  0\\0 & 1-x & 0 & -y &  0 &  0 & 0 & 0 & 0 &  x &
0\\0 & 0 &  0 &  y-1 & 1-y & 0 & 0 &  0 & 0 & -1 & 1\\0 & x-1 &  0
& 0 &  y & 0 &  0 & 0 & 0 & 0 &  -x
\end{array}\right)
\]
 This is equivalent (using the moves in Lemma \ref{pres})
to the matrix
\[\begin{pmatrix}
1-x-y & 0 &  0 & 0 \\
0 & 1-x-y & 0  & 0\\
0 & 0 & xy-x-y & 0 \\
0 & 0 & 0 & xy-x-y  \\
 0 & 0 & 0 & 0 \\
 0 & 0 & 0 & 0\end{pmatrix}.
\]
Hence $r_0\left(X\right)=1$ and \[\mathcal{A}_0 \left(X\right)=
\frac{\mathbb{Z}\left[x^{\pm1},y^{\pm1}\right]}{\left<1-x-y\right>}
\oplus
\frac{\mathbb{Z}\left[x^{\pm1},y^{\pm1}\right]}{\left<1-x-y\right>}\oplus
\frac{\mathbb{Z}\left[x^{\pm1},y^{\pm1}\right]}{\left<xy-x-y\right>}\oplus
\frac{\mathbb{Z}\left[x^{\pm1},y^{\pm1}\right]}{\left<xy-x-y\right>}
.\]Let $\psi$ be dual to a Seifert surface for one of the
components of $L$. Then $\psi$ maps one of the generators of
$H_1\left(X\right)$ to $t$ and the other to 1. The link is
symmetric so either choice will suffice.  Say $x \mapsto t$ and $y
\mapsto 1$. Choose the splitting $t \mapsto x$.  Each of the
polynomials in the latter matrix has degree 1 in
$\mathbb{K}_0\left[t^{\pm1}\right]$ since $1-x-y \mapsto
\left(1-y\right) + t$  and $xy - x - y \mapsto t\left(y-1\right)
-y$. Therefore $\delta_0\left(\psi\right)=4$ as desired. In fact,
if $\psi$ maps $x \mapsto t^m$ and $y \mapsto t^n$ then
\begin{align*}
\delta_0\left(\psi\right)= \operatorname{deg} (&t^{2n}+ t^{3m+n}+
t^{2m+2n}+ t^{2m+n}+t^{m+2*n}+ t^{m+n}+ t^{2m}+ \\
 &t^{3m}+
t^{3n}+t^{4m}+t^{4n}+t^{3m+2n}+t^{2m+3n}+t^{m+3n}+t^{4m+2n}+
\\ &t^{4m+n}+t^{3m+3n}+t^{4n+2m}+t^{4 n+m}).
\end{align*}whereas $\overline{\delta}_0\left(\psi\right)=0$.
\end{example}

Although the invariants are defined algebraically, they respect
much of the topology of the 3-manifold. We begin by considering
those 3-manifolds which fiber over $S^1$. In this case the
higher-order invariants behave in a very special manner.

\begin{proposition}
\label{fibered}If $X$ is a compact, orientable 3-manifold that
fibers over $S^{1}$ then
\[
r_{n}\left(  X\right) = 0.
\]
Let $\psi$ be dual to a fibered surface.  If $n=0$,
$\beta_{1}\left( X\right)=1$, $X \ncong S^1 \times S^2$ and $X
\ncong S^1 \times D^2$ then $\delta_{n}\left( \psi\right)=\left\|
\psi\right\| _{T}+1+\beta_{3}\left( X\right)$. Otherwise,
\[
\delta_{n}\left(  \psi\right)   = \left\|  \psi\right\|  _{T}.
\]
\end{proposition}

\begin{proof}
Consider $X\rightarrow S^{1}$ with fiber surface $F$ and $\psi$ be
the element of $H^{1}\left(  X;\mathbb{Z}\right)  $ which is dual
to $F$. The $\Gamma_{n} $-cover of $X$ factors through the
infinite cyclic cover corresponding to $\psi$ with total space
$F\times\mathbb{R}$. Hence $X_{n}$ is homeomorphic to
$F_{n}\times\mathbb{R}$ where $F_{n}$ is a regular cover of $F$.
It follows that $H_{\ast}\left(  X;\mathbb{K}_{n}\left[
t^{\pm1}\right] \right)  $ is isomorphic to $H_{\ast}\left(
F;\mathbb{K}_{n}\right) $ as a $\mathbb{K} _{n}$-module.  In
particular, $\delta_n\left(\psi\right) =
\operatorname{rk}_{\mathbb{K}_n}H_1\left(F;\mathbb{K}_n\right)$.
Moreover,  since $H_{1}\left( F;\mathbb{K}_{n}\right) $ is a
finitely generated $\mathbb{K}_{n}$-module,  $r_{n}\left( X\right)
=0$ for all $n$. That is, $H_{1}\left( X;\mathbb{K}_{n}\left[
t^{\pm1}\right] \right)  $ is a torsion module for all $n\geq0$.

We restrict to the case that $n=0$ and $\beta_{1}\left(
X\right)=1$.  We note that $F_{0}=F$ and $\mathbb{K}
_{0}=\mathbb{Q}$ so that
$\operatorname{rk}_{\mathbb{K}_{0}}H_{1}\left( F;\mathbb{K}
_{0}\right) \cong\mathbb{\beta}_{1}\left( F\right) $.   Thus
$\delta_{0}\left( \psi\right)  =\beta_{1}\left( F\right) =
-\chi\left(F\right) +1+\beta_{3}\left(  X\right)  $.  As long as
$X \ncong S^1 \times S^2$ and $X \ncong S^1 \times D^2$ the Euler
characteristic of $F$ is non-positive hence $\left\|
\psi\right\|_T=-\chi\left(F\right)$.   Therefore
$\delta_0\left(\psi\right)=\left\| \psi\right\|_T
+1+\beta_{3}\left(  X\right)$.

Note that if the Euler characteristic of $F$ is ever positive then
$\pi_1\left(F\right)=1$.  Thus we have
$H_1\left(F;\mathbb{K}_n\right)=0$ for all $n\geq 0$.  Therefore
both $\delta_n\left(\psi\right)$ and $\left\|\psi\right\|$ are
zero and hence equal for all $n \geq 0$.

Otherwise, $F_n$ factors through a non-trivial free abelian cover
of $F$.  By Lemma \ref{lem:H0}, $H_{0}\left(
F;\mathbb{K}_{n}\right)=0$.  Since $F_n$ is noncompact,
$H_{2}\left( F;\mathbb{K}_{n}\right) = 0$. It follows that
$\delta_{n}\left(
\psi\right)=\operatorname{rk}_{\mathbb{K}_{n}}H_{1}\left(
F;\mathbb{K}_{n}\right) = -\chi\left( F\right)  $.  Moreover, $F$
has non-positive Euler characteristic so $-\chi\left(
F\right)=\left\| \psi\right\| _{T} $.
\end{proof}

As with the first two examples in this section, there is a large
class of 3-manifolds which have vanishing (unrefined) higher-order
degrees. This is the class of Seifert fibered manifolds that do
not fiber over $S^1$. We remark that the condition that $X$ not
fiber over $S^1$ is necessary by the previous proposition. Some
good references on Seifert fibered manifolds are \cite[Ch. 2]{Ha},
\cite[Ch. 12]{H}, and \cite[Ch. VI]{Jaco}.

\begin{proposition}\label{pro:SFS}
Let $X$ be a compact, orientable Seifert fibered manifold that
does not fiber over $S^1$.  If $\beta_1\left(X\right) \geq 2$ or
$n\geq 1$ then
\[\overline{\delta}_{n}\left(
\psi\right) =0
\]
 for all $\psi\in H^{1}\left( X;\mathbb{Z}\right)
$.
\end{proposition}

\begin{proof}
This is most easily proven using Theorem \ref{delbarthm} and some
well known results about Seifert fibered 3-manifolds.  By Theorem
VI.34 of \cite{Jaco}, we see that any two-sided incompressible
surface in $X$ must be a disc, annulus, or a torus. Therefore the
Thurston norm of $X$ is trivial. Theorem \ref{delbarthm} implies
that $\overline{\delta}_{n}\left( \psi\right) \leq \| \psi \|_T
=0$ whenever $\beta_1\left(X\right)$ or $n \geq 1$.
\end{proof}

We end this section by showing that under the connected sum of
3-manifolds, the degrees are additive and the ranks plus 1 are
additive. The following is not at all obvious because the
fundamental groups of the spaces involved are completely
different!

\begin{proposition}\label{connectedsum}
Suppose $X=X_{1}\#X_{2}$, $\beta_{1}\left(  X_{i}\right)  \geq1$
and $\psi\in H^{1}\left(  X;\mathbb{Z}\right)  $. Then \[
r_{n}\left(  X\right)  =r_{n}\left(  X_{1}\right)  +r_{n}\left(
X_{2}\right) +1
\]
and
\[
\delta_{n}\left(  \psi\right)  =\delta_{n}\left(  \psi_{1}\right)  +\delta
_{n}\left(  \psi_{2}\right)
\]
where $\psi=\psi_{1}\oplus\psi_{2}$.
\end{proposition}

\begin{proof}
We begin by showing that $r_{n}\left( X\right) =r_{n}\left(
X_{1}\right) +r_{n}\left( X_{2}\right) +1$. Consider the following
Mayer-Vietoris sequence of $ R $-modules for any ring $ R$ with
$\mathbb{Z}\Gamma_{n}\subseteq
 R\subseteq\mathcal{K}_{n}$.  By $\Gamma_{n}$ we mean
the quotient of $G=\pi_{1}\left(  X\right)  $ by the $\left(
n+1\right)^{st}$ term of the rational derived series of $G$.
\begin{equation}\label{connectsum} 0\rightarrow H_{1}\left(
X_{1}; R\right) \oplus H_{1}\left( X_{2}; R\right)
\overset{\nu}{\rightarrow}H_{1}\left( X; R\right)
\overset{\partial_{1}} {\rightarrow}H_{0}\left( S^{2}; R\right)
 \rightarrow
\end{equation}
We note that $H_{0}\left(  S^{2}; R\right) \cong R$ since $S^{2}$
is simply connected. For $j=1,2$ let $i_{j}:G_{j}\rightarrow G$ be
the inclusion map, $\operatorname{pr}_{j}:G\rightarrow G_{j}$ be
the projection onto the $j^{th}$ factor, and
$\Gamma_{n}^{j}=\left( G_{j}\right) /\left( G_{j}\right)
_{r}^{\left(  n+1\right)  }$ where $G_{j}=\pi_{1}\left(
X_{j}\right)  $. Since $G=G_{1}\ast G_{2}$, $\operatorname{pr}_{j}
\circ i_{j}=id_{G_{j}}$. Hence the induced maps $
\Gamma_{n}^{j}\overset{\overline{i}_{j}}{\rightarrow}\Gamma_{n}\overset
{\overline{\operatorname{pr}}_{j}}{\rightarrow}\Gamma_{n}^{j} $
are also the identity making $\overline{i}_{j}$ a monomorphism.
Thus the $\Gamma_{n}^{j}$ cover of $X_{j}$, $X_{jn}$, can be
constructed as the regular cover corresponding to the map
$\phi_{n}\circ i_{j}:G_{j}
\twoheadrightarrow\operatorname{Im}\left( \phi_{n}\circ
i_{j}\right)  $. We extend $\overline{i}_{j}$ to a ring
monomorphism $\overline{i}_{j}
:\mathbb{Z}\Gamma_{n}^{j}\rightarrow\mathbb{Z}\Gamma_{n}$.

The map $G_{j}\rightarrow G/G_{r}^{\left(  n+1\right)  }$ is the
zero map if and only if $G_{j}/\left(  G_{j}\right)  _{r}^{\left(
n+1\right)  }=0$.  We assumed $\beta_{1}\left( X_{j}\right)  >0$
hence $G_{j}/\left( G_{j}\right)  _{r}^{\left( n+1\right) }\neq0$.
By Lemma \ref{lem:H0}, $H_{0}\left( X_{j};\mathcal{K}_{n}\right)
=0$ . Replacing $R$ by $\mathcal{K}_{n}$ in $\left(
\ref{connectsum}\right)  $ we have
\[
r_{n}\left(  X\right)
=\operatorname{rk}_{\mathcal{K}_{n}}H_{1}\left( X_{1};\mathcal{K}
_{n}\right) +\operatorname{rk}_{\mathcal{K}_{n}}H_{1}\left(
X_{2};\mathcal{K}_{n}\right) +1.
\]
We will show that $\operatorname{rk}_{\mathcal{K}_{n}}H_{1}\left(
X_{1};\mathcal{K}_{n}\right)
=\operatorname{rk}_{\mathcal{K}_{n}^{j}}H_{1}\left(
X_{1};\mathcal{K}_{n} ^{j}\right)  $ hence
\[
r_{n}\left(  X\right)  =r_{n}\left(  X_{1}\right)  +r_{n}\left(
X_{2}\right) +1.
\]

Let $\widetilde{X}_{jn}$ be the cover of $X_{j}$ corresponding to
$G_{j}\rightarrow\Gamma_{n}$. Then $X_{jn}^{\prime}$ is a disjoint
union of $\Gamma_{n}/\Gamma_{n}^{j}$ copies of $X_{jn}$. The
extension of $\overline{\operatorname{pr}}_{j}$ to a ring
homomorphism $\overline{\operatorname{pr}}_{j}:\mathbb{Z}
\Gamma_{n}\rightarrow\mathbb{Z}\Gamma_{n}^{j}$ gives
$\mathbb{Z}\Gamma_{n} ^{j}$ the structure as a
$\mathbb{Z}\Gamma_{n}$-bimodule. Moreover, since
$\overline{\operatorname{pr}}_{j}\circ\overline{i}_{j}=id_{\Gamma_{n}^{j}}$,
$\mathbb{\cdot }\otimes_{\mathbb{Z}\Gamma_{n}^{j}}\left(
\mathbb{Z}\Gamma_{n}\otimes
_{\mathbb{Z}\Gamma_{n}}\mathbb{Z}\Gamma_{n}^{j}\right)  $ acts
trivially on any right $\mathbb{Z}\Gamma_{n}^{j}$-module.
Therefore
\begin{align}
C_{\ast}\left(  \widetilde{X}_{jn}\right)
\otimes_{\mathbb{Z}\Gamma_{n} }\mathbb{Z}\Gamma_{n}^{j}  &
\cong\left(  C_{\ast}\left(  X_{jn}\right)
\otimes_{\mathbb{Z}\Gamma_{n}^{j}}\mathbb{Z}\Gamma_{n}\right)
\otimes
_{\mathbb{Z}\Gamma_{n}}\mathbb{Z}\Gamma_{n}^{j} \label{changecc} \\
& \cong C_{\ast}\left(  X_{jn}\right)
\otimes_{\mathbb{Z}\Gamma_{n}^{j} }\left(
\mathbb{Z}\Gamma_{n}\otimes_{\mathbb{Z}\Gamma_{n}}\mathbb{Z}
\Gamma_{n}^{j}\right) \nonumber \\
& \cong C_{\ast}\left(  X_{jn}\right). \nonumber
\end{align}
$\overline{i}_{j}:\mathbb{Z}\Gamma_{n}^{j}\rightarrow\mathbb{Z}
\Gamma_{n}$ is a monomorphism, hence we can extend
$\overline{i}_{j}$ to the right ring of quotients of
$\mathbb{Z}\Gamma_{n}^{j}$ and $\mathbb{Z}\Gamma_{n}$,
$\overline{i}_{j}:\mathcal{K}_{n}^{j}\rightarrow\mathcal{K}_{n}$.
Therefore $\mathcal{K}_{n}$ is a flat left
$\mathcal{K}_{n}^{j}$-module by the following lemma.

\begin{lemma}\label{flatlemma} Suppose that $R$ is a right
and left principal ideal domain, $S$ has no zero divisors, and
$f:R \hookrightarrow S$ is a ring monomorphism. Then $S$ is a flat
left $R$-module.
\end{lemma}
\begin{proof}
Let $s \in S$ and $r \in R$ with $s \neq 0$.  Suppose that
$f\left(r\right)s=0$. $S$ has no zero-divisors hence
$f\left(r\right)=0$.  Moreover, $f$ is a monomorphism so $r=0$.
Therefore $S$ is $R$-torsion-free.  Since $R$ is a PID, every
finitely generated torsion-free $R$-module is free hence flat.
Every module is the direct limit of its finitely generated
submodules.   Hence $S$ is the direct limit of flat modules. Thus,
by \cite[Prop. 10.3]{Ste}, $S$ is flat.
\end{proof}

We apply $-\otimes_{\mathbb{Z}\Gamma_{n}^{j}}\mathcal{K}_{n}$ to
$\left( \ref{changecc}\right)  $ to get
\[
C_{\ast}\left(  \widetilde{X}_{jn}\right)
\otimes_{\mathbb{Z}\Gamma_{n} }\mathcal{K}_{n}\cong C_{\ast}\left(
X_{jn}\right)  \otimes_{\mathbb{Z} \Gamma_{n}^{j}}\mathcal{K}_{n}.
\]
Since $C_{\ast}\left(  X_{jn}\right)
\otimes_{\mathbb{Z}\Gamma_{n}^{j} }\mathcal{K}_{n}\cong
C_{\ast}\left(  X_{jn}\right)  \otimes_{\mathbb{Z}
\Gamma_{n}^{j}}\mathcal{K}_{n}^{j}\otimes_{\mathcal{K}_{n}^{j}}\mathcal{K}
_{n}$ and $\mathcal{K}_{n}$ is a flat left
$\mathcal{K}_{n}^{j}$-module,
\[
H_{\ast}\left(  C_{\ast}\left(  \widetilde{X}_{jn}\right)  \otimes
_{\mathbb{Z}\Gamma_{n}}\mathcal{K}_{n}\right)  \cong
H_{\ast}\left(  C_{\ast }\left(  X_{jn}\right)
\otimes_{\mathbb{Z}\Gamma_{n}^{j}}\mathcal{K}_{n} ^{j}\right)
\otimes_{\mathcal{K}_{n}^{j}}\mathcal{K}_{n}.
\]
Thus $\operatorname{rk}_{\mathcal{K}_{n}^{j}}H_{1}\left(
X_{1};\mathcal{K}_{n}^{j}\right)
=\operatorname{rk}_{\mathcal{K}_{n}}H_{1}\left(
X_{1};\mathcal{K}_{n}\right)  $ as desired.

Now we show that
\begin{equation}
T_{ R}H_{1}\left(  X; R\right)  \cong T_{ R} H_{1}\left( X_{1};
R\right)  \oplus T_{ R}H_{1}\left( X_{2}; R\right)
.\label{addtorsion}
\end{equation}
First we note that $T\left(  H_{1}\left(  X_{1}; R\right) \oplus
H_{1}\left(  X_{2}; R\right)  \right)  \cong TH_{1}\left( X_{1};
R\right)  \oplus TH_{1}\left( X_{2}; R\right)  $.
Consider the restriction of $\nu$ in $\left(
\ref{connectsum}\right)$ to the torsion submodule of $H_{1}\left(
X_{1}; R\right)  \oplus H_{1}\left( X_{2}; R\right)  $,
\[
\nu_{T}:T\left(  H_{1}\left(  X_{1}; R\right)  \oplus H_{1}\left(
X_{2}; R\right)  \right)  \rightarrow TH_{1}\left( X; R \right)  .
\]
We show that $\nu_{T}$ is an isomorphism. It is immediate that
$\nu_{T}$ is a monomorphism since $\nu$ is a monomorphism. To show
that $\nu_{T}$ is surjective, let $x\in TH_{1}\left( X; R\right) $
and $0\neq r\in R$ with $xr=0$. Since $H_{0}\left( S^{2}; R\right)
$ is $ R$-torsion free, $\partial_{1}\left(  x\right) =0$ hence
there exists $y\in H_{1}\left(  X_{1}; R\right) \oplus H_{1}\left(
X_{2}; R\right)  $ such that $\nu_{T}\left( y\right) =x$.
Moreover, $\nu_{T}\left(  yr\right) =\nu_{T}\left( y\right)
r=xr=0$. Since $\nu_{T}$ is a monomorphism $yr=0$. Hence $y\in
T\left( H_{1}\left( X_{1}; R \right) \oplus H_{1}\left( X_{2};
R\right) \right)  $.

Since $H^{1}\left( X ; \mathbb{Z}\right) \cong H^{1}\left(
X_{1};\mathbb{Z}\right) \oplus H^{1}\left( X_{2};\mathbb{Z}\right)
$, $\psi $ can be uniquely written as $\psi _{1}\oplus\psi _{2}$
where $\psi _{j}\in H^{1}\left( X_{j}; \mathbb{Z} \right) $. Note
that $\psi _{j}$ need not be a primitive class in $H^{1}\left(
X\right) $. For each $j$, let $d_{j}$ be the largest divisor of
$\psi _{j}$. Hence, there exist $\psi _{j}^{\prime }$ primitive
with $d_{j}\psi _{j}^{\prime }=\psi $.    Recall that $\ker
\psi_j^{\prime} = \ker \psi_j$ and $\delta _{n}\left( \psi
_{j}\right) =d_{j}\delta _{n}\left( \psi _{j}^{\prime }\right)$.

Substitute $ R=\mathbb{Z}\Gamma_{n}\left( \mathbb{Z}\ker
\psi\right)^{-1}$ into $\left( \ref{addtorsion}\right) $.  Then
\[
\delta_{n}\left(  \psi\right)
=\operatorname{rk}_{\mathbb{K}_{n}}TH_{1}\left(
X_{1};\mathbb{Z}\Gamma_{n}\left(  \mathbb{Z}\ker\psi\right)
^{-1}\right) \oplus \operatorname{rk}_{\mathbb{K}_{n}}TH_{1}\left(
X_{2};\mathbb{Z}\Gamma_{n}\left( \mathbb{Z}\ker\psi\right)
^{-1}\right)
\]
where $\mathbb{K}_{n}=\mathbb{Z}\ker\psi\left( \mathbb{Z}\ker\psi
- \{ 0 \} \right) ^{-1}$ is the right ring of quotients of
$\mathbb{Z}\ker\psi$. Recall that if $S$ is a right divisor set
then $RS^{-1}$ exists and is the ring obtained by inverting all of
the elements in $S$. Let $R_{n}=\mathbb{Z}\Gamma_{n}\left(
\mathbb{Z}\ker\psi - \{0\} \right) ^{-1}$,
$R_{n}^{j}=\mathbb{Z}\Gamma_{n}^{j}\left( \mathbb{Z}\ker\psi
_{j}^{\prime} - \{0\} \right) ^{-1}$, and
$\mathbb{K}_{n}^{j}=\mathbb{Z}\ker\psi _{j}^{\prime} \left(
\mathbb{Z} \ker \psi _{j}^{\prime} - \{0\} \right) ^{-1}$. To
complete the proof we must show that
\[
\operatorname{rk}_{\mathbb{K}_{n}}T_{R_{n}}H_{1}\left(
X_{j};R_{n}\right) \cong
\operatorname{rk}_{\mathbb{K}_{n}^{j}}T_{R_{n}^{j}}H_{1}\left(
X_{1};R_{n}^{j}\right) .
\]

Since $\psi=\psi_{j}\circ\overline{i}_{j}$,
$\overline{i}_{j}\left(  \ker \psi_{j}\right)  \subset\ker\psi$,
we can extend $\overline{i}_{j}$ to
$\overline{i}_{j}:R_{n}^{j}\rightarrow R_{n}$. By Lemma
\ref{flatlemma}, $R_{n}$ is a flat left $R_{n}^{j}$-module.
Therefore
\begin{equation}
H_{\ast}\left(  C_{\ast}\left(  \widetilde{X}_{jn}\right)  \otimes
_{\mathbb{Z}\Gamma_{n}}R_{n}\right)  \cong H_{\ast}\left(
C_{\ast}\left( X_{jn}\right)
\otimes_{\mathbb{Z}\Gamma_{n}^{j}}R_{n}^{j}\right)
\otimes_{^{R_{n}^{j}}}R_{n}.\label{flatR}
\end{equation}

Let $M=H_{1}\left(  X_{j};R_{n}^{j}\right)$ then $M \cong
\left(R_{n}^{j}\right)^{m}\oplus T_{R_{n}^{j}}M$ since $M$ is
finitely generated and $R_{n}^{j}$ is a principal ideal domain. It
is straightforward to show that
\[
T_{R_{n}}\left(  M\otimes_{R_{n}^{j}}R_{n}\right)  \cong
T_{R_{n}^{j}} M\otimes_{R_{n}^{j}}R_{n}.
\]
Moreover, $T_{R_{n}^j} M \cong \frac{R_{n}^j}{\langle r_1 \rangle}
\oplus \cdots \oplus \frac{R_{n}^j}{\langle r_k \rangle}$  so it
suffices to show that $\operatorname{rk}_{\mathbb{K}_n}
\frac{R_n}{\langle \overline{i}_j \left(r\right) \rangle} = d_j
 \operatorname{rk}_{\mathbb{K}_n^j} \frac{R_n^j}{\langle r \rangle}$
for any nonzero $r \in R_n^j$.  Note that this would imply
\begin{align*}
\operatorname{rk}_{\mathbb{K}_n} \left(T_{R_n^j} M \otimes_{R_n^j}
R_n \right) &= \operatorname{rk}_{\mathbb{K}_n} \left(
\frac{R_n}{\langle \overline{i}_j\left(r_1\right) \rangle} \oplus
\cdots \oplus \frac{R_n}{\langle \overline{i}_j
\left(r_k\right) \rangle} \right) \\
&= \operatorname{rk}_{\mathbb{K}_n} \frac{R_n}{\langle
\overline{i}_j\left(r_1\right) \rangle} + \cdots
+\operatorname{rk}_{\mathbb{K}_n} \frac{R_n}{\langle
\overline{i}_j
\left(r_k\right) \rangle} \\
&= d_j \operatorname{rk}_{\mathbb{K}_n^j} \frac{R_n^j}{\langle r_1
\rangle} + \cdots + d_j \operatorname{rk}_{\mathbb{K}_n^j}
\frac{R_n^j}{\langle r_k \rangle}.
\end{align*}

Let $T \in \Gamma_n^j$ such that $\psi_j\left(T\right)=t^{d_j}$.
We can write any element in $\Gamma_n^j$ as $T^m \kappa$ where
$\kappa \in \ker \psi_j$. Hence $r$ can be written as a
nonconstant (Laurent) polynomial in T with coefficients in
$\mathbb{K}_j$.  We can assume that $r = a_0 + T a_1 + ... + T^q
a_q$ with $a_0 \neq 0$.

Since $\psi$ is surjective, there is an $S \in \Gamma_n$ such that
$\psi\left(S\right) = t$.  We can write any element in $\Gamma_n$
as $S^p f$ where $f \in \ker \psi$.  In particular, any element of
$\Gamma_n$ that maps $t^{d_j}$ under $\psi$ can be written as
$S^d_j f$. Since $\psi\left( \overline{i}_j \left(T\right)
\right)=\psi\left(T\right)=t^d_j$, $\overline{i}_j\left(T\right)=
S^d_j f$ for some $f \in \ker \Gamma_n$. Hence
\begin{align*}
\overline{i}_j\left(r\right) &= \overline{i}_j\left(a_0\right) +
\overline{i}_j\left(T\right) \overline{i}_j\left(a_1\right) +
\cdots + \overline{i}_j\left(T\right)^q
\overline{i}_j\left(a_q\right) \\
   &  = \overline{i}_j\left(a_0\right) + S^d_j f \overline{i}_j\left(a_1\right) + \cdots + (S^d_j
   f)^q \overline{i}_j\left(a_q\right)  \\
   &  = \overline{i}_j\left(a_0\right) + S^d_j g_1 \overline{i}_j\left(a_1\right) + \cdots + S^{d_j q} g_k \overline{i}_j\left(a_q\right)
\end{align*}
for some $g_i \in \ker \Gamma_n$.  We note that
$\overline{i}_j\left(a_i\right) \in \ker \Gamma_n$ which gives us
our desired result. This completes the proof that
$\delta_{n}\left(  \psi\right)  =\delta_{n}\left(  \psi_{1}\right)
+\delta _{n}\left(  \psi_{2}\right)$.
\end{proof}

An immediate consequence is that $r_n\left(X\right) \geq 1$
whenever the hypotheses in Proposition \ref{connectedsum} are
satisfied. In particular, we have $\overline{\delta}_{n}\left(
\psi\right)=0$ for all $\psi \in H^1\left(X;\mathbb{Z}\right)$. We
note that if $\beta_1 \left( X_j \right) =0$ for some $j$ then
$\delta_n \left( \psi\right) =\delta_{n}\left( \psi\right) $ and
$r_n \left( X\right) =r_{n}\left( X_{2}\right) $ (similarly if
$\beta _{1}\left( X_{2}\right) =0$).

\begin{corollary}
Let $X$ be a compact, orientable 3-manifold with
$r_k\left(X\right)=0$ for some $k \geq 0$.  Suppose that
$G=\pi_1\left(X\right)$ does not satisfy both
$\frac{G}{G_r^{\left(1\right)}} \cong \mathbb{Z}$ and $\frac
{G_r^{\left(1\right)} } {G_r^{\left(2\right)}}=0$.  Then there
exists an irreducible 3-manifold $Y$ with $H=\pi_1\left(Y\right)$
such that
\[
\frac{G}{G_r^{\left(n+1\right)}}=\frac{H}{H_r^{\left(n+1\right)}}
\]
for all $n \geq 0$.
\end{corollary}
\begin{proof}
We assume that $\beta_1\left(X\right) \geq 1$.  We can factor $X$
as  $X=X_1 \# \cdots \# X_l$ where each $X_i$ is prime \cite{H}.
Since $r_k\left(X\right)=0$, there is exactly one $i$ such that
$\beta_1\left(X_i\right) \neq 0$ by Proposition
\ref{connectedsum}.  Let $Y$ be the aforementioned factor and
$H=\pi_1\left(Y\right)$.  It is easy to verify that
$\frac{G}{G_r^{\left(n+1\right)}}=\frac{H}{H_r^{\left(n+1\right)}}$.
Moreover, the hypothesis on $G$ guarantees that $Y \neq S^2 \times
S^1$.  Therefore $Y$ is irreducible.
\end{proof}

\section{Rank of Torsion Modules over Skew Polynomial
Rings}\label{sec:ranktorsionmodules}

In this section we will show that the rank of a (torsion) module
presented by an $m \times m$ matrix of the form $A+tB$ (where
$A,B$ have  coefficients in $\mathbb{K}$) has rank at most m as a
$\mathbb{K}$ vector space. This is well known when $\mathbb{K}$ is
a commutative field.  In this case the rank of $M$ over
$\mathbb{K}$ is the determinant of $A+tB$ which is a polynomial
with degree less than or equal to $l$.  We will use this result in
the proof that the higher-order degrees give lower bounds for the
Thurston norm in the next section.  For a first read, the reader
may wish to only read the statements in Proposition
\ref{boundrank} and Lemma \ref{pres} before proceeding to the next
section.

Let $M$ be a right $\mathbb{K}\left[  t^{\pm1}\right]  $-module
with presentation matrix of the form $A+tB$ where $A$ and $B$ are
$l\times m$ matrices with coefficients in $\mathbb{K}$, $l$ is the
number of generators of $M$ and $\mathbb{K}$ is a (skew) field. We
denote by $TM$ the $\mathbb{K} \left[  t^{\pm1}\right]  $-torsion
submodule of $M$. Using the embedding
$\mathbb{K}\rightarrow\mathbb{K}\left[  t^{\pm1}\right]  $ by
$k\mapsto k\cdot1$ we consider $TM$ as a module over $\mathbb{K}$.

\begin{proposition}
\label{boundrank}If $M$ is a right $\mathbb{K}\left[
t^{\pm1}\right] $-module with presentation matrix of the form
$A+tB$ where $A$ and $B$ are $l\times m$ matrices with
coefficients in $\mathbb{K}$ and $\mathbb{K}$ is a division ring
then
\begin{equation}
\operatorname{rk}_{\mathbb{K}}TM\leq\min\left\{  l,m\right\} .
\end{equation}
\end{proposition}

We begin by stating when two presentation matrices of a finitely
presented right $ R$-module $H$ are equivalent.

\begin{lemma}[pp. 117-120 of \cite{zass}]
\label{pres}Two presentation matrices of H are related by a finite
sequence of the following operations.

\begin{enumerate}
\item  Interchange two rows or two columns.

\item  Multiply a row (on left) or column (on right) by a unit of
$ R$.

\item  Add to any row a $ R$-linear combination of other rows
(multiplying a row by unit of $ R$ on left) or to any column a $
R$-linear combination of other columns (multiplying a column by a
unit of $ R$ on right). \item $P\rightarrow\left(
\begin{array}
[c]{cc} P & \ast
\end{array}
\right)  $, where $\ast$ is a $ R$-linear combination of columns
of $P$.

\item $P\rightarrow\left(
\begin{array}
[c]{cc}
P & \ast\\
0 & 1
\end{array}
\right)  $, where $\ast$ is an arbitrary column.
\end{enumerate}
\end{lemma}

We will find the following lemmas useful in the proof of
Proposition \ref{boundrank}.

\begin{lemma}
\label{zeroblock}A presentation matrix of the form$\left(
\begin{array}
[c]{cc}
A_{1}+tI_{s} & A_{2}\\
A_{3} & A_{4}
\end{array}
\right)  _{l\times m}$ where $A_{i}$ has entries in $\mathbb{K}$
(a non-commutative division ring) is related (in the sense of
Lemma \ref{pres}) to a matrix of the form $\left(
\begin{array}
[c]{cc}
A_{1}^{\prime}+tI_{s} & A_{2}^{\prime}\\
A_{3}^{\prime} & 0
\end{array}
\right)  _{\left(  l-r\right)  \times\left(  m-r\right)  }$ for
some $r\geq0$.
\end{lemma}

\begin{proof}
Since $A_{4}$ is a matrix over $\mathbb{K}$ \cite[Corollary to
Theorem 16, p.43]{Ja} there are $C$ and $D$ such that
$CA_{4}D=\left(
\begin{array}
[c]{cc}
0 & 0\\
0 & I_{r}
\end{array}
\right)  $. Here, $C$ and $D$ are units in the rings of $l\times
l$ and $m\times m$ matrices with entries in $\mathbb{K}$
respectively. Hence we can get the new presentation matrix \
\begin{align*}
\left(
\begin{array}
[c]{cc}
I & 0\\
0 & C
\end{array}
\right)  \left(
\begin{array}
[c]{cc}
A_{1}+tI_{s} & A_{2}\\
A_{3} & A_{4}
\end{array}
\right)  \left(
\begin{array}
[c]{cc}
0 & I\\
D & 0
\end{array}
\right)  &=\left(
\begin{array}
[c]{cc}
A_{1}+tI_{s} & A_{2}D\\
CA_{3} & CA_{4}D
\end{array}
\right) \\
&=\left(
\begin{array}
[c]{cc}
A_{1}+tI_{s} & A_{2}D\\
CA_{3} &
\begin{array}
[c]{cc}
0 & 0\\
0 & I_{r}
\end{array}
\end{array}
\right)  .
\end{align*}
Now we can make the last $r$ rows of the matrix of the form
$\left(
\begin{array}
[c]{cc} 0 & I_{r}
\end{array}
\right)  $ by adding $\left(  \text{column }\left(  m-\left(
l-i\right) \right)  \right)  \cdot\left(  -a_{i,j}\right)  $ to
column $j$ for each nonzero entry $a_{i,j}$ in the last $r$ rows
of $CA_{3}$. In general this will change $A_{1}+tI_{s}$ to
$A_{1}^{\prime}+tI_{s}$ for some $A_{1}^{\prime }$ whose entries
lie in $\mathbb{K}$. Using operation $5$, we delete the last $r$
rows and columns to obtain our desired result.
\end{proof}

\begin{lemma}
\label{reducesize}If $A_{3}\neq0$ then the presentation matrix
$\left(
\begin{array}
[c]{cc}
A_{1}+tI_{s} & A_{2}\\
A_{3} & 0
\end{array}
\right)  $ of size $l\times m$ is related to one of the form
$\left(
\begin{array}
[c]{cc}
A_{1}^{\prime}+tI_{s-1} & A_{2}^{\prime}\\
A_{3}^{\prime} & 0
\end{array}
\right)  $ of size $\left(  l-r\right)  \times\left(  m-r\right)
$ where $r\geq1$.
\end{lemma}

\begin{proof}
Let $A=\left(
\begin{array}
[c]{cc}
A_{1} & A_{2}\\
A_{3} & 0
\end{array}
\right)  $ and $a_{k,i}$ be the $\left(  k,i\right)  $ entry of
$A$.\ By permuting rows in $A_{3}$ we can assume that the last row
has a non-zero element. Suppose that the first non-zero element in
this row occurs in the $i^{th}$ column. We can assume that this
element is $1$. Now if $a_{l,j}$ is any other non-zero entry in
the last row $\left(  i<j\leq s\right)  $ we add $\left(
\text{column }i\right)  \cdot\left(  -a_{l,j}\right)  $ to $\left(
\text{column }j\right)  $ to get a presentation with a zero in
column $j$ of the last row.\ However, this changes the $\left(
i,j\right)  $ the entry of our matrix to $\left(
a_{i,j}-a_{i,i}a_{l,j}\right) -ta_{l,j}$ which does not lie in
$\mathbb{K}$. To remedy this, we add $\left(  ta_{l,j}
t^{-1}\right)  \cdot\left(  \text{row }j\right)  $ to $\left(
\text{row }i\right)  $. Performing these two steps for all
non-zero $a_{l,j}$ gives us a matrix whose last row is of the form
$\left(  0,\ldots,0,1,0,\ldots ,0\right)  $. By cyclically
permuting columns $i$ through $m$ (so that the $i^{th}$ column
becomes the $m^{th}$ column) and using the operation of type $5$
in Lemma \ref{pres}, we see that this matrix is related to the
matrix obtained by deleting column $i$ and row $l$. We note that
all the entries in row $i$ lie in $\mathbb{K}$. For the final step
we cyclically permute rows $i$ through $s$ (so that the $i^{th}$
row becomes the $s^{th}$ row) and use Lemma \ref{zeroblock} to get
our desired result.
\end{proof}

\begin{proof}[Proof of Proposition \ref{boundrank}]
Let $P=A+tB$ be a presentation matrix of $M$. \ As in the proof of
lemma \ref{zeroblock} there are $C_{l\times l}$ and $D_{m\times
m}$ (units in the rings of matrices over $\mathbb{K}$) such that
$CBD=\left(
\begin{array}
[c]{cc}
I_{s} & 0\\
0 & 0
\end{array}
\right)  $, where $s\leq\min\left\{  l,m\right\}  $. Now if
$C^{t}\equiv tCt^{-1}$ we have
\[
C^{t}PD=C^{t}\left(  A+tB\right)  D=C^{t}AD+tCBD.
\]
Hence $M$ has a presentation matrix of the form
\begin{equation}
A+t\left(
\begin{array}
[c]{cc}
I_{s} & 0\\
0 & 0
\end{array}
\right)  _{l\times m}= \left(
\begin{tabular}
[c]{cccc}
$a_{1,1}+t$ & $\cdots$ & $a_{1,s}$ & \\
$\vdots$ & $\ddots$ & $\vdots$ & $\ast$\\
$a_{s,1}$ & $\cdots$ & $a_{s,s}+t$ & \\
& $\ast$ &  & $\ast$
\end{tabular}
\label{step1} \right) \end{equation}
 where by $A\ $is now
$C^{t}AD$, a matrix with entries in $\mathbb{K}$. We can now use
Lemma \ref{zeroblock} and Lemma \ref{reducesize} to get a new
presentation matrix $\left(
\begin{array}
[c]{cc}
A_{1}^{\prime}+tI_{s^{\prime}} & A_{2}^{\prime}\\
0 & 0
\end{array}
\right)  $ where $s^{\prime}\leq s$. It follows that $TM$ has a
presentation matrix
\begin{equation}
\left(
\begin{array}
[c]{cc} A_{1}^{\prime}+tI_{s^{\prime}} & A_{2}^{\prime}
\end{array}
\right) \label{torpres}
\end{equation}
where $A_{i}^{\prime}$ are matrices in $\mathbb{K}$. Let
$\sigma_{i}$ $\left(  1\leq i\leq s^{\prime}\right)  $ be the
generators of $TM$ corresponding to $\left(  \ref{torpres}\right)
$. We show that these generate $TM$ as a $\mathbb{K}$-module. Let
$a_{i,j}$ be the $\left( i,j\right)  $ the entry of
$A_{1}^{\prime}$ then we have the relations
$\sigma_{1}a_{1,j}+\cdots+\sigma_{j}\left( a_{j,j}+t\right)
+\cdots +\sigma_{s^{\prime}}a_{s^{\prime},j}=0$ for $j\leq
s^{\prime}$. Hence $\sigma_{i}t=- \sum \sigma_{k}a_{k,i}$ is in
the span of $\left\{  \sigma_{i}\right\} $. We prove by induction
on $n$ that $\sigma_{i}t^{n}$ is in the span of $\left\{
\sigma_{i}\right\}  $. Suppose $\sigma_{i}t^{n}= \sum
\sigma_{k}b_{k,i}$ where $b_{k,i}\in\mathbb{K}$ then
\begin{align*}
\sigma_{i}t^{n+1}&=\left(  \sum_k
\sigma_{k}b_{k,i}\right)  t\\
&= \sum_k
\sigma_{k}t\left(  t^{-1}b_{k,i}t\right) \\
&= \sum_k \left(  \sum_l \sigma_{l}b_{l,k}\right)  \left(  t^{-1}b_{k,i}t\right) \\
&= \sum_l \sigma_{l}\left(  \sum_k b_{l,k}\left(
t^{-1}b_{k,i}t\right) \right)
\end{align*}
for all $i\leq s^{\prime}$. Therefore any element $ \sum
\sigma_{i}p_{i}\left(  t\right)  $ with $p_{i}\left(  t\right)  $
$\in\mathbb{K}\left[  t^{\pm1}\right]  $ can be written as a
linear combination of $\sigma_{i}$ with coefficients in
$\mathbb{K}$. It follows that
$\operatorname{rk}_{\mathbb{K}}TM\leq s^{\prime}\leq
s\leq\min\left\{  l,m\right\}  $.
\end{proof}

\section{Relationships of $\delta_{n}$ and $\overline{\delta}_{n}$ to the
Thurston Norm}\label{sec:reltothurstonnorm}

In this section, we will prove one of the main theorems of this
paper. We show that the higher-order degrees of a 3-manifold give
lower bounds for the Thurston norm.  The result when $X$ is a knot
complement appears in \cite{Co} although it uses some of our work.

\begin{theorem}
\label{delbarthm}Let $X$ be a compact, orientable 3-manifold
(possibly with boundary). For all $\psi\in H^{1}\left(
X;\mathbb{Z}\right)  $ and $n\geq0$
\[
\overline{\delta}_{n}\left(  \psi\right)  \leq\left\| \psi\right\|
_T
\]
except for the case when $\beta_{1}\left(  X\right) =1$, $n=0$,
$X\ncong S^1 \times S^2$, and $X \ncong S^1 \times D^2$. In this
case, $\overline{\delta}_{0}\left( \psi\right) \leq\left\|
\psi\right\| _{T}+1+\beta_{3}\left( X\right)  $ whenever $\psi$ is
a generator of $H^{1}\left( X;\mathbb{Z}\right) \cong\mathbb{Z}$.
Moreover, equality holds in all cases when $\psi:\pi_{1}\left(
X\right) \twoheadrightarrow \mathbb{Z}$ can be represented by a
fibration $X\rightarrow S^{1}$.
\end{theorem}

The proof of this Theorem will follow almost directly from
Proposition \ref{alexpres} and Proposition \ref{boundrank}.
However, because of some technical details we postpone the proof
until after Corollary \ref{cor:delbar0}.  We will begin the
section by proving a more generalized (but less applicable)
version of Theorem \ref{delbarthm}.  We first introduce some
notation.

Let $X$ be a 3-manifold, $\psi\in H^{1}\left(  X;\mathbb{Z}\right)
$, $G=\pi_{1}\left(  X\right)  $, $\Gamma_{n}=G/G_{r}^{\left(
n+1\right)  }$.  Recall that if $F$ is an embedded surface dual to
$\psi$, we can consider the homology of $F $ with coefficients in
$\mathbb{K}_{n}$, where $\mathbb{K} _{n}$ is the field of
fractions of $\mathbb{Z}\Gamma_{n}^{\prime}$.  Define the
higher-order Betti numbers of $F$ to be
\[
b_{i}^{n}\left(  F\right) =\operatorname{rk}_{\mathbb{K}_{n}}
H_i\left( F;\mathbb{K}_{n}\right)  .
\]
By Remark \ref{eulerchar} we see that the Euler characteristic of
$F$ can be computed using $b_i^n$,
\begin{equation}\label{equ:addchi}
\chi\left(  F\right)  =\sum\left(  -1\right)  ^{i}b_{i}^{n}\left(
F\right)
\end{equation}
for any $n\geq0$.

Now we consider the collection of Thurston norm minimizing
surfaces dual to $\psi$, $\mathcal{F}_{\psi}$.  It is very
possible that a surface in $\mathcal{F}_{\psi}$ is highly
disconnected.  One could ask, ``What is the minimal number of
components of a surface in $\mathcal{F}_{\psi}$?''  For our
purposes, it will turn out to be important to compute the number
of components of surface in $\mathcal{F}_{\psi}$ that lift to the
$n^{th}$ order cover of $X$.  To be precise we make the following
definitions.

Let $F=\coprod F^{i}$ be a (possibly disconnected) surface.  We
define $N_n\left( F\right) $ to be the number of components of $F$
with $i_{\ast }\left( \pi_{1}\left( F^{i}\right) \right) \subseteq
G_{r}^{\left( n+1\right)  }$ and $N_n^{c}\left( F\right)  $ to be
the number of closed components of $F$ with $i_{\ast}\left(
\pi_{1}\left( F^{i}\right) \right) \subseteq G_{r}^{\left(
n+1\right)  }$. Finally, we define
\[\mathcal{N}_n\left(\psi\right)=\min_{F \in \mathcal{F}_{\psi} }\left\{ N_n\left( F\right)+N_n^{c}\left( F\right)
  \right\}.
\]

\begin{theorem}\label{lowerbound}
Let $X$ be a compact, orientable 3-manifold (possibly with
boundary).  For all $\psi\in H^{1}\left( X;\mathbb{Z}\right) $ and
$n\geq0$
$$\delta_{n}\left( \psi\right) \leq\left\| \psi\right\| _{T}+
\mathcal{N}_n\left(\psi\right).$$
\end{theorem}

\begin{proof}
Let $F$ be a Thurston norm minimizing surface dual to $\psi$ that
minimizes $N_n\left( F\right)+N_n^{c}\left( F\right)$. We remark
that a connected surface has $b_0^n\left(F\right)=0$ if and only
if the coefficient system $\phi \circ i_{\ast}:\pi_1\left(F\right)
\rightarrow G/G_r^{\left(n\right)}$ is trivial by Lemma
\ref{lem:H0}. Therefore $N_{n}\left( F\right)= b_0^n\left(
F\right)$. Similarly, we have $N_n^c\left( F\right)= b_2^n\left(
F\right)$. By $\left(\ref{equ:addchi}\right)$,
\begin{align*}
b_1^n\left(  F\right) &=-\chi\left(  F\right) +N_{n}\left(
F\right) +N_{n}^{c}\left(  F\right) \\&\leq \left\|\psi\right\|_T
+ \mathcal{N}_n\left(\psi\right).
\end{align*}

To complete the proof, we show that $\delta_n\left(\psi\right)
\leq b_1^n\left(F\right)$.  By Proposition \ref{pro:pres_A},
$\mathcal{A}_n^{\psi}\left(X\right)$ has a presentation matrix of
the form $A+tB$ of size $\left( b_1^n\left( F\right)  \times
m\right) $ where $m=rk_{\mathbb{K}_n}
H_1\left(Y;\mathbb{K}_n\right)$.  Thus, by Proposition
\ref{boundrank} we have
\[
\delta_{n}\left(  \psi\right)
=\operatorname{rk}_{\mathbb{K}_{n}}\mathcal{A}_{n}^{\psi }\left(
X\right) \leq\min\left\{  b_{1}^{n}\left(  F\right)  ,m\right\}
\leq b_{1}^{n}\left(  F\right).
\]
\end{proof}

We note that the term $\mathcal{N}_n\left(\psi\right)$ is an
invariant of the pair $\left(X,\psi\right)$.  However, in a
general, $\mathcal{N}_n\left(\psi\right)$ may be difficult to
compute.  Fortunately, in some cases, we may be able bound this
term by a constant.

Suppose that we are interested in the genera of knots or links.
More generally, suppose we are only interested in the connected
surfaces embedded in a 3-manifold.  Then it is reasonable to
measure the complexity of the surface by its first Betti number.
Using the proof of Theorem \ref{lowerbound}, we can find a lower
bound for the first Betti number of $F$ that has no ``extra
term.''
\begin{corollary}
If $F$ is any surface dual to $\psi$ then $\delta_{n}\left(  \psi\right)
\leq\beta_{1}\left(  F\right)  $.
\end{corollary}

\begin{proof}
Since $N_{n}\left(  F\right)  \leq\beta_{0}\left(  F\right)  $ and
$N_{n}^{c}\left(  F\right)  \leq\beta_{2}\left(  F\right)  $,
\begin{align*}
b_{1}^{n}\left(  F\right)  &=-\chi\left(  F\right)  +N_{n}\left(
F\right)
+N_{n}^{c}\left(  F\right) \\
&=\beta_{1}\left(  X\right)  +\left(  N_{n}\left(  F\right)
-\beta_{0}\left( F\right)  \right)  +\left(  N_{n}^{c}\left(
F\right)  -\beta_{2}\left(
F\right)  \right) \\
&\leq\beta_{1}\left(  X\right)  .
\end{align*}
Therefore, $\delta_{n}\left(  \psi\right)  \leq b_{1}^{n}\left(
F\right)  \leq\beta_{1}\left(  X\right)  $.
\end{proof}

We consider the case when $X$ is the complement of a link $L$ in
$S^3$. If $L$ has $m$ components then
$H_1\left(X;\mathbb{Z}\right) \cong \mathbb{Z}^m$ generated by the
$m$ meridians $\mu_i$.  Let $\psi_i$ be defined by
$\psi_i\left(\mu_i\right)=t^{\delta_{ij}}$. That is, $\psi_i$ is
dual to any surface that algebraically intersects the $i^{th}$
meridian once and the $j^{th}$ meridian zero times for $j \neq i$.
We will show that a Thurston norm minimizing surface dual to
$\psi_i$ can be chosen to be connected and hence we can bound the
term $\mathcal{N}_n\left(\psi_i\right)$ by 1.

\begin{corollary}\label{cor:linkbound}Let $X=S^3-L$ and $\psi_i$
be as defined above.  Then
\[\delta_n\left(  \psi_i\right)
\leq\left\|  \psi_i\right\|_T + 1
\]
for all $n \geq 0$.
\end{corollary}
\begin{proof}We show that for all $i \in \{1,\dots, m\}$ there exists a Thurston minimizing surface
$F_i$ which is connected and has non-trivial boundary.  Hence
$\mathcal{N}_n\left(\psi_i\right) \leq 1$ and
$\mathcal{N}_n^c\left(\psi_i\right)=0$. The result follows from
Theorem \ref{lowerbound}.

Let $F=\coprod F^{j}$ be a Thurston norm minimizing surface dual
to $\psi_i$.  Let $\{\alpha_k\}$ be the set of boundary components
of $F$.  Suppose that $\alpha_k$ and $\alpha_l$ are parallel and
have opposite orientation.  Then we can glue an annulus along
$\alpha_k$ and $\alpha_l$ to get a new surface whose relative
homology class and $\chi\_$ are unchanged.  Altering our surface
in this way, we can assume that there is exactly one $k_0$ such
that $\alpha_{k_0} \cdot \mu_i =
 1$ and $\alpha_k \cdot \mu_j=0$ for all $k$ whenever $j \neq i$.
Secondly, we can assume that all the components of $F$ have
boundary since every closed surface is zero in
$H_2\left(X,\partial X;\mathbb{Z}\right)$. Now, let $F_i$ be the
 connected component of $F$ having $\alpha_{k_0}$ as one of its boundary components.
 Then $F_i$ represents the
 same relative homology class as $F$ and
 $\chi\_\left(F_i\right) \leq \chi\_\left(F\right)$.  Thus
 $F_i$ is a Thurston norm minimizing surface dual to $\psi$
 which is connected and has non-trivial boundary.
\end{proof}

The \emph{ropelength} of a link is the quotient of its length by
it's thickness.  In \cite{CKS}, J. Cantarella, R. Kusner, and J.
Sullivan show that the \emph{minimal ropelength}
$R\left(L_i\right)$ of the $i^{th}$ component of a link $L=\coprod
L_i$ is bounded from below by $2\pi \left(1+\sqrt{\left\| \psi_i
\right\| _T}\right)$. Hence the higher-order degrees give
computable lower bounds for the ropelength of knots and links.

\begin{corollary}\label{cor:rl}Let $X=S^3-L$ and $\psi_i$
be as defined above.  For each $n \geq 0$,
\[ R\left(L_i\right) \geq 2\pi
\left(1+\sqrt{\delta_n\left(\psi_i\right)-1}\right).
\]Moreover, if $\beta_1\left(X\right) \geq 2$ or $n \geq 1$ (or both) then
\[ R\left(L_i\right) \geq 2\pi
\left(1+\sqrt{\overline{\delta}_n\left(\psi_i\right)}\right).
\]
\end{corollary}
\begin{proof}The first (respectively second) statement follows from the bound given in
\cite{CKS} and Corollary \ref{cor:linkbound} (respectively Theorem
\ref{delbarthm}).
\end{proof}

Although it seems that the second statement in the Corollary is
``stronger'', in practice the first statement is often more
useful.  That is, $\overline{\delta}_n=0$ whenever the rank is
positive hence gives no new information.  We exemplify this
phenomena in Example \ref{ex:ropelength}.

We would like to determine conditions that will guarantee that a
surface will not lift to the $n^{th}$-order cover of $X$.  We show
that if $\beta _{1}\left(  X\right)  \geq2$ then $r_{n}\left(
X\right)  =0$ guarantees that no homologically essential surface
can lift to the $n^{th}$-order cover of $X$. In particular, if
$r_{0}\left( X\right)  =0$ then $i_{\ast}\pi_{1}\left( F\right)
\nsubseteq G_{r}^{\left(  1\right) }$ so that $i_{\ast}
\pi_{1}\left( F\right)  \nsubseteq G_{r}^{\left( n+1\right)  }$
for all $n\geq0$ . If $\beta_{1}\left(  X\right) =1$ a surface
representing the generator of $H_{2}\left( X,\partial
X;\mathbb{Z}\right)  $ can only lift if the rational derived
series stabilizes at the first step, i.e. $G_{r}^{\left( 1\right)
}=G_{r}^{\left(  2\right) }=\cdots=G_{r}^{\left( n+1\right)  }$.

\begin{proposition}\label{positiverank}
\label{posbigrank}If there exists a compact, connected, orientable, two-sided
properly embedded surface $F\subseteq X$ with $\beta_{1}\left(  X\right)
\geq2$ such that $0\neq\left[  F\right]  \in H_{2}\left(  X,\partial
X;\mathbb{Z}\right)  $ and $i_{\ast}\pi_{1}\left(  F\right)  \subseteq
G_{r}^{\left(  n+1\right)  }$ then $r_{n}\left(  X\right)  \geq1$.
\end{proposition}

\begin{proof}
Let $Y=\overline{X\setminus\left(  F\times I\right)  }$, since
$\left[ F\right]  \neq0$ $F$ does not separate $X$. Hence $Y$ is
connected. Let $\gamma$ be a oriented simple closed curve that
intersects $F$ exactly once, then $G=\pi_{1}\left(  X\right)
=\left\langle \pi_{1}\left(  Y\right) ,\gamma|\text{ relations
from }\pi_{1}\left(  F\right)  \right\rangle $. If $\pi_{1}\left(
Y\right)  \subseteq G_{r}^{\left(  1\right)  }$ then
$G/G_{r}^{\left(  1\right)  }=\left\langle \gamma\right\rangle $
which contradicts $\beta_{1}\left(  X\right) \geq2$. This implies
that $\pi _{1}\left(  Y\right)  \nsubseteq G_{r}^{\left( 1\right)
}$ hence $\pi _{1}\left(  Y\right) \nsubseteq G_{r}^{\left(
n+1\right)  }$ for all $n\geq0$. Now we consider the
Mayer-Vietoris sequence
\begin{align*}
0\rightarrow\operatorname{Im}\left(  i_{\ast}\oplus
j_{\ast}\right) \rightarrow H_{1}&\left(  X;\mathcal{K}_{n}\right)
\rightarrow H_{0}\left( F_{-}\amalg F_{+};\mathcal{K}_{n}\right)
\rightarrow   \\ &H_{0}\left(  F\times I;\mathcal{K}_{n}\right)
\oplus H_{0}\left( Y;\mathcal{K}_{n}\right)  \rightarrow
H_{0}\left( X;\mathcal{K}_{n}\right)  \rightarrow0.
\end{align*}
Since $\pi_{1}\left(  Y\right)  \nsubseteq G_{r}^{\left(
n+1\right)  }$, $\pi_{1}\left(  Y\right)  \rightarrow
G\rightarrow\Gamma_{n}$ is a non-trivial coefficient system.
Therefore we have $H_{0}\left(  Y;\mathcal{K} _{n}\right)  =0$ and
$H_0\left(X;\mathcal{K}_n\right)=0$ by Lemma \ref{lem:H0}. We note
that $\operatorname{rk}_{\mathcal{K} _{n}}H_{0}\left(
F;\mathcal{K}_{n}\right)  =1$ since $\pi_{1}\left( F\right)
\subseteq G_{r}^{\left(  n+1\right)  }$. It follows that
\begin{align*}
r_{n}\left(  X\right)
&=\operatorname{rk}_{\mathcal{K}_{n}}H_{1}\left( X;\mathcal{K}
_{n}\right) \\
&=\operatorname{rk}_{\mathcal{K}_{n}}H_{0}\left(
F;\mathcal{K}_{n}\right) +\operatorname{rk}_{\mathcal{K}
_{n}}\operatorname{Im}\left( i_{\ast}\oplus j_{\ast}\right)
\\
&\geq1.
\end{align*}
\end{proof}In particular, if there is a non-trivial surface that lifts to the n$^{th}$
cover then $r_{n}\left(  X\right)  \geq1$.

\begin{corollary}\label{cor:delbar0}
If there exists a compact, connected, orientable, two-sided
properly embedded surface $F\subseteq X$ with $\beta_{1}\left(
X\right) \geq2$ such that $0\neq\left[  F\right]  \in H_{2}\left(
X,\partial X;\mathbb{Z}\right)  $ and $i_{\ast}\pi_{1}\left(
F\right)  \subseteq G_{r}^{\left(  n+1\right)  }$ then
$\overline{\delta}_{n}\left(  \psi\right)  =0$ for all $\psi\in
H^{1}\left( X;\mathbb{Z}\right)  .$
\end{corollary}

\begin{proof}
$r_{n}\left(  X\right)  \geq 1  $ implies that
$\overset{\_}{\delta}_{n}\left(  \psi\right)  =0$ for all $\psi\in
H^{1}\left(  X;\mathbb{Z}\right)  $.
\end{proof}

We now prove the main theorem of this section.

\begin{proof}[Proof of Theorem \ref{delbarthm}]We break the proof
up into two cases.

 Case 1: Let $X$ be a 3-manifold with
$\beta_{1}\left( X\right) \geq2$. Let $F=\cup F_{i}$ be a surface
dual to $\psi$ that is minimal with respect to $\left\|
\cdot\right\| _{T}$.  We can assume that $\left[ F_{i}\right]
\neq0$ for all $i$.   If any component of $F$, say $F_{j}$ lifts
to the $n^{th}$ rational derived cover of $X$, i.e. $\pi
_{1}\left( F_{j}\right) \subset G_{r}^{\left( n+1\right) }$ then
$\overline{\delta}_{n}\left( \psi\right)  =0$ by Corollary
\ref{cor:delbar0}.  Otherwise $\mathcal{N}_n\left(\psi\right)=0$
so by Theorem \ref{delbarthm} we have $\overline
{\delta}_{n}\left( \psi\right) \leq\delta_{n}\left( \psi\right)
\leq\left\| \psi\right\| _{T}$.

Case 2: Let $X$ be a 3-manifold with $\beta_{1}\left(  X\right)
=1$ and $\psi$ be a generator of $H^{1}\left( X;\mathbb{Z}\right)
$. Let $F=\cup F_{i}$ be a surface dual to $\psi$ that is minimal
with respect to $\left\| \cdot\right\| _{T}$.  Since
$\beta_{1}\left(  \ker \psi\right) <\infty$, we can represent $F$
as a connected surface with $\beta_{2}\left( F\right)
=\beta_{3}\left(  X\right) $ \cite[Proposition 6.1]{Mc}. Therefore
$\mathcal{N}_0\left(\psi\right)\leq 1 + \beta_3\left(X\right)$ so
by Theorem \ref{lowerbound} we have $\overline{\delta}_{0}\left(
\psi\right) =\delta_{0}\left( \psi\right)  \leq\left\|
\psi\right\| _{T}+1+\beta _{3}\left( X\right)  $.  Now suppose
$n\geq1$.  If $\pi_{1}\left( F\right) \nsubseteqq G_{r}^{\left(
2\right) }$ (hence $\pi_{1}\left( F\right) \nsubseteqq
G_{r}^{\left( n+1\right) }$) then
$\mathcal{N}_n\left(\psi\right)=0$ so the result follows from
Theorem \ref{lowerbound}.  Otherwise, by Proposition
\ref{pro:torsionis0}, $\overline {\delta}_{n}\left( \psi\right)
=\delta_{n}\left( \psi\right)=0$.  We remark that if the
higher-order degrees of $S^1 \times S^2$ and $S^1 \times D^2$ are
zero.

The last sentence in the theorem follows from the calculations in
Proposition \ref{fibered}. Note that $r_{n}\left(  X\right)  =0$
for fibered 3-manifolds so that $\overline{\delta}_{n}\left(
\psi\right) =\delta _{n}\left(  \psi\right)  $ for all $n$.
\end{proof}

To complete the proof of Theorem \ref{delbarthm}, we need prove
Proposition \ref{pro:torsionis0} which states that if a
homologically essential surface dual to $\psi$ lifts to the
$n^{th}$ order cover then $\mathcal{A}^{\psi}_i \left(X\right)=0$
for $i < n$.  This will be our main objective for the rest of this
section.

\begin{figure}[htbp]
\begin{center}
\includegraphics
{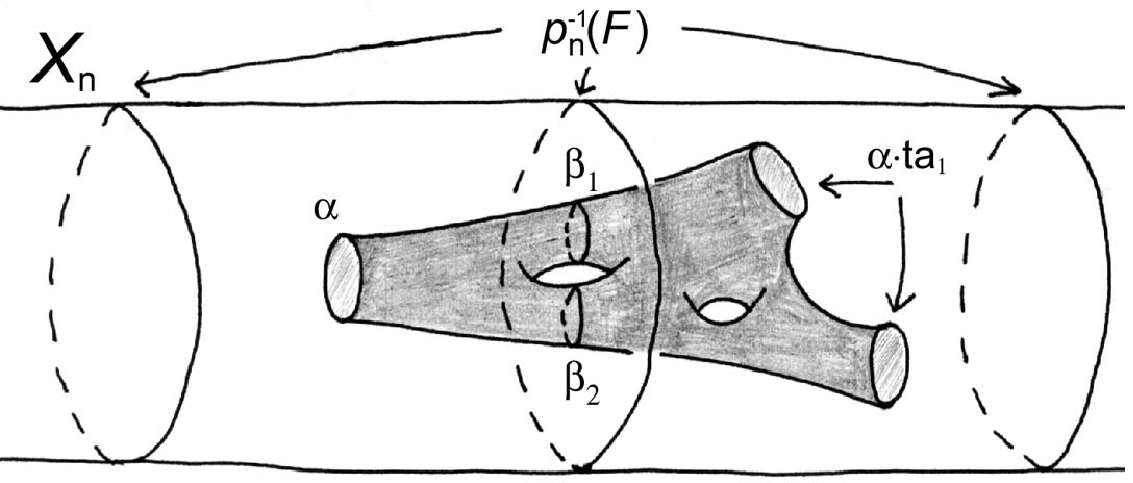}
\caption{$\mathcal{A}_n^{\psi}\left(X\right)$ is generated by
$H_{1}\left( F;\mathbb{K}_{n}\left[ t^{\pm1}\right]  \right)$}
\label{genbyF}
\end{center}
\end{figure}

 We begin by showing that $\mathcal{A}^{\psi}_n \left(X\right)$ is
generated by $H_{1}\left( F;\mathbb{K}_{n}\left[ t^{\pm1}\right]
\right)$. The idea behind the proof is simple. If $\alpha \neq 0$
is $\mathbb{K}_n\left[t^{\pm 1}\right]$-torsion then there exists
a $p\left(t\right) \in \mathbb{K}\left[t^{\pm 1}\right]$ such that
$\alpha p\left(t\right) = 0$.  Moreover, since $\mathbb{K}$ is a
(skew) field, we can assume that $p(t) = 1 + t a_1 + \cdots + t^m
a_m$ where $a_m \neq 0$. Thus $\alpha$ and $\alpha t a_1 + \cdots
+ \alpha t^m a_m$ cobound a surface, $S$ in $X_n$ (see Figure
\ref{genbyF}). Since the power of $t$ on each term of the latter
sum is positive, $S$ must intersect a lift of the surface $F$.
Hence $\alpha$ is homologous to the intersection of $S$ with the
lift of $F$. Note that in Figure \ref{genbyF}, $\alpha$ is
homologous  to $\beta_1 + \beta_2$.

\begin{lemma}\label{lem:generatedby}
$\mathcal{A}_n^{\psi}\left(X\right)
\subseteq\operatorname{Im}\left( i_{\ast}\right) $ where
\[i_{\ast} :H_{1}\left(  F;\mathbb{K}_{n}\left[  t^{\pm1}\right]
\right)  \rightarrow H_{1}\left(  X;\mathbb{K}_{n}\left[
t^{\pm1}\right]  \right).  \]
\end{lemma}

\begin{proof}
By Proposition \ref{pro:pres_A},
$\mathcal{A}_n^{\psi}\left(X\right) \subseteq
T\operatorname{Im}\left( j_{\ast}\right)$ where
$j_{\ast}:H_{1}\left( Y;\mathbb{K}_{n}\left[ t^{\pm1}\right]
\right) \overset{j_{\ast}}{\rightarrow}H_{1} \left( X
;\mathbb{K}_{n}\left[ t^{\pm 1}\right]  \right)$.  We will show
that $T\operatorname{Im}\left( j_{\ast}\right) \subseteq
\operatorname{Im}\left( i_{\ast}\right) $ which completes the
proof.

Let $\sigma_{X}\in T\operatorname{Im}\left(  j_{\ast}\right)  $
with $j_{\ast}\left( \sigma _{Y}\right)  =\sigma_{X}$. Since
$\sigma_{X}$ is $\mathbb{K}_{n}\left[ t^{\pm1}\right]$-torsion,
there exists $p\left(  t\right) \in\mathbb{K}_{n}\left[
t^{\pm1}\right] $ such that $\sigma_{X}p\left( t\right)  =0.$ We
have $j_{\ast}\left(  \sigma_{Y}p\left(  t\right) \right)
=j_{\ast}\left(  \sigma_{Y}\right)  p\left(  t\right)  =\sigma
_{X}p\left(  t\right)  =0$ so there exists $\sigma_{F}\in
H_{1}\left( F;\mathbb{K}_{n}\left[  t^{\pm1}\right]  \right)  $
such that $\eta\left( \sigma_{F}\right)  =\sigma_{Y}p\left(
t\right)  $. We can assume that $p\left(  t\right)
=1+tc_{1}+\cdots+t^{m}c_{m}$ since $\sigma_{X}p\left( t\right) =0$
if and only if $\sigma_{X}p\left(  t\right)  u=0$ for any unit
$u\in\mathbb{K}_{n}\left[  t^{\pm1}\right]  $. Now $H_{1}\left(
F;\mathbb{K}_{n}\left[  t^{\pm1}\right]  \right)  \simeq
H_{1}\left( F;\mathbb{K}_{n}\right)
\otimes_{\mathbb{K}_{n}}\mathbb{K}_{n}\left[ t^{\pm1}\right]  $
and $H_{1}\left(  Y;\mathbb{K}_{n}\left[  t^{\pm1}\right] \right)
\simeq H_{1}\left(  Y;\mathbb{K}_{n}\right)
\otimes_{\mathbb{K}_{n} }\mathbb{K}_{n}\left[  t^{\pm1}\right]  $
so every element in $H_{1}\left( F;\mathbb{K}_{n}\left[
t^{\pm1}\right]  \right)  $ (resp. $H_{1}\left(
Y;\mathbb{K}_{n}\left[  t^{\pm1}\right] \right)  )$ has the form
$\overset{\infty}{\underset{i=-\infty}{\sum}}\alpha_{i}\otimes
t^{i}$ (resp.
$\overset{\infty}{\underset{i=-\infty}{\sum}}\beta_{i}\otimes
t^{i}$) such that $\alpha_{i}\in H_{1}\left(
F;\mathbb{K}_{n}\right)  $ (resp. $\beta _{i}\in H_{1}\left(
Y;\mathbb{K}_{n}\right)  $) and there are only finitely many
non-zero $\alpha_{i}$ (resp. $\beta_{i}$). We write $\sigma
_{F}=\overset{\infty}{\underset{i=-\infty}{\sum}}\alpha_{i}\otimes
t^{i}$ and (as with $p\left(  t\right)  $) we may write
$\sigma_{Y}=\overset{k} {\underset{i=0}{\sum}}\beta_{i}\otimes
t^{i}$. Using this notation we now have
\begin{align*}
\sigma_{Y}p\left(  t\right)
&=\overset{k}{\underset{i=0}{\sum}}\left(
\beta_{i}\otimes t^{i}\right)  \overset{m}{\underset{j=0}{\sum}}t^{j}c_{j}\\
&=\overset{k+m}{\underset{l=0}{\sum}}\underset{i+j=l}{\sum}\left(
\beta
_{i}\left(  c_{j}\right)  ^{t}\otimes t^{i+j}\right) \\
&=\overset{k+m}{\underset{l=0}{\sum}}\left(
\underset{i+j=l}{\sum}\beta _{i}\left(  c_{j}\right)  ^{t}\right)
\otimes t^{l}
\end{align*}
and
\begin{align*}
\eta\left(  \sigma_{F}\right)
&=\overset{\infty}{\underset{i=-\infty}{\sum}
}\eta\left(  \alpha_{i}\otimes t^{i}\right) \\
&=\overset{\infty}{\underset{i=-\infty}{\sum}}\left(  i_{-}\right)
_{\ast }\left(  \alpha_{i}\otimes t^{i}\right)  -\left(
i_{+}\right)  _{\ast}\left(
\alpha_{i}\otimes t^{i}\right)  t\\
&=\overset{\infty}{\underset{i=-\infty}{\sum}}\left(  i_{-}\right)
_{\ast }\left(  \alpha_{i}\right)  \otimes t^{i}-\left(
i_{+}\right)  _{\ast}\left(
\alpha_{i}\right)  \otimes t^{i+1}\\
&=\overset{\infty}{\underset{i=-\infty}{\sum}}\left(  i_{-}\right)
_{\ast }\left(  \alpha_{i}\right)  \otimes t^{i}-\left(
i_{+}\right)  _{\ast}\left(
\alpha_{i-1}\right)  \otimes t^{i}\\
&=\overset{\infty}{\underset{i=-\infty}{\sum}}\left(  \left(
i_{-}\right) _{\ast}\left(  \alpha_{i}\right)  -\left(
i_{+}\right)  _{\ast}\left( \alpha_{i-1}\right)  \right)  \otimes
t^{i}.
\end{align*}

Recall that $\eta\left(  \sigma_{F}\right)  =\sigma_{Y}p\left(
t\right)  $ which implies that
$\underset{i+j=l}{\sum}\beta_{i}\left( a_{j}\right) ^{t}=\left(
i_{-}\right)  _{\ast}\left( \alpha_{l}\right)  -\left(
i_{+}\right)  _{\ast}\left( \alpha_{l-1}\right)  $ for all $0\leq
l\leq k+m$. \ In particular $c_{0}=1$ so when $l\leq k$ we can
write $\beta_{l}$ as a combination of $\beta_{i}$ and $\left(
i_{-}\right)  _{\ast}\left( \alpha_{l}\right)  +\left(
i_{+}\right)  _{\ast}\left( \alpha_{l-1}\right) $ with $i<l$. That
is,
\[
\beta_{l}\otimes t^{l}=\left(  i_{-}\right)  _{\ast}\left(
\alpha_{l}\right) -\left(  i_{+}\right)  _{\ast}\left(
\alpha_{l-1}\right)  \otimes
t^{l}-\underset{i+j=l,i<l}{\sum}\beta_{i}\left(  c_{j}\right)
^{t}\otimes t^{l}.
\]

We will prove by induction that $j_{\ast}\left(  \beta_{l}\otimes
t^{l}\right)  \in\operatorname{Im}\left(  i_{\ast}\right)  $ for
each $l$ implying that $\sigma_{X}=j_{\ast}\left(
\sigma_{Y}\right)  =j_{\ast}\left( \sum\nolimits_{0\leq l\leq
k}\beta_{l}\otimes t^{l}\right)  \in \operatorname{Im}\left(
i_{\ast}\right)  $ which completes the proof. We first note that
\begin{align*}
j_{\ast}\left(  \left(  i_{-}\right)  _{\ast}\left(
\alpha_{l}\right)  +\left(  i_{+}\right)  _{\ast}\left(
\alpha_{l-1}\right) \otimes t^{l}\right)  &=j_{\ast}\left(  \left(
i_{-}\right)  _{\ast}\left( \alpha_{l}\otimes t^{l}\right) +\left(
i_{+}\right)  _{\ast}\left( \alpha_{l-1}\otimes t^{l}\right)
\right)  \\&=i_{\ast}\left(  \alpha_{l}\otimes
t^{l}-\alpha_{l-1}\otimes t^{l}\right)  \\&=i_{\ast}\left(
\alpha_{l} -\alpha_{l-1}\otimes t^{l}\right).\end{align*} It
follows that $\beta_{0} \otimes1=\left(  i_{-}\right)
_{\ast}\left( \alpha_{0}\right) -\left( i_{+}\right) _{\ast}\left(
\alpha_{-1}\right) \otimes1=i_{\ast}\left(
\alpha_{0}-\alpha_{-1}\otimes1\right) \in\operatorname{Im}\left(
i_{\ast }\right)  $. Now we assume that $\beta_{i}\otimes
t^{i}=i_{\ast}\left( \gamma_{p}\right)  $ for all $i\leq l-1$ so
that
\begin{align*}
j_{\ast}\left(  \beta_{l}\otimes t^{l}\right)  &=j_{\ast}\left(
\left( i_{-}\right)  _{\ast}\left(  \alpha_{l}\right)  -\left(
i_{+}\right)  _{\ast }\left(  \alpha_{l-1}\right)  \otimes
t^{l}\right)  -\underset{i+j=l,i<l} {\sum}j_{\ast}\left(
\beta_{i}\left(  c_{j}\right)  ^{t}\otimes t^{l}\right)
\\
&=i_{\ast}\left(  \alpha_{l}-\alpha_{l-1}\otimes t^{l}\right)
-\underset
{i+j=l,i<l}{\sum}j_{\ast}\left(  \beta_{i}\otimes t^{l}\right)  c_{j}\\
&=i_{\ast}\left(  \alpha_{l}-\alpha_{l-1}\otimes t^{l}\right)
-\underset
{i+j=l,i<l}{\sum}i_{\ast}\left(  \gamma_{i}\right)  c_{j}\\
&\in\operatorname{Im}\left(  i_{\ast}\right)  .
\end{align*}
\end{proof}

We can use this to show that if $\psi$ is dual to a union of
surfaces in $X$ whose fundamental groups all include into the
$\left(  n+1\right)  ^{st}$ rational derived subgroup of
$G=\pi_{1}\left(  X\right)  $ then $\delta _{i}\left(  \psi\right)
=0$ for $i\leq n-1$.
\begin{proposition}\label{pro:torsionis0}
If there exists a union of properly embedded surfaces $F=\cup
F_{j}$ in $X$ with $\left[  F\right]  \in H_{2}\left(  X,\partial
X;\mathbb{Z}\right)  $ dual to $\psi\in H^{1}\left(
X;\mathbb{Z}\right)  $ such that for all $j$, $\pi_{1}\left(
F_{j}\right)  \subseteq G_{r}^{\left(  n+1\right)  }$ then
$\mathcal{A}_i^{\psi}\left(X\right) =0$ whenever $0\leq i\leq
n-1$.
\end{proposition}

\begin{proof}
Consider the following diagram of abelian groups:
\begin{equation}
\begin{diagram} H_{1}\left( F_{\Gamma_i} \right) & \rTo^{i_*} & H_{1}
\left( X_{\Gamma_{i}}\right) \\ \dTo^{\rho_F} & & \dTo_{\rho_X} \\
H_{1} \left( F;\mathbb{K}_{i}\left[ t^{\pm1}\right]\right) &
\rTo^{i_*} & H_{1} \left( X;\mathbb{K}_{i}\left[
t^{\pm1}\right]\right) \\ \end{diagram} \label{comm}\vspace{12pt}
\end{equation}
where $i_{\ast}$ is induced by the inclusion map $i:F\rightarrow
X$ and $\rho_{F}\left(  \left[  \sigma\right]  \right)  =\left[
\sigma \otimes1\right]  $ (similarly for $\rho_{X}$). First we
observe that if $i\leq n-1$, every class in $\pi_{1}\left(
F_{\Gamma_{i}}\right)  $ gets mapped into $G_{r}^{\left(
n+1\right)  }\subseteq G_{r}^{\left(  i+2\right) }$ by
$i_{\ast}\circ p_{\ast}$ hence is zero in $G_{r}^{\left(
i+1\right) }\left/  G_{r}^{\left(  i+2\right)  }\right.
=H_{1}\left(  X_{\Gamma_{i} }\right)  \left/  \left\{
\mathbb{Z}\!-\!\operatorname{torsion}\right\}  \right.  $. \
Therefore $i_{\ast}:H_{1}\left(  F_{\Gamma_{i}}\right) \rightarrow
H_{1}\left(  X_{\Gamma_{i}}\right)  $ maps $H_{1}\left(
F_{\Gamma_{i} }\right)  $ to the $\mathbb{Z}$-torsion subgroup of
$H_{1}\left( X_{\Gamma_{i}}\right)  $. Since $H_{1}\left(
X;K_{i}\left[ t^{\pm 1}\right]  \right)  $ is $\mathbb{Z}$-torsion
free, $\rho_{X}\circ i_{\ast} =0$.
\[
\begin{diagram} F_{\Gamma_n} & \rTo& X_{\Gamma_n} \\ \dTo^p & & \dTo
_p \\ F & \rTo& X \\ \end{diagram} \vspace{12pt}\]
 Since $\left(
\ref{comm}\right)  $ commutes, the image of $\rho_{F}$ goes to
zero under $i_{\ast}$. Therefore if $\left[ \sigma\otimes p\left(
t\right)  \right]  $ is an element of $H_{1}\left(
F;\mathbb{K}_{i}\left[ t^{\pm1}\right]  \right)  $,
$i_{\ast}\left(  \left[  \sigma\otimes p\left( t\right)  \right]
\right)  =i_{\ast}\left(  \left[  \sigma\otimes1\right] \right)
p\left(  t\right)  =0p\left(  t\right)  =0$ ($i_{\ast}$ is a
$\mathbb{K}_{i}\left[  t^{\pm1}\right]  $-module homomorphism). By
Lemma \ref{lem:generatedby}, $\mathcal{A}_i^{\psi}\left(X\right)$
is generated by $\operatorname{Im}\left( i_{\ast }\right)  =0$
hence $\mathcal{A}_i^{\psi}\left(X\right)=0$.
\end{proof}

\begin{corollary}
Let $X$ be a 3-manifold with $\beta_{1}\left(  X\right)  =1$ and
$F$ a surface dual to a generator of $H^{1}\left(
X;\mathbb{Z}\right)  $. If $\pi_{1}\left(  F\right)  \subseteq
G_{r}^{\left(  2\right)  }$ then $G_{r}^{\left(  i\right)
}=G_{r}^{\left(  i+1\right)  }$ for all $i\geq1$.
\end{corollary}

\begin{proof}
Since $\beta_{1}\left(  X\right)  =1$, $r_n\left(X\right) =0$ by
Proposition \ref{upperboundrank}. That is, $H_{1}\left(
X;\mathbb{K}_{i}\left[ t^{\pm1}\right] \right)$ is a torsion
module.   When $i=0$, $\mathbb{K} _{0}=\mathbb{Q}$ so that
\[TH_{1}\left( X;\mathbb{K}_{0}\left[ t^{\pm 1}\right]  \right)
=H_{1}\left( X;\mathbb{Q}\left[ t^{\pm1}\right] \right)
=H_{1}\left( X_{\Gamma_{0}}\right) \otimes\mathbb{Q}=G_{r}^{\left(
1\right) }/G_{r}^{\left( 2\right) }\otimes\mathbb{Q}.\] If
$\pi_{1}\left(  F\right) \subset G_{r}^{\left( 2\right)  }$,
Proposition \ref{pro:torsionis0} implies $TH_{1}\left(
X;\mathbb{K}_{0}\left[ t^{\pm1}\right] \right) =0$.   Since
$G_{r}^{\left( 1\right) }/G_{r}^{\left( 2\right) }$ is
$\mathbb{Z}$-torsion free, $G_{r}^{\left( 1\right) }/G_{r}^{\left(
2\right) } \rightarrow G_{r}^{\left( 1\right) }/G_{r}^{\left(
2\right) } \otimes\mathbb{Q}$ sending $g \mapsto g \otimes 1$ is a
monomorphism.
 Therefore $G_{r}^{\left( i\right) }=G_{r}^{\left( i+1\right) }$ for all
$i\geq1$.
\end{proof}

\section{Realization Theorem}\label{realizationsection}
We are ready to prove that the invariants $\delta_n$ give much
more information than the classical invariants. In fact, we subtly
alter 3-manifolds to obtain new 3-manifolds with striking
behavior.  Cochran proves this result when
$\beta_1\left(X\right)=1$ \cite{Co}.

\begin{theorem}
\label{increase}For each $m\geq1$ and $\mu\geq2$ there exists a
3-manifold $X$ with $\beta_{1}\left(  X\right)  =\mu$ such that
\[
\left\|  \psi\right\|  _{A}=\delta_{0}\left(  \psi\right)
<\delta_{1}\left( \psi\right)  <\cdots<\delta_{m}\left(
\psi\right)  \leq\left\|  \psi\right\| _{T}
\]
for all $\psi\in H^{1}\left( X;\mathbb{Z}\right)$.  Moreover, $X$
can be chosen so that it is closed, irreducible and has the same
classical Alexander module as a 3-manifold that fibers over $S^1$.
\end{theorem}

The proof of this will be an application of the following more
technical theorem.  We will postpone the proof until later in the
section.  Theorem \ref{realization} is a tool that will allow us
to subtly alter 3-manifolds in order to construct new 3-manifolds
whose degrees are unchanged up to the $n^{th}$ stage but increase
at the $n^{th}$ stage.

\begin{theorem}
[Realization Theorem]\label{realization}Let $X$ be a compact,
orientable 3-manifold with $G=\pi_{1}\left(  X\right)  $ and
$G_{r}^{\left( n\right)  } /G_{r}^{\left(  n+1\right)  }\neq0$ for
some $n\geq0$.  Let $\left[x\right]$ be a primitive class in
$H_1\left(X;\mathbb{Z}\right)$.  Then for any positive integer
$k$, there exists a 3-manifold $X\hspace{-1 pt}\left( n,k\right)$
homology cobordant to $X$ such that \newline $(1)$ \[
\frac{G}{G_{r}^{\left(  i+1\right)  }} \cong
\frac{H}{H_{r}^{\left(  i+1\right)  }} \text{\quad for \quad}
    0 \leq i \leq n-1\]
    where $H=\pi_1\left(X\left(n,k\right)\right)$ and
\newline
$(2)$
    \[
\delta_{n}^{X\left(  n,k\right)  }\left(  \psi\right)
\geq\delta_{n} ^{X}\left(  \psi\right)  +k\left|  p\right|  .
\]
for any $\psi \in H^{1}\left(  X\left( n,k\right)
;\mathbb{Z}\right) $ with $\psi\left( x\right) =t^{p}$.
\end{theorem}
\begin{proof}
Let $X$ be a compact 3-manifold with $G_{r}^{\left(  n\right)  }
/G_{r}^{\left(  n+1\right)  }\neq0$. $G=G_{r}^{\left(  1\right) }$
implies that $G_{r}^{\left(  n\right)  }=G_{r}^{\left( n+1\right)
}$ hence our hypothesis guarantees that $\beta_{1}\left( X\right)
\geq1$. Since $\left[x\right]$ is a primitive class in
$H_1\left(X;\mathbb{Z}\right)$, we can present $G$ as
\[G \cong
\left\langle x_{1},\ldots,x_{\mu},y_{1},\ldots
,y_{l}|R_{1},\ldots,R_{m}\right\rangle \]
 where $y_{i}\in
G_{r}^{\left( 1\right)  }$, $x_1=x$, and $\{ \left[x_1\right],
\ldots,\left[x_{\mu}\right]\}$ is a basis for $G/G_{r}^{\left(
1\right) }$.

Begin by adding a 1-handle to $X\times I$ to obtain a 4-manifold
$V$ with boundary $\left(  X\sqcup X^{\prime}\right) \cup\left(
\partial X\times I\right)  $ where $X^{\prime}$ is obtained from
$X$ by taking the connect sum with $S^{1}\times S^{2}$. Then
$\pi_{1}\left( V\right) \cong\pi_{1}\left( X^{\prime}\right)
\cong G\ast\left\langle z\right\rangle $ where $z$ is the
generator of $\pi_{1}\left( S^{1}\times S^{2}\right) $. Choose a
non-trivial element $B\in G_{r}^{\left(  n\right) }-G_{r}^{\left(
n+1\right)  }$, and let $w=zx^{-1}$ and $\alpha=\left[
A_{k},B\right]  $ where $A_{k}$ is defined inductively as
\begin{align*}
A_{1}&=w\\
A_{k}&=\left[  A_{k-1},x\right]  \text{ for }k\geq2.
\end{align*}
Now add a 2-handle to $V$ along a curve $c$ (any framing) embedded
in $X^{\prime}$ representing $w\left[  x,\alpha\right]  $ to
obtain a 4-manifold $W$ with boundary $\left(X\sqcup-X\left(
n,k\right) \right)  \cup\left(
\partial X\times I\right)  $. Let $E=\pi_{1}\left(  W\right)  $, $H=\pi
_{1}\left(  X\left(  n,k\right)  \right)  $ and denote by $i$ and $j$ the
inclusion maps of $X$ and $X\left(  n,k\right)  $ into $W$ respectively.

Adding the 2-handle to $X^{\prime}$ kills the element $w\left[
x,\alpha \right]  $ in $G\ast\left\langle z\right\rangle \cong
G\ast\left\langle w\right\rangle $ so $E\cong\left\langle
G,w|w\left[  x,\alpha\right] \right\rangle =\left\langle
x_{1},\ldots,x_{\mu},y_{1},\ldots,y_{l}
,w|R_{1},\ldots,R_{m},w\left[  x,\alpha\right]  \right\rangle $.
We see that $X\left(  n,k\right)  $ is the 3-manifold obtained by
performing Dehn surgery (with integer surgery coefficient
corresponding to the framing of the 2-handle) along the curve $c$.
Let $\gamma$ be the meridian curve to $c$ in $X\left( n,k\right)
$. The dual handle decomposition of $W$ rel $X\left( n,k\right) $
is obtained by adding to $X\left( n,k\right)  $ a 0-framed
2-handle along $\gamma$ and 3-handle. This gives us
$E\cong\left\langle H|\gamma\right\rangle $.

We show that
\begin{equation}\label{equ:samegroups2}\frac{G}{\left[
G_{r}^{\left( i\right) },G_{r}^{\left( i\right)  }\right]
}\overset{\cong}{\longrightarrow}\frac{E}{\left[ E_{r}^{\left(
i\right)  },E_{r}^{\left(  i\right)  }\right]  }\overset{\cong
}{\longleftarrow}\frac{H}{\left[  H_{r}^{\left(  i\right)
},H_{r}^{\left( i\right)  }\right]  }\end{equation} for $0\leq
i\leq n$. Using Lemma \ref{derivedquotient}, this will imply that
\begin{equation} \label{equ:samegroups}
\frac{G}{G_{r}^{\left(  i+1\right)
}}\overset{\cong}{\longrightarrow}\frac {E}{E_{r}^{\left(
i+1\right)  }}\overset{\cong}{\longleftarrow}\frac{H}
{H_{r}^{\left(  i+1\right)  }}.
\end{equation}
There is a surjective map $\operatorname{pr}:\left\langle
G,w|w\left[ x,\alpha\right]  \right\rangle \twoheadrightarrow
G\,$\ defined by killing $w$ so that $\operatorname{pr}\circ
i_{\ast}=id_{G}$. Consider the induced maps
\[
\frac{G}{\left[  G_{r}^{\left(  i\right)  },G_{r}^{\left( i\right)
}\right]
}\overset{\overline{i}_{\ast}}{\longrightarrow}\frac{E}{\left[
E_{r}^{\left( i\right)  },E_{r}^{\left(  i\right)  }\right]
}\overset{\overline
{\operatorname{pr}}}{\twoheadrightarrow}\frac{G}{\left[
G_{r}^{\left( i\right)  },G_{r}^{\left(  i\right)  }\right] }.
\]

We will show by induction that $w\in\left[  E_{r}^{\left( i\right)
} ,E_{r}^{\left(  i\right)  }\right]  $ for $0\leq i\leq n$. Since
$w=\left[ \left[  A_{k},B\right]  ,x\right]  $, it is clear that
$w\in\left[ E_{r}^{\left(  0\right)  },E_{r}^{\left( 0\right)
}\right]  $. Now suppose that $w\in\left[  E_{r}^{\left(
i-1\right)  },E_{r}^{\left(  i-1\right) }\right]  $ for some
$i\leq n$. Since $A_{k}=\left[  A_{k-1},x\right]  $ and $A_{1}=w$,
we have $A_{k}\in\left[  E_{r}^{\left(  i-1\right)  }
,E_{r}^{\left( i-1\right)  }\right]  \subseteq E_{r}^{\left(
i\right)  }$. \ Moreover, since $B\in G_{r}^{\left(  n\right)  }$,
we have $B\in E_{r}^{\left(  n\right)  }\subseteq E_{r}^{\left(
i\right)  }$ for $i\leq n$. \ Therefore $\left[  A_{k},B\right]
\in\left[ E_{r}^{\left(  i\right) },E_{r}^{\left(  i\right)
}\right]  $ for $i\leq n$ and hence $w \in\left[  E_{r}^{\left(
i\right)  } ,E_{r}^{\left(  i\right)  }\right]  \subseteq
E_{r}^{\left( i+1\right)  }$.  It follows that
$\overline{\operatorname{pr}}$ is an isomorphism.  Since
$\overline{\operatorname{pr}}\circ\overline{i}_{\ast}$ is an
isomorphism, $\overline{i}_{\ast}$is an isomorphism for $0\leq
i\leq n$.

Now consider the maps $\frac{H}{\left[  H_{r}^{\left(  i\right)  }
,H_{r}^{\left(  i\right)  }\right]  }\overset{\overline{j}_{\ast}
}{\longrightarrow}\frac{E}{\left[  E_{r}^{\left(  i\right)
},E_{r}^{\left( i\right)  }\right]  }$ where now we are
considering $E$ as the group $\left\langle H|\gamma\right\rangle
$. We will show that $\gamma\in\left[ H_{r}^{\left(  i\right)
},H_{r}^{\left(  i\right)  }\right]  $ for $0\leq i\leq n $ hence
the above map will be an isomorphism. Recall that $X\left(
n,k\right)  $ can be obtained from $X$ by first doing Dehn surgery
on a 0-framed unlinked trivial knot in $X$ to get the manifold
$X^{\prime }=X\#\left(  S^{1}\times S^{2}\right)  $ and then
performing Dehn surgery along a curve $c$ representing $w\left[
x,\alpha\right]  $ in $X^{\prime}$. Let $Y=X^{\prime}-N\left(
c\right)  $ be the 3-manifold obtained by removing a regular
neighborhood of $c$ in $X^{\prime}$. We use the notations
$P=\pi_{1}\left(  Y\right)  $, $K=\pi_{1}\left(  X^{\prime}\right)
\cong G\ast\left\langle z\right\rangle $, and $l:Y\rightarrow
X^{\prime}$ be the inclusion map. Let $\gamma$ be the meridian of
$c$ based at $x_0$ as in Figure \ref{cobounddisc}. We show that
$\gamma\in\left[ P_{r}^{\left(  i\right)  } ,P_{r}^{\left(
i\right)  }\right]  $ for $0\leq i\leq n$
 which implies that
$\gamma\in\left[ H_{r}^{\left(  i\right)  },H_{r}^{\left( i\right)
}\right]  $ since $P\twoheadrightarrow H$.

\begin{figure}[htbp]
\begin{center}
\includegraphics{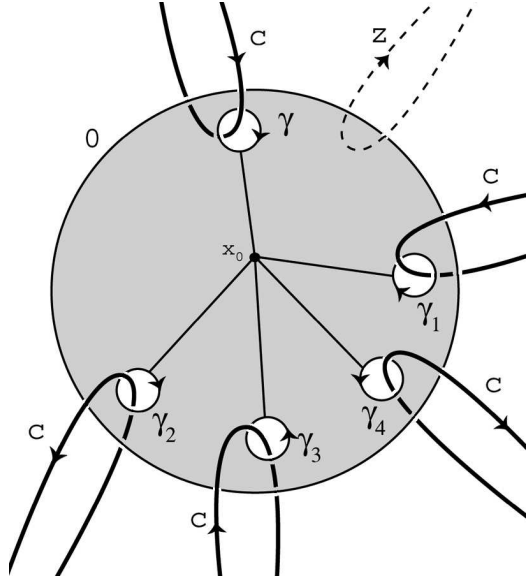}
\caption{$\gamma$ and the $\gamma_i^{\pm1}$s cobound a punctured
sphere} \label{cobounddisc}
\end{center}
\end{figure}

To begin, we will show that
\begin{equation}
\gamma=\left[  \gamma,u_{1}\right]  ^{v_{1}}\cdots\left[
\gamma,u_{2^{k} }\right]  ^{v_{2^{k}}}\left[
\lambda_{1},\lambda_{2}\right]  ^{w_{1}} \cdots\left[
\lambda_{m-1},\lambda_{m}\right]  ^{w_{m/2}}\label{formofgamma}
\end{equation}
where $l_{\ast}\left(  u_{j}\right)  \in K_{r}^{\left(  n\right)
}$ and $\lambda_{j}\in Ncl\left\langle \gamma\right\rangle
=\ker\left(  l_{\ast }:P\twoheadrightarrow K\right)  $. Using
this, we will show that the induced map $l_{\ast}:P_{r}^{\left(
i\right)  }\twoheadrightarrow K_{r}^{\left( i\right)  }$ is
surjective for $0\leq i\leq n+1$.
Assuming these two statements to be true for
now, we shall prove by induction on $i$ that $\gamma\in\left[
P_{r}^{\left( i\right) },P_{r}^{\left( i\right) }\right]  $ for
$0\leq i\leq n$ as desired. It is clear from $\left(
\ref{formofgamma}\right)  $ that $\gamma\in\left[ P_{r}^{\left(
0\right) },P_{r}^{\left( 0\right) }\right]  $. Now suppose that
$\gamma\in\left[ P_{r}^{\left( i-1\right) },P_{r}^{\left(
i-1\right)  }\right] \subseteq P_{r}^{\left( i\right)  }$ for some
$i\leq n$.
Then $\lambda_{j}\in
P_{r}^{\left( i\right) }$ for all $j$. Since
$l_{\ast}:P_{r}^{\left(  i\right) }\twoheadrightarrow
K_{r}^{\left(  i\right)  }$ is surjective and $l_{\ast }\left(
u_{j}\right)  \in K_{r}^{\left(  n\right)  }$, it follows that
$u_{j}=p_{j}\Lambda_{j}$ where $p_{j}\in P_{r}^{\left(  i\right)
}$ and $\Lambda_{j}\in Ncl\left\langle \gamma\right\rangle
\subseteq P_{r}^{\left( i\right)  }$. Thus by $\left(
\ref{formofgamma}\right)$, $\gamma\in\left[ P_{r}^{\left( i\right)
},P_{r}^{\left( i\right) }\right]  $ for $0\leq i\leq n$.

Let $c$ be a curve representing the element
\begin{align*}
w\left[  x,\alpha\right]  &=z\alpha x^{-1}\alpha^{-1}\\
&=zA_{k}BA_{k}^{-1}B^{-1}x^{-1}BA_{k}B^{-1}A_{k}^{-1}.
\end{align*}
For simplicity, we prove the case when $c$ intersects the cosphere
(belt sphere) of the 1-handle attached to $X\times I$ exactly
$1+2^{k+1}$ times. The proof can be modified for the case where
$c$ intersects the cosphere more than $1+2^{k+1}$ times. The
$2^{k+1}$ intersections are a result of the $2^{k-1}$ occurrences
of $z$ and $z^{-1}$ in $A_{k}$ and the first intersection is a
result the first $z$ that occurs in $z\alpha x^{-1}\alpha^{-1}$.
Note that the only occurrences of $z$ in $\alpha
x^{-1}\alpha^{-1}$ show up in $A_{k}$ since $B$ is an element of
$G$.  Let $\gamma_{j}$ be the meridian of $c$ based at $x_0$
corresponding to the $j^{th}$ occurrence of $z$ or $z^{-1}$ in
$\alpha x^{-1}\alpha^{-1}$ as shown in Figure \ref{cobounddisc}.

Before proceeding, we sketch the idea of the next part of the
proof. First we note that $\gamma$ is a product of
$\gamma_j^{\pm1}$. We can pair each $\gamma_j$ with $\gamma_{2^k
+1 -j}$ since they bound an annulus as in Figure
\ref{fig:coboundannulus}.  Thus they are related by $\gamma_{2^k
+1 -j}=u^{-1}\gamma_{j}u$ and hence $\gamma_{j} \gamma_{2^k +1
-j}^{-1}=\left[\gamma_j,u^{-1}\right]$.  We show that $l_{\ast}
\left(u\right)$ is in the $n^{th}$ term of the derived series of
$K$. Since $\gamma_j$ is a conjugate of $\gamma$, $\gamma$ is a
product of commutators which can be written as
$\left[\gamma,u\right]^{v}$ with $l_{\ast} \left(u\right) \in
K_r^{\left(n\right)}$.

\begin{figure}[htbp]
\begin{center}
\includegraphics{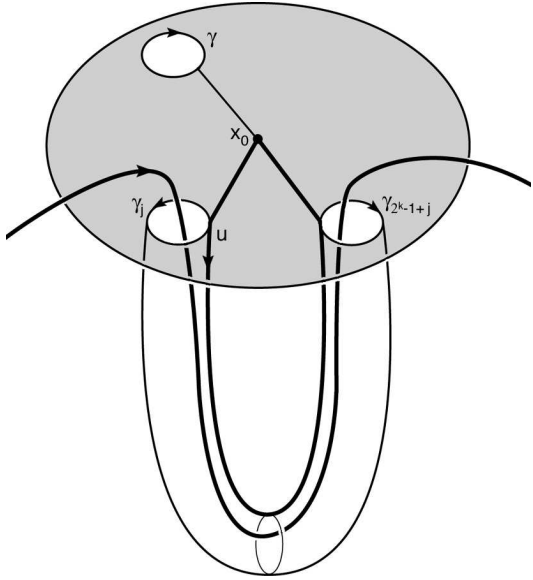}
\caption{$\gamma_j$ and $\gamma_{2^k +1 -j}$ differ by a
conjugation} \label{fig:coboundannulus}
\end{center}
\end{figure}

Let $ a^{b}=b^{-1}ab$. Since the longitude of the unknot is
trivial in $P$, we see that $\gamma$ is equal to a product of the
$\gamma_{j}^{\pm1}$ as in Figure \ref{cobounddisc}.  Moreover, we
can order the $\gamma_{j}$ as we choose since switching
$\gamma_{i}$ and $\gamma_{j}$ only changes the element by a
commutator of elements of $Ncl\left\langle \gamma\right\rangle $.
For example,
\begin{align*}
\lambda_{1}\gamma_{i}\gamma_{j}\lambda_{2}&=\lambda_{1}\gamma_{j}\gamma
_{i}\left[  \gamma_{i}^{-1},\gamma_{j}^{-1}\right]  \lambda_{2}\\
&=\lambda_{1}\gamma_{j}\gamma_{i}\lambda_{2}\left[
\gamma_{i}^{-1},\gamma _{j}^{-1}\right] ^{\lambda_{2}}.
\end{align*}
Hence we see that
\[
\gamma=\left(\prod_{j=1}^{2^{k-1}}
\gamma_{j}^{\pm1}\gamma_{2^{k}+1-j}^{\pm 1}\right)\left(
\prod_{j=1}^{2^{k-1}} \gamma_{2^{k}+j}^{\pm1}\gamma
_{2^{k+1}+1-j}^{\pm1}\right)  \left[
\lambda_{1},\lambda_{2}\right]  ^{w_{1} }\cdots\left[
\lambda_{m-1},\lambda_{m}\right]  ^{w_{m/2}}
\]
where $\lambda_{j}\in Ncl\left\langle \gamma\right\rangle $. The
chosen ordering will become clear in the next paragraph.

If the $j^{th}$ occurrence of a $z^{\pm1}$is a $z$ then
$\gamma_{j} =\gamma^{p_{j}}$ where $l_{\ast}\left(  p_{j}\right)
=z\omega_{j}$ and $\omega_{j}$ is the word that occurs in $\alpha
x^{-1}\alpha^{-1}$ up to but not including the $j^{th}$ $z$.
Whereas, if the $j^{th}$ occurrence of a $z^{\pm1}$is a $z^{-1}$
then $\gamma_{j}=\gamma^{p_{j}}$ where $l_{\ast }\left(
p_{j}\right)  =z\omega_{j}z^{-1}$ where $\omega_{j}$ is the word
that occurs in $\alpha x^{-1}\alpha^{-1}$ up to but not including
the $j^{th}$ $z^{-1}$. Now we consider the case $1\leq
j\leq2^{k-1}$. The $j^{th}$ occurrence of $z^{\pm1}$ in $\alpha
x^{-1}\alpha^{-1}=A_{k}
BA_{k}^{-1}B^{-1}x^{-1}BA_{k}B^{-1}A_{k}^{-1}$ occurs in the first
$A_{k}$ and the $\left(  2^{k}-j+1\right)  ^{th}$ occurrence of
$z^{\pm1} $ occurs in the first $A_{k}^{-1}$ as the opposite power
as the $j^{th}$ occurrence. Hence if $p_{j}=z\omega_{j}$ then
$p_{2^{k}+1-j}=zA_{k}BA_{k}^{-1}\omega_{j}$ and if
$p_{j}=z\omega_{j}z^{-1}$ then
$p_{2^{k}+1-j}=zA_{k}BA_{k}^{-1}\omega _{j}z^{-1}$. Moreover, the
term $\gamma_{j}^{\pm1}\gamma_{2^{k}+1-j}^{\pm 1}$ in the formula
above will always be of the form
$\left(\gamma_{j}\gamma_{2^{k}+1-j}^{-1}\right)^{\pm 1}$.
Similarly, the $\left(2^k + j\right)^{th}$ occurrence of
$z^{\pm1}$ in $\alpha x^{-1}\alpha^{-1}$ occurs in the second
$A_{k}$ and the $\left( 2^{k+1}+1-j\right) ^{th}$ occurrence of
$z^{\pm1}$ occurs in the second $A_{k}^{-1}$ as the opposite power
as the $\left(2^k+j\right)^{th}$ occurrence.   Thus, if
$p_{2^{k}+j}=z\alpha x^{-1}B\omega_{j}$ then
$p_{2^{k+1}+1-j}=z\alpha x^{-1}\alpha^{-1}\omega_{j}$, if
$p_{2^{k} +j}=z\alpha x^{-1}B\omega_{j}z^{-1}$ then
$p_{2^{k+1}+1-j}=z\alpha x^{-1}\alpha^{-1}\omega_{j}z^{-1}$, and
$\gamma_{2^{k}+j}^{\pm1}\gamma
_{2^{k+1}+1-j}^{\pm1}=\left(\gamma_{2^{k}+j}\gamma
_{2^{k+1}+1-j}^{-1}\right)^{\pm1}$. In either case we see that for
$1\leq j\leq2^{k-1}$
\begin{align*}
\gamma_{j}\gamma_{2^{k}+1-j}^{-1}
&=\gamma^{p_{j}}\left(\gamma^{-1}\right)^{p_{2^{k}+1-j}}\\
&=\left[  \gamma,p_{j}p_{2^{k}+1-j}^{-1}\right]  ^{p_{j}}
\end{align*}
and
\begin{align*}
\gamma_{2^{k}+j}\gamma_{2^{k+1}+1-j}^{-1}
&=\gamma^{p_{2^{k}+j}}\left(\gamma^{-1}\right)  ^{p_{2^{k+1}+1-j}}\\
&=\left[  \gamma,p_{2^{k}+j}p_{2^{k+1}+1-j}^{-1}\right]
^{p_{2^{k}+j}}
\end{align*}
where
$l_{\ast}\left(  p_{j}p_{2^{k}+1-j}^{-1}\right)
=l_{\ast}\left(
\left(B^{-1}\right)^{\left(zA_{k}\right)^{-1}}\right) \in
K_{r}^{\left( n\right) }
$
and
$
l_{\ast}\left(  p_{2^{k}+j}p_{2^{k+1}+1-j}^{-1}\right)
=l_{\ast}\left(  \left(  B\alpha\right) ^{\left(z\alpha
x^{-1}\right)^{-1}}\right) \in K_{r}^{\left(  n\right)  }.
$
Therefore
\[
\gamma=\left[  \gamma,u_{1}\right]  ^{v_{1}}\cdots\left[
\gamma,u_{2^{k} }\right]  ^{v_{2^{k}}}\left[
\lambda_{1},\lambda_{2}\right]  ^{w_{1}} \cdots\left[
\lambda_{m-1},\lambda_{m}\right]  ^{w_{m/2}}
\]
where $l_{\ast}\left(  u_{j}\right)  \in K_{r}^{\left(  n\right)
}$ and $\lambda_{j}\in Ncl\left\langle \gamma\right\rangle $ as
desired.

Before proceeding, we note that $\gamma$ can be simplified to the
form
\begin{equation}
\gamma=\left[  \gamma,u_{1}\right]  ^{v_{1}}\cdots\left[
\gamma,u_{2^{k} +m}\right]  ^{v_{2^{k}+m}}\label{simpler}
\end{equation}
where $l_{\ast}\left(  u_{j}\right)  \in K_{r}^{\left(  n\right)
}$, since if $\lambda_{1},\lambda_{2}\in Ncl\left\langle \gamma
\right\rangle $ then $\left[ \lambda_{1},\lambda_{2}\right]  $ is
a product of elements of the form $\left[ \gamma,u\right]  ^{v}$
where $l_{\ast}\left(  u\right)  =1$. This is easily verified
using the relation
\[
\left[  ab,c\right]  =\left[  b,c\right]  ^{a}\left[  a,c\right] .
\]

We prove by induction that $l_{\ast}:P_{r}^{\left(  i\right)
}\rightarrow K_{r}^{\left(  i\right)  }$ is surjective for $0\leq
i\leq n+1$. It is clear that $l_{\ast}:P_{r}^{\left(  0\right)
}\rightarrow K_{r}^{\left(  0\right) }$ is surjective. Now assume
that $l_{\ast}:P_{r}^{\left(  i\right) }\twoheadrightarrow
K_{r}^{\left(  i\right)  }$ for $i\leq n$ and let $g\in
K_{r}^{\left(  i+1\right)  }$. We note that if
$G\twoheadrightarrow H$ is surjective then $\left[  G,G\right]
\twoheadrightarrow\left[  H,H\right]  $ is surjective. Therefore
it suffices to consider $g$ such that $g^{k} \in\left[
K_{r}^{\left(  i\right)  },K_{r}^{\left(  i\right) }\right]  $ for
some $k\neq0$. \ Since $P_{r}^{\left(  i\right)
}\twoheadrightarrow K_{r}^{\left(  i\right)  }$ is surjective,
$l_{\ast}:\left[  P_{r}^{\left( i\right)  },P_{r}^{\left( i\right)
}\right]  \twoheadrightarrow\left[ K_{r}^{\left( i\right)
},K_{r}^{\left(  i\right)  }\right]  $ is surjective and hence
there exists an $f\in\left[  P_{r}^{\left(  i\right)  }
,P_{r}^{\left(  i\right)  }\right]  $ such that $l_{\ast}\left(
f\right) =g^{k}$. Moreover, since $P\twoheadrightarrow K$ there
exists a $p\in P$ such that $l_{\ast}\left(  p\right)  =g$. It
follows that $l_{\ast}\left( f\right)  =l_{\ast}\left(
p^{k}\right)  $ and hence $p^{k}=f\lambda$ where $\lambda\in
Ncl\left\langle \gamma\right\rangle $. Since $i\leq n$,
$l_{\ast}\left(  u_{j}\right)  \in K_{r}^{\left(  i\right)  }$.
Hence, by the induction hypothesis, there exist $q_{j}\in
P_{r}^{\left(  i\right)  }$ with $u_{j}=q_{j}\Lambda_{j}$ for
$\Lambda_{j}\in Ncl\left\langle \gamma \right\rangle $. Using
(\ref{simpler}) we have
\[
\gamma=\left[  \gamma,q_{1}\Lambda_{1}\right]
^{v_{1}}\cdots\left[ \gamma,q_{2^{k}+m}\Lambda_{2^{k}+m}\right]
^{v_{2^{k}+m}}
\]
hence $\gamma\in\left[  P_{r}^{\left(  i\right)  },P_{r}^{\left(
i\right) }\right]  $. Since $\lambda\in Ncl\left\langle
\gamma\right\rangle $, $\lambda\in\left[  P_{r}^{\left(  i\right)
},P_{r}^{\left(  i\right) }\right]  $ so that
$p^{k}=f\lambda\in\left[  P_{r}^{\left(  i\right)  }
,P_{r}^{\left(  i\right)  }\right]  $. Thus $p\in P_{r}^{\left(
i+1\right) }$ and $l_{\ast}\left(  p\right)  =g$ which implies
that $l_{\ast} :P_{r}^{\left(  i+1\right)  }\rightarrow
K_{r}^{\left(  i+1\right)  }$ is surjective for $i\leq n$.   This
concludes the proof of (\ref{equ:samegroups}).

The isomorphisms in (\ref{equ:samegroups2}) and
(\ref{equ:samegroups}) imply the following three statements.
First, we can obtain the $G/G_{r}^{\left( n+1\right) }$ and
$H/H_{r}^{\left( n+1\right) }$-regular covers of $X$ and $X\left(
n,k\right)  $ respectively by restricting to the boundary of the
$E/E_{r}^{\left(  n+1\right) }$-regular cover of $W$. Secondly,
when $i=0$ the inclusion maps $i$ and $j $ induce isomorphisms on
$H_{1}\left( -,\mathbb{Z}\right) $ hence on $H^{1}\left(
-,\mathbb{Z}\right) $. In fact, $i$ and $j$ induce isomorphisms on
all integral homology groups. Thus, we can consider the
homomorphism $\psi$ as a homomorphism on $E$ and $H$ as well as
$G$. Lastly,
\[
\frac{G_{r}^{\left(  i\right)  }}{\left[  G_{r}^{\left(  i\right)
} ,G_{r}^{\left(  i\right)  }\right]
}\overset{\cong}{\longrightarrow} \frac{E_{r}^{\left(  i\right)
}}{\left[  E_{r}^{\left(  i\right)  } ,E_{r}^{\left(  i\right)
}\right]  }\overset{\cong}{\longleftarrow} \frac{H_{r}^{\left(
i\right)  }}{\left[  H_{r}^{\left(  i\right)  } ,H_{r}^{\left(
i\right)  }\right]  }
\]
for $0\leq i\leq n-1$.  In particular this implies that for any
$\psi\in H^{1}\left(  X\left(  n,k\right)  ,\mathbb{Z}\right)  $
and $0\leq i\leq n-1,$ $\delta_{i}^{X\left(  n,k\right)  }\left(
\psi\right)  =\delta_{i}^{W}\left( \psi\right)
=\delta_{i}^{X}\left(  \psi\right)  $.

We use the presentation given by the Fox Free Calculus (Section
\ref{sec:foxcalc}) to compute $\delta_{i}^{W}\left(  \psi\right)
$. Let $F=F\left\langle
x_{1},\ldots,x_{\mu},y_{1},\ldots,y_{l},w\right\rangle $, $\chi
:F\twoheadrightarrow G$.\ Recall that
\[
\frac{\partial\left[  C,D\right]  }{\partial w}=\left(  1-\left[
C,D\right] D\right)  \frac{\partial C}{\partial w}+\left(
C-\left[  C,D\right]  \right) \frac{\partial D}{\partial w}
\]
for any $C,D\in F$. We compute $\frac{\partial A_{k}}{\partial
w}=\left( 1-A_{k}x\right)  \frac{\partial A_{k-1}}{\partial w}$
and $\frac{\partial A_{1}}{\partial w}=1$ so
\[
\frac{\partial A_{k}}{\partial w}=\left(  1-A_{k}x\right)
\cdots\left( 1-A_{2}x\right)  .
\]
It follows that
\begin{align*}
\frac{\partial w\left[  x,\alpha\right]  }{\partial w}&=1+w\left(
x-\left[
x,\alpha\right]  \right)  \frac{\partial\alpha}{\partial w}\\
&=1+w\left(  x-\left[  x,\alpha\right]  \right)  \left(  \left(
1-\alpha B\right)  \frac{\partial A_{k}}{\partial w}+\left(
A_{k}-\alpha\right)
\frac{\partial B}{\partial w}\right) \\
&=1+w\left(  x-\left[  x,\alpha\right]  \right) \left( \left(
1-\alpha B\right)  \left(  1-A_{k}x\right) \cdots\left(
1-A_{2}x\right)  +\left(  A_{k}-\alpha\right) \frac{\partial
B}{\partial w}\right)  .
\end{align*}
Similarly we compute $\frac{\partial w\left[  x,\alpha\right]
}{\partial x}$ and $\frac{\partial w\left[  x,\alpha\right]
}{\partial v}$ when $v\in\left\{
x_2,\ldots,x_{\mu},y_{1},\ldots,y_{l}\right\}:$
\begin{align*} \frac{\partial w\left[
x,\alpha\right]  }{\partial x}&=w\left[  \left( 1-\left[
x,\alpha\right]  \right)  +\left(  x-\left[  x,\alpha\right]
\right)  \left(  \left(  1-\alpha B\right)  \frac{\partial
A_{k}}{\partial x}+\left(  A_{k}-\alpha\right)  \frac{\partial
B}{\partial x}\right)  \right]
\text{,}\\
\frac{\partial w\left[  x,\alpha\right]  }{\partial v}&=w\left(
x-\left[ x,\alpha\right]  \right)  \left(  \left(  1-\alpha
B\right)  \frac{\partial A_{k}}{\partial v}+\left(
A_{k}-\alpha\right)  \frac{\partial B}{\partial v}\right) .
\end{align*}
We note that $\frac{\partial A_{k}}{\partial v}=0$ since $A_{k}$
does not involve $v$. Moreover $\frac{\partial A_{k}}{\partial
x}^{\chi i_{\ast }\phi_{n}^{E}}=0$ since $\frac{\partial
A_{1}}{\partial x}=0$ and
\[
\frac{\partial A_{k}}{\partial x}=\left(  1-\left[
A_{k-1},x\right]  \right) \frac{\partial A_{k-1}}{\partial
x}+\left(  A_{k-1}-\left[  A_{k-1},x\right] \right)  .
\]
Using the involution and projecting to $\mathbb{Z}\left[
E/E_{r}^{\left( n+1\right) }\right]  $ we get
\begin{align*}
\overline{\frac{\partial w\left[  x,\alpha\right]  }{\partial
w}}^{\chi i_{\ast}\phi _{n}^{E}}&= \overline{1+\left(  x-1\right)
\left( 1-B\right) \left(  1-x\right)
^{k-1}} \\
 \overline{\frac{\partial w\left[  x,\alpha\right]  }{\partial
x_{i}}}^{\chi i_{\ast} \phi_{n}^{E}}&=\overline{\frac{\partial
w\left[ x,\alpha\right] }{\partial y_{i}}}^{\chi
i_{\ast}\phi_{n}^{E}}=0
\end{align*}
Thus
\begin{equation}\label{equ:matrix}
\left(\!
\begin{array}
{c@{\,}c} \left( \overline{ \frac{\partial R_{j}}{\partial
x_{i}}}\right) ^{^{\chi i_{\ast}
\phi_{n}^{E}}} & 0 \\
0 & \overline{ 1+\left(  x-1\right) \left( 1-B\right) \left(
1-x\right) ^{k-1}}
\end{array}\!
\right)  \otimes_{\mathbb{Z}\left[  E/E_{r}^{\left(  n+1\right)
}\right] }id_{\mathbb{K}_{n}^{E}\left[  t^{\pm1}\right]  }
\end{equation}
is a presentation of $H_{1}\left(  W,\ast;\mathbb{K}_{n}^{E}\left[
t^{\pm 1}\right]  \right)  $.

Let $\psi$ be a primitive class in $H^1\left(W;\mathbb{Z}\right)$
with $\psi\left(x\right)=t^p$ and let $\xi$ be a splitting of
$\overline{\psi}_n : E/E_r \twoheadrightarrow \mathbb{Z}$. We
rewrite $a=1+\left(  x-1\right) \left( 1-B\right) \left(
1-x\right) ^{k-1}$ as a polynomial in $t$.  The lowest and highest
degree terms of $a$ are $B$ and $$t^{kp}
\left(\xi\left(t\right)^{-kp}x^k \
x^{1-k}\left(1-B\right)x^{k-1}\right)$$ respectively. Our
assumption that $B\notin G_{r}^{\left( n+1\right) }$ (hence
$B\notin E_{r}^{\left( n+1\right)  }$ guarantees that $B-1$ is a
unit in $\mathbb{K}_{n}^{E}\left[ t^{\pm1}\right] $).  Therefore
the degree of $a$ is $k$. Moreover, $\operatorname{deg}
\overline{\sum{t^{i}a_i}} =\operatorname{deg} \sum{t^{i}a_i}$
hence
 \[\operatorname{deg} \overline{1+\left(  x-1\right) \left( 1-B\right)
\left( 1-x\right) ^{k-1}} = kp. \]  Lastly, $H_{1}\left(
X,\ast;\mathbb{K}_{n}^{G}\left[ t^{\pm1}\right] \right)  $ is
presented as $\left( \overline{ \frac{\partial R_{j}}{\partial
x_{i}}}\right)  ^{\chi i\phi_{n}^{G} }\otimes_{\mathbb{Z}\left[
G/G_{r}^{\left( n+1\right)  }\right] }id_{\mathbb{K}_{n}^{G}\left[
t^{\pm1}\right] }$ therefore
\[
\delta_{n}^{W}\left(  \psi\right)  =\delta_{n}^{X}\left(
\psi\right) +kp.
\]

To finish the proof we will show that $\delta_{n}^{X\left(
n,k\right) }\left( \psi\right)  \geq\delta_{n}^{W}\left(
\psi\right)  $ for any $\psi$. Since $\left( W,X\left( n,k\right)
\right) $ has only 2 and 3-handles $H_{1}\left( W,X\left(
n,k\right)  ; R \right) =0$. By Lemma \ref{homcob}, $H_{2}\left(
W,X\left( n,k\right)  ; R\right) $ is $ R$-torsion. We have the
following long exact sequence of pairs:
\[
\rightarrow TH_{2}\left(  W,X\left(  n,k\right)  ; R\right)
\overset{\partial}{\rightarrow}H_{1}\left(  X\left(  n,k\right) ;
R\right)  \overset{j_{\ast}}{\rightarrow}H_{1}\left( W; R\right)
\rightarrow0
\]
Since $j_{\ast}\left(  TH_{1}\left(  X\left(  n,k\right) ;
R\right) \right)  \subseteq TH_{1}\left( W; R\right)  $ we can
consider the homomorphism
\[
TH_{1}\left(  X\left(  n,k\right)  ; R\right)
\overset{j_{\ast}}{\rightarrow}TH_{1}\left(  W; R\right) .\] We
show that this map is surjective. Let $\sigma\in TH_{1}\left( W;
R\right)  $ and $\theta\in H_{1}\left( X\left( n,k\right) ;
R\right)  $ such that $j_{\ast}\left( \theta\right) =\sigma$.
There exists $r\in R$ such that $\sigma r=0$ so $j_{\ast}\left(
\theta r\right)  =j_{\ast}\left( \theta\right) r=\sigma r=0$. By
exactness, this implies that $\theta
r\in\operatorname{Im}\partial$.  Hence there exists $\upsilon\in
TH_{2}\left(  W,X\left(  n,k\right) ; R\right) $ with
$\partial\left(  \upsilon\right) =\theta r$ and $\upsilon s=0$ for
some non-zero $s\in R$. Therefore $\theta\left( rs\right) =\left(
\theta r\right) s=\partial\left( \upsilon\right) s=\partial\left(
\upsilon s\right) =\partial\left(  0\right)  =0$ which implies
$\theta\in TH_{1}\left(  X\left(  n,k\right) ; R\right) $ Lastly,
let $ R=\mathbb{K}_{n}\left[  t^{\pm1}\right] $ then $TH_{1}\left(
X\left(  n,k\right)  ;\mathbb{K}_{n}\left[ t^{\pm1}\right] \right)
$ and $TH_{1}\left( W;\mathbb{K}_{n}\left[ t^{\pm1}\right] \right)
$ can be considered as free $\mathbb{K}_{n}$-modules with finite
rank. We have $j_{\ast}$ surjective so
\[
\delta_{n}^{X\left(  n,k\right) }\left( \psi\right)
=\operatorname{rk}_{\mathbb{K}_{n}}TH_{1}\left( X\left( n,k\right)
;\mathbb{K}_{n}\left[  t^{\pm1}\right] \right) \geq
\operatorname{rk}_{\mathbb{K}_{n} }TH_{1}\left(
W;\mathbb{K}_{n}\left[ t^{\pm1}\right]  \right)  =\delta
_{n}^{W}\left(  \psi\right).
\]
\end{proof}

Before proceeding, we will construct a specific example.  We will
begin with zero surgery on the trivial link with 2 components and
subtly alter the manifold to increase $\delta_1$.

\begin{example}A specific example of $X\left(1,1\right)$ when $X=S^1 \times S^2 \# S^1 \times S^2$.
Let $X$ be 0-surgery on the 2-component trivial link. Then $X=S^1 \times S^2 \# S^1 \times S^2$ with $\pi_1$
generated by $x$ and $y$ (Figure \ref{startingX}).
\begin{figure}[htbp]
\begin{center}
\includegraphics{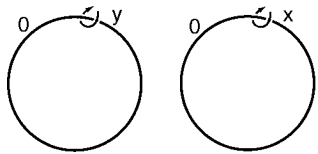}
\caption{$X=S^1\times S^2 \# S^1 \times S^2$} \label{startingX}
\end{center}
\end{figure}
Let $B=\left[x,y\right]$ and construct $X\left(1,1\right)$ as in the theorem by doing $k$-framed surgery on $c$ (see Figure \ref{anexample}).
\begin{figure}[htbp]
\begin{center}
\includegraphics
{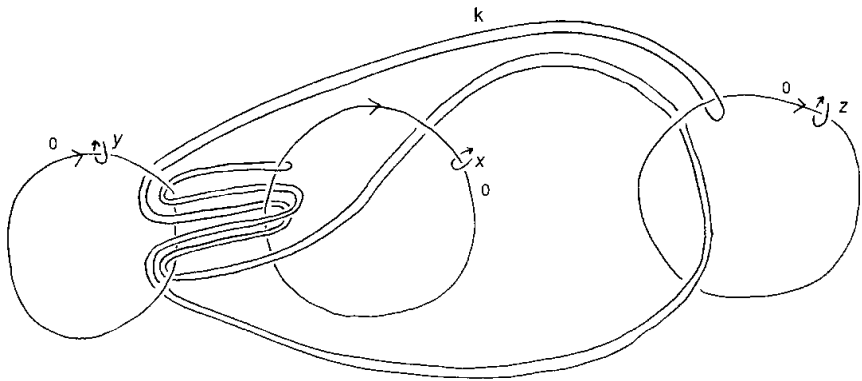} \caption{$X\left(1,1\right)$ is our resulting
manifold} \label{anexample}
\end{center}
\end{figure}
For all $\psi \in H^1 \left( X \left(1,1\right)\right)$,
\[
\delta_{0}^{X\left(  1,1\right)  }\left(  \psi\right)
=\delta_{0}^{X}\left( \psi\right)=0.
\]
Moreover, we have
\[
\delta_{1}^{X\left(  1,1\right)  }\left(  \psi_{x}\right)
\geq \delta_{1}^{X}\left(  \psi_{x}\right) + 1 =1.
\]
\end{example}

We can assume that the manifolds that we have constructed to be irreducible.

\begin{proposition}\label{irreducible}If $X$ is irreducible and
$\partial$-irreducible then $X\left(n,k\right)$ is irreducible and
$\partial$-irreducible.
\end{proposition}

\begin{proof}
Recall that $X\left(n,k\right)$ can be constructed from $X$ by
first taking a connected sum with $S^1 \times S^2$ and then doing
integer surgery on a curve $c$.  Recall that the first homology
class of $c$ was equal to $x^{-1}z$ where $x$ was a generator of
$H_1\left(X\right)$ and $z$ was the generator of $S^1 \times S^2$.
Let $M=\left( X \# S^1 \times S^2 \right)$.  We will show that
$M-c$ is irreducible and $\partial$-irreducible.  A theorem of M.
Scharlemann \cite[p. 481]{S} implies that $X\left(n,k\right)$ is
irreducible.  It is clear that $X\left(n,k\right)$ is
$\partial$-irreducible.

First we show that $M-c$ is $\partial$-irreducible.  Since $M$ is
$\partial$-irreducible, it suffices to show that $\ker i_\ast :
\pi_1\left(\partial c\right) \rightarrow \pi_1\left(M\right)$ is
trivial.  Any curve on the boundary of $c$ that is parallel to
$m\neq 0$ copies of $c$ is nontrivial in $\pi_1\left(M\right)$
since it is nonzero in homology.  Any other curve on $\partial c$
is homotopic to the meridian of $c$ which we showed to be
nontrivial in the proof of Theorem \ref{realization}.

Let $S$ be a nonseparating 2-sphere in $M$ that represents the
class $\{pt\}\times S^2$ and $N=M-S$.  Choose $S$ so that it
minimizes $\#\left(S \cap c\right)$. Note that $N=M-\left(B_1
\sqcup B_2\right)$ where $B_1$ and $B_2$ are disjoint $3$-balls in
$M$. After isotoping $c$ to make it transverse to $S$, let $P$ be
the punctured $2$-sphere in $M-c$ obtained by puncturing $S$ at
each intersection point with $c$. Let $M'$ be the manifold
obtained by cutting $M-c$ along $P$. $M'$ has two copies of $P$ in
it's boundary, denote these punctured $2$-spheres $P_1$ and $P_2$
(so that $P_i \subset \partial B_i$). We will show that $M'$ is
irreducible and $P_1$ is incompressible in $M'$.  It will follow
that $M-c$ is irreducible.

Suppose that $M-c$ is reducible and let $\Sigma$ be a 2-sphere in
$M-c$ that does not bound a 3-ball and minimizes $\# \left( \Sigma
\cap P\right)$.  Since $M'$ is irreducible, we have $\# \left(
\Sigma \cap P\right) \geq 1$.  Consider the intersection of
$\Sigma$ and $P$ and let $\alpha$ be an innermost circle on
$\Sigma$.  Then $\alpha$ bounds a disc $D$ in $M'$.  Since $P_1$
is incompressible in $M'$, $\alpha$ bounds a disc $E$ in $P$.  $D
\cup E$ is an embedded 2-sphere in $M'$ so it bounds a 3-ball $B$
in $M'$.  We use $B$ to isotope $\Sigma$ in $M-c$ to get rid of
the intersection $\alpha$.  This contradicts the minimality of $\#
\left( \Sigma \cap P\right)$.  Thus $M-c$ is irreducible.

We note that $M'$ is homeomorphic to $M-f\left(W\right)$ where $W$
is a wedge of spheres and $f:W \rightarrow M$ is an embedding of
$W$ into $M$.  Suppose that $M'$ is reducible and let $\Sigma$ be
an embedded 2-sphere in $M'$ that does not bound a 3-ball.  Since
$M$ is irreducible, $\Sigma$ bounds a 3-ball $B$ in $M$.  Hence
$M=B \cup_\Sigma V$.  $f\left(W\right)$ is connected and
$f\left(W\right) \cap \Sigma = \varnothing$ so either
$f\left(W\right) \subset B$ or $f\left(W\right) \subset V$.
However, the homology class of $c$ is equal to $x^{-1}z$ hence
$f\left(W\right) \nsubseteq B$.  Therefore $f\left(W\right)
\subset V$ hence $\Sigma$ must bound a ball in $M'$, a
contradiction. Thus $M'$ is irreducible.

Suppose that $P_1$ is compressible in $M'$.  Let $\alpha$ be an
curve on $P_1$ that bounds an embedded disc $D$ in $M'$. $\alpha$
bounds the discs $E_1$ and $E_2$ on $S$.  Since $M$ is
irreducible, $D \cup E_1$ bound a 3-ball $B$ in $M$.  If $B \cap
B_2 = \varnothing$ then either $D \cup E_1$ or $D \cup E_2$ bounds
a 3-ball $B'$ in $N$.  We can use $B'$ to isotope $c$ and reduce
the number of intersections of $c$ with $E_1$ or $E_2$ (see Figure
\ref{boundsball}) .
\begin{figure}[htbp]
\begin{center}
\includegraphics{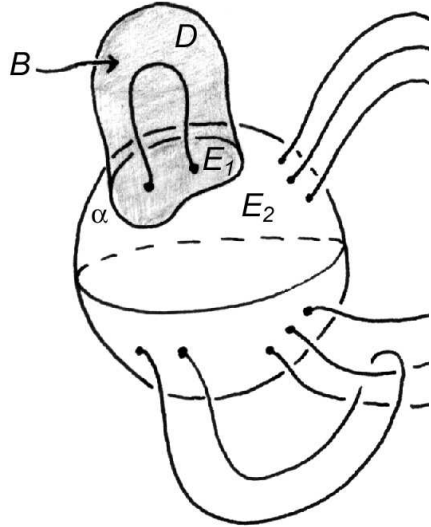}
\caption{$D$ is a compression disc} \label{boundsball}
\end{center}
\end{figure}
This contradicts the minimality of $\#\left(S \cap c\right)$.

Now suppose that $B \cap B_2 \neq \varnothing$.  Using $B$ we can
assume that either $D \cup E_1$ or $D \cup E_2$ bounds a 3-ball
$B'$ in $M-B_1$.  Let $S'=\partial B'$.  We note that $\#\left(S' \cap
c \right) < \#\left(S \cap c\right)$ since $c$ intersects $E_1$ and $E_2$.
Moreover, $S'$ and $\partial B_2$ cobound an embedded $S^2 \times
I$ in $N$.
\begin{figure}[htbp]
\begin{center}
\includegraphics
{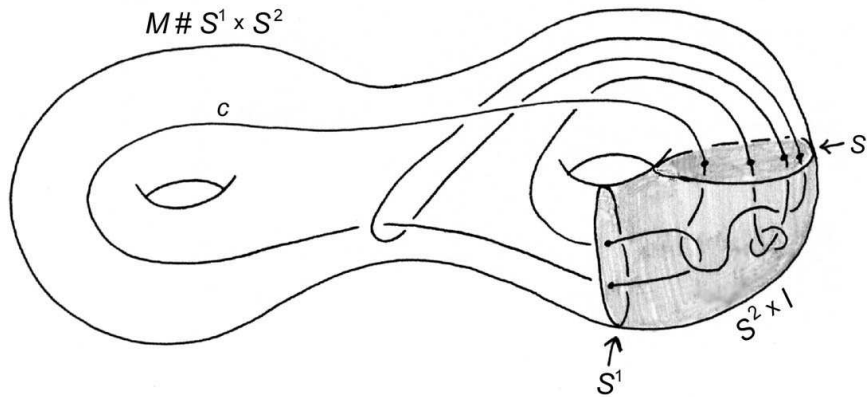}
\caption{$S$ and $S'$ cobound an $S^2 \times I$} \label{cobound}
\end{center}
\end{figure}
Therefore $S'$ and $S$ cobound an embedded $S^2 \times
I$ in $M$ (see Figure \ref{cobound}) which can be used to isotope
the curve $c$ to reduce the number of intersections with $S$. This
contradicts our minimality condition hence $P_1$ is incompressible
in $M'$.  This completes the proof that $M-c$ is irreducible.
\end{proof}

\begin{lemma}
\label{homcob}The manifold $X\left(  n,k\right)  $ is
$\mathcal{K}_{i}$-homology cobordant to $X$ for $i\leq n.$ That
is, $i:X\rightarrow W$ and $j:X\left(  n,k\right)  \rightarrow W$
induce isomorphisms on homology with $\mathcal{K}_{i}$
coefficients where $W$ is the cobordism between $X$ and $X\left(
n,k\right)$.
\end{lemma}

\begin{proof}
Consider the relative chain complex $C_{\ast}\left(  W,X\left(
n,k\right) \right)  $. Since $W$ is obtained from $X\left(
n,k\right)  $ by adding only 2 and 3-handles, $W$ is homotopy
equivalent to a cell complex obtained by adding a single 2 and
3-cell. Hence we can assume $C_{2}\left(  W,X\left( n,k\right)
\right)  =C_{3}\left(  W,X\left(  n,k\right)  \right)
\cong\mathbb{Z} $ and $C_{j}\left(  W,X\left(  n,k\right)  \right)
=0$ for all other $j$. For all $i\leq n$ we lift the cells of
$\left(  W,X\left( n,k\right)  \right)  $ to form the chain
complex of $\left(  \widetilde {W},\widetilde{X}\left(  n,k\right)
\right)  $
\[
0\rightarrow C_{3}\left(  \widetilde{W},\widetilde{X}\left(  n,k\right)
\right)  \otimes\mathcal{K}_{i}\overset{\overset{\sim}{\partial}_{3}\otimes
id}{\rightarrow}C_{2}\left(  \widetilde{W},\widetilde{X}\left(  n,k\right)
\right)  \otimes\mathcal{K}_{i}\rightarrow0
\]
where $\left(  \widetilde{W},\widetilde{X}\left(  n,k\right)
\right)  $ are the regular $E/E_{r}^{\left(  i+1\right)  }$-cover
of $\left(  W,X\left( n,k\right)  \right)  $. Since
$\overset{\sim}{\partial}_{3}\otimes
id:\mathcal{K}_{i}\rightarrow\mathcal{K}_{i}$,
$\overset{\sim}{\partial} _{3}\otimes id$ is an isomorphism if and
only if $\partial_{3}\otimes id\neq0$ if and only if
$\partial_{3}\left( \sigma\right)  \neq0$ for some $\sigma$. \ But
since $H_{\ast}\left(  W,X\left(  n,k\right)  \right)  =0$,
$\partial_{3}:$ $C_{3}\left(  W,X\left(  n,k\right)  \right)
\rightarrow C_{2}\left(  W,X\left(  n,k\right)  \right)  $ is not
the zero map, hence $\overset{\sim}{\partial_{3}}$ is not the zero
map. $\ $Therefore $H_{\ast }\left(  W,X\left(  n,k\right)
;\mathcal{K}_{i}\right)  =0$ which gives us
$j_{\ast}:H_{\ast}\left(  X\left(  n,k\right)
;\mathcal{K}_{i}\right) \overset{\cong}{\rightarrow}H_{\ast}\left(
W;\mathcal{K}_{i}\right)  $. The proof that
$i_{\ast}:H_{\ast}\left(  X;\mathcal{K}_{i}\right)  \overset{\cong
}{\rightarrow}H_{\ast}\left(  W;\mathcal{K}_{i}\right)  $ follows
almost verbatim except that $\left(  W,X\right)  $ has only cells
in dimension 1 and 2. \end{proof}

\begin{lemma}
\label{checknotzero}For each $n\geq1$, if $\delta_{n-1}\left(
\psi\right)  \neq0$ for some $\psi\in H^{1}\left(
X;\mathbb{Z}\right)  $ then $G_{r}^{\left( n\right)
}/G_{r}^{\left(  n+1\right)  }\neq0$.
\end{lemma}

\begin{proof}
If $\delta_{n-1}\left(  \psi\right)  \neq0$ for some $\psi$ then
the rank of $H_{1}\left(  X_{n-1}\right)  $ as an abelian group is
at least $1$, hence $H_{1}\left(  X_{n-1}\right)  /\left\{
\mathbb{Z}\!-\!torsion\right\}  \neq0$.  But, Lemma
\ref{derivedquotient} gives $G_{r}^{\left(  n\right)
}/G_{r}^{\left( n+1\right) }=
H_{1}\left( X_{n-1}\right) /\left\{
\mathbb{Z}\!-\!\operatorname{torsion}\right\}  $ so $G_{r}^{\left(
n\right) }/G_{r}^{\left(  n+1\right)  }\neq0$.
\end{proof}

We now proceed with the proof of Theorem \ref{increase}.

\begin{proof}[Proof of Theorem \ref{increase}]
Let $X_{0}$ be a 3-manifold with $r_{0}\left(  X_{0}\right)  =0$,
$G=\pi _{1}\left(  X_{0}\right)  $, $\mu=\beta_{1}\left(
X_{0}\right)  $ and whose universal abelian cover has non-trivial
$\beta_{1}$. For example, we could let $X_{0}$ be any 3-manifold
with $\beta_{1}\left(  X_{0}\right)  =\mu$ which fibers over
$S^{1}$ with fiber a surface of genus $g\geq2$ (see Proposition
\ref{fibered} ) .  Let $\left\{ x_{1},\ldots,x_{\mu}\right\} $ be
a basis of $H_{1}\left( X_{0} ;\mathbb{Z}\right)  /\left\{
\mathbb{Z}\!-\!\operatorname{torsion}\right\}  $ and $\left\{
\psi_{x_{1}},\ldots,\psi_{x_{\mu}}\right\}  $ the (Hom) dual basis
of $H^{1}\left(  X_{0};\mathbb{Z}\right)  $. Since
$\beta_{1}\left(  \left(  X_{0}\right)  _{\Gamma_{0}}\right)
>0$, $G_{r}^{\left(  1\right)  }/G_{r}^{\left(  2\right)  }\neq0$
we can use Theorem \ref{realization} with $k=1$ to construct a new
manifold $X_{1}$ with $\delta_{0}^{X_{0}}\left(  \psi\right)
=\delta_{0}^{X_{1}}\left( \psi\right)  $ and
$\delta_{1}^{X_{0}}\left(  \psi\right)  <\delta_{1}^{X_{1} }\left(
\psi\right)  $ for all $\psi\in H^{1}\left(
X_{1};\mathbb{Z}\right) $. We do this by first constructing
$X_{1}^{1}$ from $X_{0}$ to accomplish $\delta_{0}^{X_{0}}\left(
\psi\right)  =\delta_{0}^{X_{1}^{1}}\left( \psi\right)  $ for all
$\psi\in H^{1}\left(  X_{1}^{1};\mathbb{Z}\right)  $ and
$\delta_{1}^{X_{0}}\left(  \psi_{x_{1}}\right)
<\delta_{1}^{X_{0}}\left( \psi_{x_{1}}\right)
+1\leq\delta_{1}^{X_{1}^{1}}\left(  \psi_{x_{1}}\right) $ \ We
note that Theorem \ref{realization} guarantees
$\delta_{1}^{X_{1}^{1} }\left( \psi\right)
\geq\delta_{1}^{X_{0}}\left(  \psi\right)  $ for all other $\psi$.
Now we continue this for all other basis elements $\psi
_{x_{2}},\ldots,\psi_{x_{\mu}}$ to get a 3-manifold
$X_{1}=X_{1}^{\mu}$ with $\delta_{0}^{X_{0}}\left(  \psi\right)
=\delta_{0}^{X_{1}}\left( \psi\right)  $ and
$\delta_{0}^{X_{0}}\left(  \psi\right)  <\delta_{0}^{X_{1} }\left(
\psi\right)  $ for all $\psi\in H^{1}\left(
X_{1};\mathbb{Z}\right) $.

In particular, $\delta_{1}^{X_{1}}\left(  \psi\right)  >0$ so by
Lemma \ref{checknotzero} $G_{r}^{\left(  2\right)  }/G_{r}^{\left(
3\right)  } \neq0$. Hence we can construct $X_{2}$ with
$\delta_{i}^{X_{2}}\left( \psi\right)  =\delta_{i}^{X_{1}}\left(
\psi\right)  $ when $i\leq1$ and $\delta_{2}^{X_{1}}\left(
\psi\right)  <\delta_{2}^{X_{2}}\left( \psi\right)  $ for all
$\psi\in H^{1}\left(  X_{2};\mathbb{Z}\right)  $. We continue this
process until we obtain a 3-manifold $X=X_{m}$ with $\delta
_{0}^{X}\left(  \psi\right)  <\delta_{1}^{X}\left(  \psi\right)
<\cdots<\delta_{m}^{X}\left(  \psi\right)  $ for all $\psi\in
H^{1}\left( X;\mathbb{Z}\right)  $. Since $r_{0}\left(
X_{0}\right)  =0$, Lemma \ref{homcob} guarantees that $r_{0}\left(
X\right)  =0$ hence $\left\| \psi\right\|  _{A}=\delta_{0}\left(
\psi\right)  $ and $\delta_{m}\left( \psi\right)  \leq\left\|
\psi\right\|  _{T}$.

If we choose $X_0$ to be closed, then $X$ will be closed. Finally,
to guarantee that $X$ is irreducible, it suffices to choose $X_0$
irreducible by Proposition \ref{irreducible}.
\end{proof}

We note that since $r_{i}\left(  X\right)  =0$, $\delta_{i}\left(
\psi\right)  =\overline{\delta}_{i}\left(  \psi\right)  $ hence we
could have stated this theorem in terms of $\overline{\delta}_{i}$
as well as $\delta _{i}$.

\section{Applications}\label{applications}
We show that the higher-order degrees give new computable
algebraic obstructions to 3-manifolds fibering over $S^1$ even
when the classical Alexander module fails. Moreover, using the
work of P. Kronheimer, T. Mrowka, and S. Vidussi we are able to
show that the higher-order degrees give new computable algebraic
obstructions a 4-manifold of the form $X \times S^1$ admitting a
symplectic structure, even when the Seiberg-Witten invariants
fail.

\subsection{Fibered 3-manifolds}

Recall that if $X$ is a compact, orientable 3-manifold that fibers
over $S^{1}$ then by Proposition \ref{fibered}, the higher-order
ranks must be zero.  Moreover, if $\beta_1\left(X\right) \geq 2$
and $\psi$ is dual to a fibered surface then
$\delta_n\left(\psi\right)$ must be equal to the Thurston norm for
all $n$ and hence are constant as a function of $n$.  We define
the following function of $\psi$.  Let $d_{ij} :
H^1\left(X;\mathbb{Z}\right) \rightarrow \mathbb{Z}$ be defined by
$d_{ij} = \delta_i - \delta_j$ for $i,j \geq 0$. Note that $d_{ij}
= 0$ if and only if $\delta_i = \delta_j$ for all $\psi \in
H^1\left(X;\mathbb{Z}\right)$.

\begin{theorem}
\label{obsfib}Let $X$ be a compact, orientable 3-manifold.
 If at least one of the following conditions is satisfied then $X$ does not
 fiber over $S^1$.  \newline
$(1)$ $r_n\left(X\right)\neq 0$ for some $n \geq 0$,
\newline
$(2)$ $\beta_{1}\left( X\right) \geq2$ and there exists $i,j\geq
0$ such that $d_{ij}\left(\psi\right) \neq 0$ for all $\psi\in
H^{1}\left( X;\mathbb{Z}\right)$,
\newline
$(3)$ $\beta_{1}\left( X\right) =1$ and $d_{ij}\left(\psi\right)
\neq 0$ for some $i,j\geq 1$ and $\psi \in
H^1\left(X;\mathbb{Z}\right)$,
\newline$(4)$ $\beta_{1}\left( X\right) =1$, $X\ncong S^1 \times S^2$, $X \ncong S^1 \times D^2$ and $d_{0j}\left( \psi\right)\neq 1+\beta_{3}\left( X\right)$ for some
$j\geq 1$ where $\psi$ is a generator of
$H^1\left(X;\mathbb{Z}\right)$.
\end{theorem}
\begin{proof}We consider each of the cases separately.\newline
(1) This follow immediately from
Proposition~\ref{fibered}.\newline To prove that each of last
three conditions implies that $X$ does not fiber over $S^1$ we can
assume that $r_n\left(X\right)=0$ for all $n \geq 0$. Otherwise
the conclusion would be (vacuously) true since $X$ would satisfy
condition (1). Hence $\delta_n=\overline{\delta}_n$ for all $n
\geq 0$ by Remark~\ref{rem:del_eq_delbar}.\newline (2) If $X$
fibers over $S^1$ and $\beta_1\left(X\right) \geq 2$ then for all
$n \geq 0$, $\delta_n\left(\psi\right)=\left\| \psi \right\|_T$
for some $\psi \in H^1\left(X;\mathbb{Z}\right)$ by
Proposition~\ref{fibered}. Hence $\delta_n\left(\psi\right)$ is a
constant function of $n$. In particular,
$d_{ij}\left(\psi\right)=0$ for all $i,j\geq 0$ which contradicts
our hypothesis. \newline(3) If $X$ fibers over $S^1$ and
$\beta_1\left(X\right)=1$ then for all $n \geq 1$ and $\psi \in
H^1\left(X;\mathbb{Z}\right)$, $\delta_n\left(\psi\right)=\left\|
\psi \right\|_T$ by Proposition~\ref{fibered}. Hence
$\delta_n\left(\psi\right)$ is a constant function of $n$ for $n
\geq 1$. In particular, $d_{ij}\left(\psi\right)=0$ for all
$i,j\geq 1$ and $\psi \in H^1\left(X;\mathbb{Z}\right)$ which
contradicts our hypothesis.
\newline (4) If $X$ fibers over $S^1$,
$\beta_1\left(X\right)=1$, $X\ncong S^1 \times S^2$, $X\ncong S^1
\times D^2$ and $\psi$ is a generator of
$H^1\left(X;\mathbb{Z}\right)$ then by Propostion~\ref{fibered},
$\delta_0\left(\psi\right)=\left\|\psi\right\|_T+
1+\beta_3\left(X\right)$. The rest of the proof is similar to the
previous two cases.
\end{proof}

The previously known algebraic obstructions to a 3-manifold
fibering over $S^{1}$ are that the Alexander module
$H_1\left(X;\mathbb{Z}\Gamma_0\right)$ is finitely generated and
(when $\beta_{1}\left( X\right)  =1$) the Alexander polynomial is
monic. If $\beta_{1}\left( X\right) \geq2$, the Alexander module
being finitely generated implies that $r_{0}\left( X\right) =0$.

Consider the 3-manifolds in Theorem \ref{increase}.  We note that
$\delta_1 > \delta_0$ hence they cannot fiber over $S^1$.
Moreover, they can be chosen to have the same Alexander module as
those of a 3-manifold that fibers over $S^1$ as remarked in the
first paragraph in the proof of Theorem \ref{increase}.

\begin{corollary}\label{cor:notfiber}For each $\mu\geq 1$, Theorem
\ref{increase} gives an infinite family of
closed irreducible 3-manifolds $X$ where
$\beta_1\left(X\right)=\mu$, $X$ does not fiber over $S^1$, and
$X$ cannot be distinguished from a fibered 3-manifold using the
classical Alexander module.
\end{corollary}

\subsection{Symplectic 4-manifolds of the form $X \times S^1$}

We now turn our attention to symplectic 4-manifolds of the form
$X\times S^{1}$. It is well known that if $X$ is a closed
3-manifold that fibers over $S^1$ then $X \times S^1$ admits a
symplectic structure.  Taubes conjectures the converse to be true.

\begin{conjecture}
[Taubes]\label{taubes}Let $X$ be a 3-manifold such that $X\times
S^{1}$ admits a symplectic structure. Then $X$ admits a fibration
over $S^{1}$.
\end{conjecture}

Using the work of Meng-Taubes and Kronheimer-Mrowka, Vidussi
\cite{Vi} has recently given a proof of McMullen's inequality
using Seiberg-Witten theory. This generalizes the work of
Kronheimer \cite{K2} who dealt with the case that $X$ is the
0-surgery on a knot. Moreover, Vidussi shows that if $X\times
S^{1}$ admits a symplectic structure (and $\beta_1\left(X\right)
\geq 2$) then the Alexander and Thurston norms of $X$ coincide on
a cone over a face of the Thurston norm ball of $X$, supporting
the conjecture of Taubes.

\begin{theorem}
[\cite{K2,V,Vi}]\label{vid} Let $X$ be an closed, irreducible
3-manifold such that $X \times S^1$ admits a symplectic structure.
If $\beta_1\left(X\right) \geq 2$ there exists a $\psi \in
H^1\left(X;\mathbb{Z}\right)$ such that $\left\|\psi\right\|_A =
\left\|\psi\right\|_T$. If $\beta _{1}\left( X\right) =1$ then for
any generator $\psi$ of $H^1\left(X;\mathbb{Z}\right)$,
$\left\|\psi\right\|_A = \left\|\psi\right\|_T + 2$.
\end{theorem}

Consequently, we show that the higher-order degrees of a
3-manifold $X$ give new computable algebraic obstructions to a
4-manifold of the form $X\times S^{1}$ admitting a symplectic
structure.

\begin{theorem}
\label{sympl}Let $X$ be a closed irreducible 3-manifold.  If at
least one of the following conditions is satisfied then $X \times
S^1$ does not admit a symplectic structure.
\newline$(1)$ $\beta_1\left(X\right) \geq 2$ and there exists an $n \geq 1$ such that
$\overline{\delta}_{n}\left( \psi\right)
>\overline{\delta}_{0}\left(  \psi\right) $ for all $\psi\in H^{1}\left(  X;\mathbb{Z}\right)$.
\newline$(2)$ $\beta_1\left(X\right) =1$,  $\psi$ is a generator of $H^{1}\left(  X;\mathbb{Z}\right)$,  and
$\overline{\delta}_{n}\left( \psi\right)
>\overline{\delta}_{0}\left(  \psi\right) - 2$ for some $n \geq 1$.
\end{theorem}
\begin{proof}If $\beta_2\left(X\right) \geq 2$, $n \geq 1$, and $X \times
S^1$ admits a symplectic structure then by Theorem
\ref{delbarthm}, Theorem \ref{vid}, and Proposition
\ref{pro:del0}, $\overline{\delta}_n\left(\psi\right) \leq
\left\|\psi\right\|_T = \overline{\delta}_0\left(\psi\right)$ for
some $\psi \in H^1\left(X;\mathbb{Z}\right)$.  If
$\beta_2\left(X\right)=1$, $n \geq 1$,  $\psi$ is a generator of
$H^1\left(X;\mathbb{Z}\right)$ and $X \times S^1$ admits a
symplectic structure then by Theorems \ref{delbarthm} and
\ref{vid}, $\overline{\delta}_n\left(\psi\right) \leq
\left\|\psi\right\|_T = \overline{\delta}_0\left(\psi\right)-2$.
\end{proof}

Thus, Theorem \ref{increase} gives examples of 4-manifolds of the
form $X\times S^{1}$ which do not admit a symplectic structure but
cannot be distinguished from a symplectic 4-manifold using the
invariants of Seiberg-Witten theory.

\begin{corollary}\label{cor:notsympl}For each $\mu\geq 1$, Theorem \ref{increase} gives an infinite family of
4-manifolds $X \times S^1$ where $\beta_1\left(X\right)=\mu$, $X
\times S^1$ does not admit a symplectic structure, and $X$ cannot
be distinguished from fibered 3-manifold using the classical
Alexander module.
\end{corollary}

We note that the conditions in Theorem \ref{sympl} are (strictly)
stronger that the conditions in Theorem \ref{obsfib}. The cause of
this discrepancy is our lack of knowledge of the behavior of
higher-order degrees when $X \times S^1$ admits a symplectic
structure.  We make the following conjecture, supporting the
conjecture of Taubes.

\begin{conjecture}\label{con:allequal}
If $X$ is a closed, orientable, irreducible 3-manifold such that
$X\times S^{1}$ admits a symplectic structure then there exists a
$\psi \in H^1\left(X;\mathbb{Z}\right)$ such that
$\overline{\delta}_{n}\left( \psi \right) =\left\| \psi\right\|_T$
for all $n \geq 1$.
 \end{conjecture}

 More interesting would be the possibility of finding an a symplectic
4-manifold of the form $X\times S^{1}$ such that
$\overline{\delta}_{1}\left( \psi\right)
<\overline{\delta}_{0}\left( \psi\right)$ for all $\psi\in
H^{1}\left(  X;\mathbb{Z}\right)$; giving a counterexample to the
conjecture of Taubes \ref{taubes}.

We conclude with the remark that Conjecture \ref{con:allequal} is
true when $X$ is a knot complement in $S^3$, \cite[Theorem
9.5]{Co}. The proof of this relies on the fact that the
higher-order degrees are non-decreasing in $n$ \cite[Theorem
5.4]{Co} and are bounded by the Thurston norm. More precisely,
Cochran proves that $\delta_0\left(\psi\right) -1 \leq
\delta_1\left(\psi\right) \leq \cdots \leq
\delta_n\left(\psi\right) \cdots$ whenever $X = S^3 -K$ and $\psi$
is a generator $H^1\left(X;\mathbb{Z}\right)$. Moreover, the proof
of Theorem 5.4 in \cite{Co} can be modified to prove that
higher-order degrees are non-decreasing (in $n$) when $X$ is any
finite CW-complex homotopy equivalent to a 2-complex with Euler
characteristic zero. Hence Conjecture \ref{con:allequal} is also
true for any 3-manifold with \textbf{non-empty} toroidal boundary.

\vspace{10 pt} \noindent Shelly L. Harvey
\\Department of Mathematics,
\\Rice University, Houston, TX 77005
\\ email: \verb"shelly@math.rice.edu"

\end{document}